\documentclass[12pt,leqno]{amsbook}
\usepackage{latexsym,amssymb,amscd,fancyhdr,eqname,smartref,graphicx,amsopn}
\usepackage[Sonny]{fncychap}

\def\NZQ{\Bbb}               
\def\NN{{\NZQ N}}

\def\ZZ{{\NZQ Z}}

\def\frk{\frak}               

\def\pp{{\frk p}}

\def\qq{{\frk q}}

\def\mm{{\frk m}}

\def\B'c{{\mathcal{B'}}}
\def\U'c{{\mathcal{U'}}}

\def\xb{{\bold x}}

\def\opn#1#2{\def#1{\operatorname{#2}}} 
\def\opn#1#2{\def#1{\operatorname{#2}}}
\opn\chara{char}
\opn\length{\ell}
\opn\projdim{proj\,dim}
\opn\injdim{inj\,dim}
\opn\ini{in}
\opn\rank{rank}
\opn\depth{depth}
\opn\height{ht}
\opn\bigheight{bight}
\opn\embdim{emb\,dim}
\opn\codim{codim}

\opn\Tr{Tr}
\opn\bigrank{big\,rank}
\opn\superheight{superheight}\opn\lcm{lcm}
\opn\trdeg{tr\,deg}%
\opn\reg{reg}
\opn\lreg{lreg}
\opn\set{set}
\opn\supp{Supp}
\opn\shad{Shad}
\opn\del{del}
\opn\succ{succ}
\opn\pred{pred}
\opn\div{div}
\opn\Div{Div}
\opn\cl{cl}
\opn\Cl{Cl}
\opn\Spec{Spec}
\opn\Supp{Supp}
\opn\supp{supp}
\opn\Sing{Sing}
\opn\Ass{Ass}
\opn\Mon{Mon}
\opn\Min{Min}
\opn\Ann{Ann}
\opn\Rad{Rad}
\opn\Soc{Soc}
\opn\depth{depth}
\opn\Ker{Ker}
\opn\Coker{Coker}
\opn\Im{Im}
\opn\Hom{Hom}
\opn\Tor{Tor}
\opn\Ext{Ext}
\opn\End{End}
\opn\Aut{Aut}
\opn\id{id}

\opn\nat{nat}
\opn\GL{GL}
\opn\SL{SL}
\opn\mod{mod}
\opn\ord{ord}
\opn\id{id}
\opn\chr{char}
\opn\ara{ara}
\opn\indeg{indeg}
\opn\cd{cd}
\opn\link{link}
\opn\aff{aff}
\opn\con{conv}
\opn\relint{relint}
\opn\st{st}
\opn\lk{lk}
\opn\cn{cn}
\opn\core{core}
\opn\vol{vol}
\opn\gr{gr}

\def\pot#1#2{#1[\kern-0.28ex[#2]\kern-0.28ex]}

\opn\dirlim{\underrightarrow{\lim}}
\opn\invlim{\underleftarrow{\lim}}
\def\pnt{{\raise0.5mm\hbox{\large\bf.}}}

\def\twoline#1#2{\aoverb{\scriptstyle {#1}}{\scriptstyle {#2}}}
\newcommand{\aoverb}[2]{{\genfrac{}{}{0pt}{1}{#1}{#2}}}

\let\to=\rightarrow

\def\Implies{\ifmmode\Longrightarrow \else
     \unskip${}\Longrightarrow{}$\ignorespaces\fi}
\def\implies{\ifmmode\Rightarrow \else
     \unskip${}\Rightarrow{}$\ignorespaces\fi}
\def\iff{\ifmmode\Longleftrightarrow \else
     \unskip${}\Longleftrightarrow{}$\ignorespaces\fi}

\let\:=\colon

\newtheorem{Theorem1}{\bfseries Theorem}[section]
\newtheorem{Lemma1}[Theorem1]{\bfseries Lemma}
\newtheorem{Corollary1}[Theorem1]{\bfseries Corollary}
\newtheorem{Proposition1}[Theorem1]{\bfseries Proposition}
\newtheorem{Remark1}[Theorem1]{\bfseries Remark}

\newtheorem{Example1}[Theorem1]{\bfseries Example}
\newtheorem{Examples1}[Theorem1]{\bfseries Examples}
\newtheorem{Definition1}[Theorem1]{\bfseries Definition}

\newtheorem{Conjecture1}[Theorem1]{\bfseries Conjecture}
\newtheorem{Question1}[Theorem1]{\bfseries Question}

\newtheorem{Theorem111}{\bfseries Theorem}

\newtheorem{Proposition111}[Theorem111]{\bfseries Proposition}

\let\epsilon=\varepsilon
\let\phi=\varphi
\let\kappa=\varkappa
\numberwithin{section}{chapter}
\numberwithin{equation}{chapter}
\numberwithin{figure}{chapter}
\title{}
\author{}
\date{}

\begin{document}
\begin{titlepage}

\pagestyle{empty}
\vspace{0in}
\begin{center}
\large\textbf{     UNIVERSITATEA "OVIDIUS" CONSTAN\c TA} 
	
\noindent\large\textbf{FACULTATEA DE MATEMATIC\u A \c SI INFORMATIC\u A } 	

\large\textbf{\c SCOALA DOCTORAL\u A}

\vspace{1.5in}

\huge{\textbf{TEZ\u A DE DOCTORAT}}
\end{center}
\vspace{2in}
\begin{flushleft}
	\large{\textbf{$\ \ $CONDUC\u ATOR \c STIIN\c TIFIC}
	
		\textbf{PROF. UNIV. DR. MIRELA \c STEF\u ANESCU}	}
\end{flushleft}

\vspace{0.5in}
\begin{flushright}
	\[
\]
\large{\textbf{DOCTORAND$\ \ $\qquad\qquad}

	\textbf{OANA-\c STEFANIA OLTEANU}}
\end{flushright}

\vspace{0.5in}
\vspace{0.8in}
\begin{center}
\large{\textbf{CONSTAN\c TA\\ 2011}}
\end{center}
\end{titlepage}
\newpage
\large\textbf{    } 

\noindent\large\textbf{ } 	

\large\textbf{ }

\vspace{2in}
\begin{center}
\LARGE\textbf{{INVARIANTS OF SOME CLASSES OF MONOMIAL IDEALS}}
\end{center}
\vspace{2in}
\begin{flushleft}
	\large{\textbf{$\ \ $}
	
		\textbf{}	}
\end{flushleft}

\vspace{0.5in}
\begin{flushright}
\large{\textbf{}

	\textbf{}}
\end{flushright}
\normalsize

\pagestyle{empty}
\chapter*{Acknowledgments}
\thispagestyle{empty}
I am very grateful to my advisor, Professor Mirela \c Stef\u anescu, for her guidance provided during my Ph.D. study and research. 

I would like to thank Professor J\"urgen Herzog for all the mathematical discussions and advices offered.

I am also extremely grateful to Professor Naoki Terai, for all his support and advices offered during his visits to our university.

I wish to express my gratitude to Professor Viviana Ene, whose patience and guidance were invaluable to the completion of my Ph.D. study.

Finally, my heartfelt gratitude goes to my parents and to my two sisters, for the permanent encouragement and support, without whom none of this would have been possible. 

This thesis would not have been possible without their help.

\nopagebreak
\pagestyle{empty}
\rm{\bf\tableofcontents}
\frontmatter
\addcontentsline{toc}{chapter}{Preface} 
 
\chapter*{Preface}
\pagestyle{plain}
Commutative algebra has an important role in the present development of mathematics. It is also at the foundation of new trends in the modern research. By combining different techniques, several branches developed along years, such as algebraic geometry, algebraic combinatorics or algebraic statistics. 


Monomial ideals are at the intersection of commutative algebra and combinatorics. Passing, by using Gr\"obner bases theory, from a polynomial ideal to its initial ideal, many important properties of the original ideal are preserved.

The radical monomial ideals which are ideals generated by squarefree monomials, have a beautiful combinatorial interpretation in terms of simplicial complexes. Namely, with any simplicial complex one may associate a squarefree monomial ideal, generated by the squarefree monomials corresponding to the minimal non-faces of the simplicial complex. This ideal is called the Stanley--Reisner ideal of the simplicial complex. This process is self--dual, namely with a squarefree monomial ideal $I\subset S=k[x_1,\ldots,x_n]$, we may associate a simplicial complex $\Delta$ whose Stanley--Reisner ideal coincides with $I$. By this process, the combinatorial properties of the simplicial complexes can be described by using algebraic methods on the Stanley--Reisner ideals, and many invariants of the squarefree monomial ideals can be studied by combinatorial methods applied to the associated simplicial complexes. Another useful concept is the Alexander duality. A famous result concerning Alexander duality is the Eagon--Reiner Theorem \cite{ER}, which relates homological data of the Stanley--Reisner ideal of a simplicial complex to combinatorial properties of its Alexander dual. At the same time, R. Stanley gave a nonpure generalization of the notion of Cohen--Macaulay simplicial complex, by defining the sequentially Cohen--Macaulay simplicial complexes. Shortly after, J. Herzog and T. Hibi related the sequentially Cohen--Macaulayness property to the componentwise linearity concept via the Alexander duality.

In this thesis we are interested in describing some homological invariants of certain classes of monomial ideals. We will pay attention to the squarefree and non-squarefree lexsegment ideals. Their relevance is given by the major role played by them in the study of Hilbert function. It is known that initial lexsegment ideals possess the maximal graded Betti numbers among all the graded ideals with a given Hilbert function. 

The aim of this thesis is to present the original results obtained in this field. These results are contained in the following papers:

$\bullet$ \cite{OOS} A. Olteanu, O. Olteanu, L. Sorrenti, \textit{Gotzmann lexsegment ideals}, Le Matematiche, {\bf 63}(2008) Fasc. II, 229--241;

$\bullet$ \cite{EOT} V. Ene, O. Olteanu, N. Terai, \textit{Arithmetical rank of lexsegment edge ideals}, Bull. Math. Soc. Sci. Math. Roumanie, {\bf 53}(101) no.4  (2010), 315--327;

$\bullet$ \cite{OO} O. Olteanu, \textit{Classes of sequentially Cohen--Macaulay squarefree lexsegment ideals}, Algebra Colloquium, accepted;

$\bullet$ \cite{EKOT} V. Ene, K. Kimura, O. Olteanu, N. Terai, \textit{Invariants of squarefree lexsegment ideals of degree $3$}, work in progress.

The thesis is structured in four chapters, as follows. 

The \textbf{first chapter} represents a brief introduction of the notions and concepts which are intensively used along this thesis. We first recall the definitions of the Krull dimension and depth of a module. We also present Cohen--Macaulay modules and the sequentially Cohen--Macaulay concept.

The monomial ideals play a key role in the topic of this thesis, thus it is payed a special attention to recalling many important concepts and properties concerning them, like primary decomposition, linear quotients and resolutions.

One of the combinatorial approaches in the study of monomials ideals uses the Hilbert function. There are known two basic results in this sense: Macaulay's theorem and the Gotzmann persistence theorem. We recall both of them and we present some basic properties of Gotzmann ideals. Moreover, we prove some identities, which appears in \cite{OOS}, involving the operators related to the binomial expansion of an integer (Lemma \ref{b(d)=c(d)} and Lemma \ref{c^d=0 dnd c<egal b}, Section $1.2.5$).

Another important combinatorial technique used to obtain properties of monomial ideals is given by the polarization. By the process of polarization, we may pass from a monomial ideal to a squarefree monomial ideal which preserves many invariants of the original ideal. This shows the importance of the squarefree monomial ideals. An important tool used in studying squarefree monomial ideals is given by simplicial complexes, which are combinatorial objects. We mention concepts and notions concerning the simplicial complexes, which will be used in this thesis.  

Given a polynomial ideal $I\subset S=k[x_1,\ldots,x_n]$ it is interesting to find the minimal number of polynomials which generates $\sqrt I$. This number is called the arithmetical rank of $I$. This problem comes from algebraic geometry. The arithmetical rank of $I$ is the minimal number of equations which defines the algebraic variety $V(I)$. We briefly recall some important bounds for arithmetical rank in the last section of the first chapter.

Squarefree lexsegment ideals are the main subject of the \textbf{second chapter}. Let $S=k[x_1,\ldots,x_n]$ be the polynomial ring in $n$ variables over a field $k$. We order the monomials in $S$ lexicographically with $x_1>x_2>\ldots >x_n$. For an arbitrary integer $q\geq 2$, we denote by $\Mon^s_q(S)$ the set of all squarefree monomials of degree $q$ in $S$. A \textit{squarefree lexsegment set of degree $q$} determined by the monomials $u,\ v\in \Mon_q^s(S)$, $u\geq_{lex} v$, is a subset of $\Mon^s_q(S)$ of the form $L(u,v)=\{w\in\Mon_q^s(S):u\geq_{lex}w\geq_{lex} v\}$. An ideal generated by a squarefree lexsegment set is called \textit{squarefree lexsegment ideal}. In particular, one can define \textit{initial} and \textit{final squarefree lexsegment sets} to be sets of the form $L^i(v)=\{w\in\Mon_q^s(S):w\geq_{lex} v\}$, respectively $L^f(u)=\{w\in\Mon_q^s(S):u\geq_{lex}w\}$ and \textit{initial} and \textit{final squarefree lexsegment ideals}, accordingly. 

The concept of squarefree lexsegment ideal was introduced in \cite{AHH} by A. Aramova, J. Herzog and T. Hibi, where the notion was associated with the nowadays concept of initial squarefree lexsegment ideals, but it was also studied in \cite{ADH}, \cite{AHH}, \cite{B} and \cite{BoS}.  


In this chapter we are interested in studying the completely squarefree lexsegment ideals. In the first step, in Section $2.1$, we explicitly compute the minimal primary decomposition for initial and final squarefree lexsegment ideals, results which are contained in the paper \textit{Classes of sequentially Cohen--Macaulay squarefree lexsegment ideals}, \cite{OO}. The theorem shows that the primary decomposition of initial squarefree lexsegment ideals can be written easily, just by looking at the ends of the lexsegment.

\begin{Theorem111}\label{prim dec initial}\cite{OO}
Let $I\subset k[x_1,\ldots,x_n]$ be the initial squarefree lexsegment ideal generated in degree $q$, determined by the monomial $v=x_{j_1}\cdots x_{j_q}$, with $2\leq j_1<\ldots<j_q\leq n$. Consider the sets $A_t=[j_t]\setminus\{j_1,\ldots,j_{t-1}\}$, for $1\leq t\leq q$. Then $I$ has the minimal primary decomposition of the form:
$$I=\left(\bigcap\limits_{t=1}^{q}(x_i\ :\ i\in A_t)\right)\cap\left(\bigcap\limits_\twoline{F\subset[n],\ |F|=q-1}{F\cap A_t\neq\emptyset,\ \forall t}P_{F^c}\right).$$  
\end{Theorem111}

\begin{Theorem111}\label{prim dec final}\cite{OO}
Let $I\subset S$, $I\neq I_{n,q}$, be the final squarefree lexsegment ideal generated in degree $q$, determined by the monomial $u=x_1x_{i_2}\cdots x_{i_q}$, $2\leq i_2<\ldots<i_q\leq n$. Let $F=\{i_2,\ldots,i_q\}$ and $x_F=\prod\limits_{i\in F}x_i$. Then $I$ has the minimal primary decomposition of the form:
$$I=\left(\bigcap\limits_\twoline{G\subset[n],\ |G|=n-q+1}{x_G\geq_{lex}x_{F^c}}P_G\right)\cap\left(\bigcap\limits_\twoline{G\subset[n]\setminus\{1\},\ |G|=n-q+1}{x_{G\setminus\min(G)}\geq_{lex}x_{F^c\setminus\{1\}}}P_G\right)\cap\left(\bigcap\limits_\twoline{G\subset[n],\ |G|=n-q}{x_{F^c\setminus\{1\}}>_{lex}x_G}P_G\right).$$
\end{Theorem111}

The above theorems also allows us to derive formulae for two main invariants of initial and final squarefree lexsegment ideals: the Krull dimension and the depth (Corollary \ref{dim}, Corollary \ref{dim final}). Their description depends on the ends of the lexsegment and on the degree of the monomials of the lexsegment set. By using the standard decompositions of initial/final squarefree lexsegment ideals, we obtain the multiplicity of the corresponding Stanley--Reisner ring (Corollary \ref{multiplicity initial}, respectively Corollary \ref{dim final} (c)). 

It is known that that any completely squarefree lexsegment ideal $I=(L(u,v))$ can be written as the intersection of $(L^i(v))$ with $(L^f(u))$. Using this fact, we are able to compute the standard primary decomposition for completely squarefree lexsegment ideals (Theorem \ref{prim dec completely}). As consequences, we obtain formulae for the Krull dimension and multiplicity of $S/I$, where $I$ is a completely squarefree lexsegment ideal (Corollary \ref{dim+multiplicity completely}).

R. Stanley \cite[Section III.2]{St} defined the concept of sequentially Cohen--Macaulay simplicial complex, a nonpure generalization of the Cohen--Macaulay simplicial complex. A simplicial complex $\Delta$ is \textit{sequentially Cohen--Macaulay} if all its pure skeletons are Cohen--Macaulay. It is known \cite{HH} that the associated Stanley--Reisner ideal $I_\Delta$ is sequentially Cohen--Macaulay, that is, $S/I_{\Delta}$ is a sequentially Cohen--Macaulay module, if and only if $I^\vee=I_{\Delta^\vee}$ is componentwise linear, which means that for all $d\geq 0$, the ideal $I^\vee_{\langle d\rangle}$ generated by all degree $d$ elements in $I^\vee$ has a linear resolution. Here we denoted as usual by $\Delta^\vee$ the Alexander dual of $\Delta$.

As an application of the minimal primary decomposition, in Section $2.2$ we characterize all the completely squarefree lexsegment ideals which are sequentially Cohen--Macaulay (Proposition \ref{sqCM initial}, Proposition \ref{sqCM final} and Theorem \ref{sqCM completely}). 

In \textbf{Chapter $3$} we are interested in computing some other homological invariants for arbitrary squarefree lexsegment ideals. From the previous chapter, it turns out that their computation becomes difficult, even for completely squarefree lexsegment ideals. We begin Chapter $3$ by giving some bounds for the depth of a squarefree lexsegment ideal. Next, we study squarefree lexsgment ideals generated in small degrees. We start by analyzing the $2-$degree case, when a squarefree lexsegment ideal is called a lexsegment edge ideal. We study these ideals and we derive formulae for the Krull dimension, the depth and the regularity. All these invariants have a nice description since their formulae depend only on the ends of the lexsegments, and can be found in \textit{Arithmetical rank of lexsegment edge ideals}, joint work with V. Ene and N. Terai.

\begin{Proposition111}\cite{EOT}
Let $I=(L(u,v))$ be a lexsegment edge ideal which is neither initial nor final and is determined by $u=x_1x_i$ and $v=x_jx_r.$ Then 
$\dim(S/I)=n-j.$
\end{Proposition111}

\begin{Proposition111}\cite{EOT}
Let $I=(L(u,v))$ be a lexsegment edge ideal where $u=x_1x_i, v=x_jx_r, j\geq 2.$ Then
\[
\reg(I)=\left\{ 
\begin{array}{ll}
	3, & \text{ if } i\geq j+2 \text{ and } x_n\not| v \\
	2, & { otherwise.}
\end{array}\right.
\]
\end{Proposition111}

We point out that for the Krull dimension and the regularity, we give proofs different from the original paper.

For lexsegment edge ideals we succeeded to compute the arithmetical rank.  

\begin{Theorem111}\label{aralei}\cite{EOT}
Let $I=(L(u,v))$ be a lexsegment edge ideal. Then 
\[
\ara(I)=\projdim_S(S/I).
\]
\end{Theorem111}

In Section $3.3$ we are interested in analyzing the degree $3$ case, and to observe some similarities with the results obtained in degree $2$. We succeeded to describe the depth of a squarefree lexsegment ideal generated in degree $3$. However, its computation involved various techniques and it turned out to be much more difficult than in degree two case.

\begin{Theorem111}\label{main}\cite{EKOT}
Let $u=x_1x_{i_2}x_{i_3}$ and $v=x_{j_1}x_{j_2}x_{j_3}$ where $j_1\geq 2,$ be two squarefree monomials of degree $3$ and $I=(L(u,v))$ the squarefree 
lexsegment ideal determined by them. Then:
\begin{itemize}
	\item [(a)] $\depth(S/I)=2 $ if  $x_{i_2-1}x_{i_3-1}x_n\geq_{lex} v.$ 
	\item [(b)] $\depth(S/I)=4$  if  $i_2=4$, $i_3\geq 6$ and $j_1=2,\ j_2=3,\  j_3<i_3-1$ or $i_2\geq 5$ and $j_1=2,\ j_2=3,\ i_2-1\leq j_3\leq n-1.$
	\item [(c)] $\depth(S/I)=i_2-j_3+3$ if  $i_2 > 4, j_2=2,$ and $j_3\leq i_2-1.$
	\item [(d)] $\depth(S/I)=3$ in all the other cases. 
\end{itemize}
\end{Theorem111}

The results obtained for the squarefree lexsegment ideals generated in degrees $2$ and $3$ allowed us to formulate some conjectures concerning the depth and the Castelnuovo--Mumford regularity for an arbitrary squarefree lexsegment ideal, which are listed in the end of this chapter.

The \textbf{last chapter} is devoted to the lexsegment ideals in the non-squarefree case. A monomial ideal $I\subset S$ is called a lexsegment ideal if, for each degree $j$, if $I_j\neq 0$, then $I_j$ is generated by a lexsegment set of degree $j$, that is a set of monomials of degree $j$ of the form
	\[\mathcal{L}_j(u,v)=\{w\in\Mon_j(S)\mid u\geq_{lex}w\geq_{lex}v\}.
	\]
for some monomials $u,v$ of degree $j$, $u\geq_{lex}v$.

We prove that, for this class of ideals, being componentwise linear is equivalent to having componentwise linear quotients. Next we consider a smaller class of lexsegment ideals that we call componentwise lexsegment ideals. Let $d$ be the least degree of the minimal monomial generators of the ideal $I$, such that $\mathcal{L}_d(u,v)$ generates $I_d$. For $I$ a componentwise lexsegment ideal, we require that for any $j\geq d+1$, the $j-$degree component $I_j$ is generated over $k$ by the lexsegment set $\mathcal{L}_j(x_1^{j-d}u,x_n^{j-d}v)$. In other words, each higher component $I_j$ is generated over $k$ by the lexsegment set of degree $j$ determined by the largest and the smallest monomial in the lexicographic order in the shadow of the lexsegment set which generates the previous component $I_{j-1}$. We remark that not every lexsegment ideal is componentwise lexsegment. For instance, one may consider the lexsegment ideal $I=(\mathcal L(x_1x_n^{d-1},x_2^d))\subset S$, which it is not a componentwise lexsegment ideal.
 
For a componentwise lexsegment ideal $I$, we show that the property of being componentwise linear is equivalent to the condition that $I_{\langle d\rangle}$ has a linear resolution, where $I_d$ is the first non-zero component of $I$. 

\begin{Theorem111}\label{lex equivalences}\cite{OOS} Let $I$ be a componentwise lexsegment ideal and $d\geq 1$ the lowest degree of the minimal monomial generators of $I$. Let $u,v\in\Mon_d(S)$, $x_1|u$ be such that $I_{\langle d\rangle}=(\mathcal L(u,v))$. The following conditions are equivalent:

\begin{itemize}
	\item[(a)] $I$ is a componentwise linear ideal.
	\item[(b)] $I_{\langle d\rangle}$ has a linear resolution.
	\item[(c)] $I_{\langle d\rangle}$ has linear quotients.
	\item[(d)] $I$ has componentwise linear quotients.
\end{itemize}
\end{Theorem111}

In the last part, we aim at characterizing the lexsegment ideals generated in one degree which are Gotzmann. Initial lexsegment ideals generated in one degree are obviously Gotzmann.

Arbitrary lexsegment ideals generated in one degree which have linear resolutions have been characterized in \cite{ADH}. Their characterization distinguishes between completely lexsegment ideals and those which are not completely lexsegment ideals. In order to characterize the Gotzmann property of a lexsegment ideal generated in one degree, we also make this distinguish. The Gotzmann lexsegment ideals are described in Sections $4.2$ and $4.3$ of this thesis.

\begin{Theorem111}\label{compG}\cite{OOS} Let $u,v\in \Mon_d(S)$, $x_1\mid u$ such that $I=(\mathcal L(u,v))$ is a completely lexsegment ideal of $S$ which is not an initial lexsegment ideal. Let $j$ be the exponent of the variable $x_n$ in $v$ and $a=|\Mon_d(S)\setminus \mathcal L^i(u)|$. The following statements are equivalent:
\begin{itemize}
	\item [(a)] $I$ is a Gotzmann ideal.
	\item [(b)] $a\geq{n+d-1\choose d}-(j+1)$.
\end{itemize}
\end{Theorem111}

\begin{Theorem111}\cite{OOS} Let $u=x_t^{a_t}\cdots x_{n}^{a_n}$, $v=x_t^{b_t}\cdots x_{n}^{b_n}$ be two monomials of degree $d$, $u>_{lex} v$, $a_t\neq0$, $t\geq1$ and $I=(\mathcal L(u,v))$ a lexsegment ideal which is not complete. Then $I$ is a Gotzmann ideal in $S$ if and only if $I=m(x_{l},x_{l+1},\ldots,x_{l+p})$ for some $t\leq l\leq n$, some $1\leq p\leq n-l$ and a monomial $m$.
\end{Theorem111}

The results of this chapter are contained in the paper \textit{Gotzmann lexsegment ideals} \cite{OOS}, joint work with A. Olteanu and L. Sorrenti.

\mainmatter
\pagestyle{fancy}

\chapter{Preliminaries}
In this chapter we collect some fundamental notions and results needed throughout this thesis. We do not give any proof for the well--known results, but we indicate the precise references for all the statements. There will be also two original results which will be accompanied by their proof and references (Lemma \ref{b(d)=c(d)} and Lemma \ref{c^d=0 dnd c<egal b}).

\section{Cohen--Macaulay rings}
In this section, we recall some basic and well--known facts of commutative algebra. We will analyze only some aspects of the theory, while for the proofs we refer the reader to \cite{BH}, \cite{Ei}, \cite{HeHi}, \cite{Mi}, \cite{Vi}. Through this section we make the assumptions that all the rings are commutative, with unity, Noetherian, and all the modules are finitely generated. 

Let $R$ be a ring and $\pp_0\subset \pp_1\subset\ldots\subset \pp_n$ be a chain of prime ideals of $R$, that is a finite strictly increasing sequence of prime ideals. The \textit{length} of the above sequence is $n$. The \textit{height} of a prime ideal $\pp$ of $R$, denoted by $\height(\pp)$, is the supremum of the lengths of all chains of prime ideals 
$$\pp_0\subset \pp_1\subset\ldots\subset \pp_n=\pp$$
which ends at $\pp$. 

For an arbitrary proper ideal $I$ of $R$, the \textit{height} of $I$ is defined as 
$$\height(I)=\min\{\height(\pp):\  \pp\mbox{ is a prime ideal of }R\mbox{ such that }\pp\supset I\}.$$ 


\begin{Definition1}\rm
Let $R$ be a ring. The \textit{Krull dimension} of $R$, denoted by $\dim(R)$, is the supremum of the lengths of all chains of prime ideals in $R$.
\end{Definition1}

\begin{Example1}\rm 
(a) The Krull dimension of any field is zero. 

(b) The polynomial ring $k[x_1,\ldots,x_n]$ is of Krull dimension $n$ and the ideal $I=(x_1,\ldots,x_i)$ is a prime ideal of $\height(I)=i$, for $1\leq i\leq n$.
\end{Example1}

From definition one can derive two immediate consequences. Firstly, one may remark that $\dim(R_\pp)=\height(\pp)$, where $R_\pp$ denotes the localization of $R$ at the multiplicative set $R\setminus \pp$. Secondly, for any ideal $I$ of $R$, one has 
$$\dim(R/I)+\height(I)\leq\dim(R).$$ 

We recall that the \textit{Krull dimension} of an $R-$module $M$ is 
$$\dim(M)=\dim(R/\Ann_R(M)),$$
where $\Ann_R(M)=\{x\in R: xM=0\}$ is the \textit{annihilator} of $M$.

\begin{Definition1}\rm\cite{Ei}
Let $R$ be a ring and $M$ be an $R-$module. An element $x$ of $R$ is a \textit{nonzerodivisor on $M$} if $xm\neq 0$, for all nonzero $m\in M$. If, in addition,
$xM\neq M$, then the element $x$ is said to be \textit{regular on $M$}, or \textit{$M-$regular}. In other words, $x$ is a \textit{nonzerodivisor on $M$} if $\Ann_R(m)=\{x\in R: xm=0\}$ is zero, for all nonzero $m\in M$, or the multiplicative map $M\stackrel{x}{\rightarrow}M$ is injective. 

A \textit{zerodivisor} on $M$ is an element $x$ of $R$ such that there is a nonzero element $m\in M$ with $xm=0$, equivalently the multiplicative map $M\stackrel{x}{\rightarrow}M$ is not injective.
\end{Definition1} 
 
The set of zerodivisors on $M$, where $M$ is a module, is closely related to the set of associated prime ideals $\pp$ of $M$. For an $R-$module $M$, the set of \textit{associated primes of $M$}, denoted by $\Ass_R(M)$, is the set of all prime ideals $\pp\subset R$ such that $\pp=\Ann_R(m)$, for some $m\in M$. The next result gathers all the main properties of the set of associated primes of a module.

\begin{Theorem1}\cite{Ei}
Let $R$ be a Noetherian ring and let $M$ be a finitely generated nonzero $R-$module. 

$(a)$ $\Ass_R(M)$ is a finite, nonempty set of primes, each containing $\Ann_R(M)$. The set $\Ass_R(M)$ includes all the primes minimal among the prime ideals containing $\Ann_R(M)$.
 
$(b)$ The union of the associated primes of $M$ is the set of zerodivisors on $M$.  
\end{Theorem1}

Let $I$ be an ideal of $R$. We recall that an ideal $\pp$ is a \textit{minimal prime ideal} of $I$ if and only if $\pp\supset I$, and there is no prime ideal $\qq\supset I$ which is properly contained in $\pp$. We denote by $\Min(I)$ the set of minimal primes of the ideal $I$. By the above theorem, $\Ass_R(R/I)$ is a finite set containing all minimal prime ideals of $I$.

Recall that an ideal $I$ in a Noetherian ring $R$ is \textit{$\pp-$primary} if $\Ass_R(R/I)=\{\pp\}$.

The notion of $M-$regular element may be extended to the following fundamental definition: 

\begin{Definition1}\rm
A sequence $\xb=x_1, \ldots,x_d$ of elements of $R$ is a \textit{regular sequence} on an $R-$module $M$ if $x_i$ is regular on $M/(x_1, \ldots, x_{i-1})M$, for each $i=1,\ldots, d$.
\end{Definition1}

Equivalently, one may say that the sequence $\xb=x_1, \ldots,x_d$ of elements of $R$ is regular on $M$ if, for every $1\leq i\leq d$, the element $x_i$ is a nonzerodivisor on $M/(x_1,\ldots, x_{i-1})M$ and $\xb M\neq M$. If $\xb$ is regular on $R$ we simply say that $\xb$ is a \textit{regular sequence} in $R$.

An ideal $I\subset R$ generated by a regular sequence in $R$ is called a \textit{complete intersection ideal}.

Since $R$ is Noetherian, any $M-$sequence $\xb=x_1,\ldots,x_d,\ldots$ is finite, since the ascending chain of ideals 
$$(x_1)\subseteq (x_1,x_2)\subseteq\ldots\subseteq (x_1,\ldots,x_d,\ldots)\subseteq\ldots$$
becomes stationary.

Let $I\subset R$ be an ideal with $IM\neq M$. An $M-$sequence $\xb=x_1,\ldots,x_d$ is called \textit{maximal} in $I$ if it cannot be extended to a longer $M-$regular sequence in $I$.

The next result tells us that any regular sequence can be extended to a maximal one.

\begin{Proposition1}\cite{Vi}
Let $M$ be an $R-$module and $I$ be an ideal of $R$ such that $IM\neq M$. If $\xb = x_1,\ldots,x_d$ is an $M-$regular sequence in $I$, then $\xb$ can be extended to a maximal $M-$regular sequence in $I$. 
\end{Proposition1}

All maximal regular sequences have the same length, as it is stated by Rees theorem:

\begin{Theorem1}[Rees, \cite{Mi}]
Let $R$ be a Noetherian ring, $M$ a finitely generated $R-$module, and $I$ an ideal such that $IM\neq M$. All maximal $M-$regular sequences in $I$ have the same length, namely
$$n=\min\{i: \Ext^i_R(R/I,M)\neq 0\}.$$
\end{Theorem1}

\begin{Definition1}\rm
Let $M\neq(0)$ be a finitely generated module over a Noetherian local ring $(R,\mm)$. The \textit{depth of $M$}, denoted by $\depth(M)$, is the length of any maximal regular sequence on $M$ 
which is contained in $\mm$. 
\end{Definition1}

As a consequence of the previous results, in a local ring, the depth may be expressed in terms of the $\Ext$ functor.

\begin{Corollary1}\cite{Mi}
Let $(R,\mm)$ be a Noetherian local ring. If $M$ is a finitely generated $R-$module, then $\depth(M)=\min\{i:\Ext^i_R(R/\mm,M)\neq 0\}$.
\end{Corollary1}

An upper bound for the depth is given in terms of the associated prime ideals.

\begin{Proposition1}\cite{BH}\label{depth-bound}
Let $M$ be a finitely generated module over a Noetherian local ring $(R,\mm)$. Then 
$$\depth(M)\leq \dim(R/\pp), \mbox{ for all } \pp\in\Ass_R(M).$$
In particular, the inequality $$\depth(M)\leq\dim(M)$$
holds.
\end{Proposition1}

An interesting class of modules is the one for which the Krull dimension equals the depth. They are of fundamental importance in commutative algebra and its applications to algebraic geometry and combinatorics.

\begin{Definition1}\rm
Let $M$ be a finitely generated module over the Noetherian local ring $(R,\mm)$. Then $M$ is called \textit{Cohen--Macaulay} if $\dim(M)=\depth(M)$ or if $M=(0)$.

A local ring $(R,\mm)$ is \textit{Cohen--Macaulay} if $R$ is Cohen--Macaulay as an $R-$module.

An ideal $I$ of the local ring $(R,\mm)$ is a \textit{Cohen--Macaulay ideal} if $R/I$ is a Cohen--Macaulay $R-$module.
\end{Definition1}

\begin{Example1}\rm
A classical example of Cohen--Macaulay ring is the polynomial ring $k[x_1,\ldots,x_n]$.
\end{Example1}

\begin{Examples1}\rm 
$(1)$ Let $I\subset S=k[x_1,\ldots,x_5]$ be the ideal $I=(x_1x_2x_3,x_4x_5^2)$. It is easy to see that $\height(I)=2$, thus $\dim(S/I)=3$. One may easily see that $\depth(S/I)=3$, hence the ideal is Cohen--Macaulay.

$(2)$ Let $I=(x_1x_3^2x_4,x_1x_4^2,x_2x_3)$ be an ideal in $S=k[x_1,\ldots,x_4]$. One has $\dim(S/I)=2$ and $\depth(S/I)=1$, hence $I$ is not a Cohen--Macaulay ideal. 
\end{Examples1}

We recall the behaviour of the depth when we take exact sequences of modules. This result will be used several times throughout this thesis.

\begin{Proposition1}[Depth Lemma, \cite{Vi}]
Let $$0\longrightarrow N \longrightarrow M \longrightarrow L\longrightarrow 0$$
be a short exact sequence of modules over a local ring $(R,\mm)$.
\begin{itemize}
	\item [(a)] If $\depth(M)<\depth(L)$, then $\depth(N)=\depth(L)$.
	\item [(b)] If $\depth(M)=\depth(L)$, then $\depth(N)\geq\depth(L)$.
	\item [(c)] If $\depth(M)>\depth(L)$, then $\depth(N)=\depth(L)+1$.
\end{itemize}
\end{Proposition1}

An important result which is used in computing the depth is the Auslander--Buchsbaum theorem.
\begin{Theorem1}[Auslander--Buchsbaum, \cite{Mi}]
Let $(R,\mm)$ be a Noetherian local ring. If $M$ is a nonzero finitely generated module of finite projective dimension, then
$$\projdim_R(M)+\depth(M)=\depth(R).$$
\end{Theorem1}

There is a generalization of the notion of Cohen--Macaulay module due to R. Stanley~\cite{St}. We first need to recall the notions of graded rings and modules.

Let $(G,+)$ be a commutative monoid and $R$ be a ring.

\begin{Definition1}\rm The ring $R$ is \textit{$G-$graded} if it satisfies the following conditions: 
\begin{enumerate}
	\item $R=\bigoplus\limits_{g\in G}R_g$ (direct sum of abelian groups);
	\item $R_{g}R_{h}\subseteq R_{g+h}$, for all $g,\ h\in G$.
\end{enumerate}
\end{Definition1}

\begin{Definition1}\rm Let $R$ be a $G-$graded ring and $M$ be an $R-$module. The module $M$ is called \textit{$G-$graded} if the following hold:
\begin{itemize}
	\item[i)]$M=\bigoplus\limits_{g\in G}M_g$ (direct sum of abelian groups);   
	\item[ii)]$R_{g}M_{h}\subseteq M_{g+h}$, for all $g,\ h\in G$.
\end{itemize}
\end{Definition1}

One calls $M_i$ the \textit{$i-$th homogeneous (or graded) component} of $M$. 

Throughout this thesis we are going to use two gradings on the polynomial ring $S=k[x_1,\ldots,x_n]$.

\begin{Example1}[The standard grading]\rm The polynomial ring $S=k[x_1,\ldots,x_n]$ over a field $k$ is $\mathbb N-$graded, where the $n-$graded component is the $k-$vector space of homogeneous polynomials of degree $n$.
\end{Example1}

\begin{Example1}[The $\mathbb N^n-$grading]\rm The polynomial ring $S=k[x_1,\ldots,x_n]$ over a field $k$ is $\mathbb N^n-$graded, with $S_\alpha=k x^\alpha$ the component of degree $\alpha$, that is the $k-$vector space of basis $x^\alpha$. 
\end{Example1}

Let $R$ be a $G-$graded ring. An ideal $I\subset R$ is called a \textit{graded (or homogeneous) ideal} if it is generated by homogeneous elements. Equivalently, $I$ is a graded ideal if for all $f\in I$, all the homogeneous components of $f$ are in $I$. Monomial ideals are graded ideals of the polynomial ring with respect to the standard grading and also with respect to the $\mathbb N^n-$grading.

Next, we recall the definition of sequentially Cohen--Macaulay module, \cite{St}.
\begin{Definition1}\rm\cite{St}
Let $M$ be a finitely generated $\mathbb Z-$graded module over $S=k[x_1,\ldots,x_n]$. We say that $M$ is \textit{sequentially Cohen--Macaulay} if there exists a finite filtration
$$0=M_0\subset M_1\subset\ldots\subset M_r=M$$
of $M$ by graded submodules $M_i$ satisfying the two conditions:
\begin{itemize}
	\item [(a)] Each quotient $M_i/M_{i-1}$ is Cohen--Macaulay;
	\item [(b)] $\dim(M_1/M_0)<\dim(M_2/M_1)<\ldots<\dim(M_r/M_{r-1})$.
\end{itemize}
\end{Definition1}

The submodules appearing in the filtration of a sequentially Cohen--Macaulay module $M$ were described firstly by P. Schenzel \cite{Sch}, and later by S. Faridi \cite{Fa1}. 

\section{Monomial ideals}
Monomial ideals represent a bridge between commutative algebra and combinatorics. Their importance is given, firstly, by the fact that the monomial ideals appears as initial ideals. By using Gr\"obner bases, one may study arbitrary polynomial ideals via their initial ideals, which are monomial ideals. Many invariants of the polynomial ideals are preserved when passing to the initial ideal or they are bounded by the invariants of the initial ideal. The advantage of passing to initial ideals comes from the fact that monomial ideals have a rich combinatorial structure. The polarization process associates to a monomial ideal a squarefree monomial ideal, which has a beautiful combinatorial interpretation in terms of simplicial complexes. Most of the homological invariants are preserved by the polarization process.

Some results of this section can be found in \cite{HH} in a more general case. However, we present them in the case of monomial ideals.

In the following, we consider a field $k$ and $S=k[x_1,\ldots,x_n]$ the polynomial ring in $n$ variables. For a monomial $u\in S$, we set $\max(u)=\max(\supp(u))$ and $\min(u)=\min(\supp(u))$, where $\supp(u)=\{i:x_i\mid u\}$. We denote by $\frak m=(x_1,\ldots,x_n)$ the maximal graded ideal of $S$. The set of minimal monomial generators of the monomial ideal $I$ is denoted by $G(I)$. 

We recall that a monomial ideal is prime if and only if it is generated by a subset of $\{x_1,\ldots,x_n\}$.

\subsection{Primary decomposition for monomial ideals}
In number theory, the fundamental theorem of arithmetic states that any integer greater than $1$ can be written as a unique product (up to the ordering of factors) of prime numbers. The primary decomposition is intended to be a generalization of this statement. 

The problem of decomposing an ideal into primary ideals is fundamental in commutative algebra. It proves the algebraic foundation for decomposing an algebraic variety into its irreducible components. 
Many authors have tried to develop several algorithms that could be and were implemented in computer algebra systems. Nowadays, in order to compute primary decompositions using modern methods, one may benefit of the support offered by several computer algebra systems \cite{CoCoA}, \cite{GPS}. 

Our aim is to recall some remarkable properties of the primary decomposition for monomial ideals (see for instance \cite{HeHi} and \cite{Vi}). More results concerning the primary decomposition may be found, for example, in \cite{St}, \cite{Vi}.


\begin{Proposition1}\cite{Vi}
Let $I\subset S$ be a monomial ideal. Then every associated prime ideal of $I$ is generated by a subset of variables. 
\end{Proposition1}

A presentation of an ideal $I$ as intersection $I=\bigcap\limits_{i=1}^{r}\qq_i$, where each $\qq_i$ is a primary ideal, is called a \textit{primary decomposition} of $I$. Let $\Ass_S(S/\qq_i)=\{\pp_i\}$. 

\begin{Definition1}\rm
The primary decomposition of $I$ is called \textit{irredundant}, if none of the $\qq_i$ can be omitted in the intersection $I=\bigcap\limits_{i=1}^{r}\qq_i$, and $\pp_i\neq\pp_j$, for all $i\neq j$. If $I=\bigcap\limits_{i=1}^{r}\qq_i$ is an irredundant primary decomposition of $I$, then the $\qq_i$ is called the \textit{$\pp_i-$primary component} of $I$ and $\Ass_S(S/I) =\{\pp_1,\ldots,\pp_r\}$.
\end{Definition1}


The primary monomial ideals were described and their particular form makes them easily to be recognized.

\begin{Proposition1}\cite{Vi}\label{primary}
A monomial ideal $\qq$ is primary if and only if, up to a permutation of variables, it has the form $\qq=(x_1^{a_1},\ldots,x_r^{a_r},x^{b_1},\ldots,x^{b_s})$, where $a_i>0$, for all $1\leq i\leq r$ and $\bigcup\limits_{1\leq j\leq s}\supp(x^{b_j})\subseteq\{1,\ldots,r\}$.
\end{Proposition1}

The class of primary monomial ideals contains a smaller class, namely the irreducible ideals. A monomial ideal $I\subset S$ is called \textit{irreducible} if it cannot be written as $I=J\cap K$, where $J$ and $K$ are monomial ideals in $S$ with the property that $J\supsetneq I$, $K\supsetneq I$. A monomial ideal $I\subset S$ which is not irreducible, it is called \textit{reducible}. Since any irreducible monomial ideal is a primary ideal, we expect to obtain a similar description of irreducible monomial ideals, as for the primary case.

\begin{Proposition1}\cite{Vi}\label{irreducible}
A monomial ideal $I$ is irreducible if and only if, up to a permutation of variables, it can be written as $I=(x_1^{a_1},\ldots,x_r^{a_r})$, where $a_i>0$, for all $1\leq i\leq r$.
\end{Proposition1}

In other words, the irreducible monomial ideals are generated by pure powers of variables.

As a consequence of the above results, it is easy to observe that: 

\begin{Proposition1}\cite{HeHi}\label{p-primary}
The irreducible monomial ideal $(x_{i_1}^{a_1},\ldots,x_{i_s}^{a_s})$ is \newline $(x_{i_1},\ldots,x_{i_s})-$primary.
\end{Proposition1}

There is a fundamental result concerning the decomposition of a monomial ideal into irredundant intersection of irreducible ideals. 

\begin{Theorem1}\cite{Vi}
If $I$ is a monomial ideal in the polynomial ring $S$, then there is a unique irredundant decomposition $I=\qq_1\cap\ldots\cap \qq_r$ such that $\qq_i$ are irreducible monomial ideals. 
\end{Theorem1}

As a consequence, it results that the decomposition of an ideal into irreducible ideals is a primary decomposition. The only problem is that it may not be an irredundant decomposition. In order to get an irredundant decomposition, we let the $\pp$-primary component of $I$ to be the intersection of all irreducible ideals $\qq_i$ which appears in the unique irredundant decomposition of $I$, having $\Ass_S(S/\qq_i)=\{\pp\}$.

One may note that a primary decomposition of a monomial ideal may not be unique. Even so, the primary decomposition obtained from an irredundant intersection of irreducible ideals, by the procedure explained above, is unique. We call it the \textit{standard primary decomposition}.

\begin{Example1}\rm
Let $I=(x^2,xy^3z,z^2y^3)\subset S=k[x,y,z]$ be a monomial ideal. The irredundant presentation of the ideal $I$ as intersection of irreducible ideals is 
$$I=(x^2,y^3)\cap(x,z^2)\cap(x^2,z).$$
One has $\Ass_S(S/(x^2,z))=\Ass_S(S/(x,z^2))=\{(x,z)\}$. We intersect them and we obtain $(x^2,z)\cap(x,z^2)=(x^2,z^2,xz)$ which is $(x,z)-$primary. Therefore, the standard primary decomposition of $I$ is 
$$I=(x^2,y^3)\cap(x^2,z^2,xz).$$
\end{Example1}

In the case of squarefree monomial ideals, the irreducible squarefree monomial ideals are generated by subsets of variables. One may conclude that any squarefree monomial ideal is just an intersection of prime ideals. Moreover, for a squarefree monomial ideal $I$, one has $\Ass_S(S/I)=\Min(I)$.

Algorithms for computing a primary decomposition of a monomial ideal were developed, see for instance \cite{HeHi} and \cite{Vi}.

	\[
\]
\subsection{Ideals with linear quotients}

The class of ideals with linear quotients was introduced by J. Herzog and Y. Takayama \cite{HT} and it is closely related to algebraic combinatorics.

\begin{Definition1}\rm
The monomial ideal $I$ of $S$ has \textit{linear quotients} if the monomials from the minimal monomial set of generators of $I$ can be ordered as $u_1,\ldots, u_s$ such that for all $2\leq i\leq s$ the colon ideals $(u_1,\ldots,u_{i-1}):u_i$ are generated by variables.
\end{Definition1}

There is a very useful equivalent characterization which can be found, for instance in \cite{HH}.

\begin{Lemma1}
A monomial ideal $I\subset S$ has linear quotients with respect to the minimal monomial generators $u_1,\ldots,u_s$ of $I$ if and only if for all $j < i$ there exist an integer $t < i$ and an integer $1\leq l\leq n$ such that
$$\frac{u_t}{\gcd(u_t,u_i)}=x_l\ \mbox{and}\ x_l\ \mbox{divides}\ \frac{u_j}{\gcd(u_j,u_i)}.$$
\end{Lemma1} 

Stable ideals are the most common examples of ideals with linear quotients.
 
The notion of stable ideal was introduced by S. Eliahou and M. Kervaire \cite{EK}. We recall the definition.

\begin{Definition1}\rm
A monomial ideal $I$ of $S$ is called \textit{stable} if for all monomials $w \in I$ and for all positive integers $i$, with $i<\max(w)$, it holds
$$\frac{x_iw}{x_{\max(w)}}\in I.$$	
\end{Definition1}

\begin{Proposition1}\cite{EK}\rm\label{stable min mon gen set}
A monomial ideal $I$ of $S$ is stable if and only if, for every minimal monomial generator $w\in G(I)$ and for all positive integers $i$ with $i<\max(w)$, it holds
$$\frac{x_iw}{x_{\max(w)}}\in I.$$
\end{Proposition1}

In other words, in order to establish the stability of a monomial ideal, we have to check the definition only on the set of minimal monomial generators. The invariants of a stable ideal, such as the Castelnuovo--Mumford regularity, the depth or the Betti numbers, are completely described in \cite{EK}. In \cite{AH}, A. Aramova and J. Herzog computed the resolution of a stable ideal by using the Koszul homology.

There is a squarefree correspondent of stable ideals. The notion of squarefree stable ideal was defined by A. Aramova, J. Herzog and T. Hibi in \cite{AHH} as follows:

\begin{Definition1}\rm
A squarefree monomial ideal $I\subset S$ is called \textit{squarefree stable} if for all monomials $w \in G(I)$ and for all positive integers $j$ with $j<\max(u)$ such that $j\notin\supp(u)$, the monomial $x_ju/x_{\max(u)}\in I$.
\end{Definition1}

Some invariants of squarefree stable ideals are described in \cite{AHH}. It turns out that important homological invariants have similar behaviour as in the case of stable ideals. In \cite{AHH} the explicit minimal free resolution of a squarefree stable ideal is constructed and one may see that it has the same formal structure as the classical Eliahou--Kervaire resolution \cite{EK} of stable monomial ideals. 

In general, the product and the sum of two ideals with linear quotients need not to have again linear quotients.  
The following example is given by A. Conca and J. Herzog in \cite{CH}.

\begin{Example1}\rm\cite{CH}
We consider $S=k[a,b,c,d]$ and the monomial ideals $I=(b,c)$ and $J=(a^2b,abc,bcd, cd^2)$. Then $J$ has linear quotients, and $I$ is generated by a subset of the variables, but the product $IJ$ has no linear quotients (not even a linear resolution).
\end{Example1}

\begin{Example1}\rm
Let $I=(x^3y,x^2y^2,y^3)$ and $J=(x^2,xy^3,y^4)$ be two ideals in $S=k[x,y]$. Even if both ideals have linear quotients, their sum $I+J=(x^2,y^3)$ does not have linear quotients.
\end{Example1}

In \cite{JZ}, A.S. Jahan and X. Zheng proved that:

\begin{Lemma1}\cite{JZ}
Let $I\subset S$ be a monomial ideal. If $I$ has linear quotients, then $\mm I$ has linear quotients, where $\mm=(x_1,\ldots,x_n)$ is the graded maximal ideal of $S$.
\end{Lemma1}

The converse is not true, as it follows from the next example.

\begin{Example1}\rm
Let $I =(a^2, bc)\subset k[a,b,c]$. Then the ideal 
$$\mm I=(a^3, a^2b,a^2c,abc,b^2c,bc^2)$$
has linear quotients with respect to the given order of the minimal monomial generators, but $I$ has no linear quotients.
\end{Example1}

For monomial ideals, one can generalize the notion of linear quotients \cite{JZ}.

Let $I$ be a monomial ideal of $S$. We denote by $I_{\langle j\rangle}$, the ideal generated by all homogeneous monomials of degree $j$ belonging to $I$.

\begin{Definition1}\rm\cite{JZ}
A monomial ideal $I\subset S$ has \textit{componentwise linear quotients} if, for all $j$, the ideal $I_{\langle j\rangle}$ has linear quotients.
\end{Definition1}

\begin{Example1}\rm
We consider the ideal $I=(a^2,b)\subset k[a,b,c]$. Then $I_{\langle 1\rangle}=(b)$ is generated by a variable and $I_{\langle 2\rangle}=(a^2,ab,b^2,bc)$ has linear quotients in the given order. Thus the ideal $I$ has componentwise linear quotients.
\end{Example1}

The connection between linear quotients and componentwise linear quotients is given in \cite{JZ}.

\begin{Theorem1}\cite{JZ}\label{lin-quot-comp}
Let $I\subset S$ be a monomial ideal. If $I$ has linear quotients, then $I$ has componentwise linear quotients.
\end{Theorem1}

	\[
\]

\subsection{Ideals with a linear resolution}
Let $S=k[x_1,\ldots,x_n]$ be the polynomial ring in $n$ variables over a field $k$. Following \cite{EH}, we recall some invariants which can be expressed by the graded Betti numbers of a finitely generated graded $S-$module.

Let $M$ be a finitely generated graded $S-$module. It is known that there exists a unique (up to isomorphism) minimal graded free resolution of $M$:
$$0\rightarrow \bigoplus\limits_{j}S(-j)^{\beta_{s,j}}\rightarrow\ldots\rightarrow \bigoplus\limits_{j}S(-j)^{\beta_{1,j}}\rightarrow \bigoplus\limits_{j}S(-j)^{\beta_{0,j}}\rightarrow M\rightarrow 0.$$ 

One has $\beta_{i,j}=\dim_k\Tor_i(k,M)_j$, for all $i,\ j$. The numbers $\beta_{i,j}$ are called the \textit{graded Betti numbers} of $M$. 

We have
$$\projdim_S(M)=\max\{i: \beta_{i,j}\neq 0, \mbox{ for some }j\}.$$ 
The number $$\reg(M)=\max\{j:\beta_{i,i+j}\neq 0, \mbox{ for some }i\}$$ is called the \textit{Castelnuovo--Mumford regularity} of $M$.

The module $M$ has a \textit{$d-$linear resolution} if the minimal graded free resolution
of $I$ is of the form
$$0\rightarrow S(-d-s)^{\beta_s}\rightarrow\ldots\rightarrow S(-d-1)^{\beta_1}\rightarrow S(-d)^{\beta_0}\rightarrow M\rightarrow0.$$  

Equivalently, $M$ has a $d-$linear resolution if and only if $\reg(M)=d$.

There are known some remarkable results concerning the behaviour of the regularity with respect to exact sequences, see for instance \cite{CH}.

Consider the short exact sequence of graded $S-$modules 
$$0\rightarrow N\rightarrow M \rightarrow P\rightarrow 0.$$
This yields the long exact sequence:
$$\ldots\rightarrow \Tor_{i+1}(P,k)\rightarrow \Tor_i(N,k)\rightarrow \Tor_i(M,k)\rightarrow \Tor_i(P,k)\rightarrow\ldots.$$

It follows that
\[
\begin{array}{ccc}
\reg(M)&\leq&\max\{\reg(P),\reg(N)\}\qquad\\
\reg(N)&\leq&\max\{\reg(M),\reg(P)+1\}\\
\reg(P)&\leq&\ \max\{\reg(N)-1,\reg(M)\}.\\
\end{array}
\]

In particular, one may obtain the results given in \cite{Ro} and \cite{O}:

\begin{Proposition1}\cite{Ro}, \cite{O}
Let $S=k[x_1,\ldots,x_n]$ be the polynomial ring and
$$0\rightarrow M'\rightarrow M \rightarrow M''\rightarrow 0$$
be an exact sequence of graded $S-$modules. 
\begin{itemize} 
\item [(a)] If $M'$ and $M''$ have $d-$linear resolutions, then $M$ has a $d-$linear resolution.
\item [(b)] If $M'$ has a $(d+1)-$linear resolution and $M$ has a $d-$linear resolution, then $M''$ has a $d-$linear resolution.
\end{itemize}
\end{Proposition1}

Monomial ideals with linear quotients generated in one degree are examples of ideals with a linear resolution. In particular, (squarefree) stable ideals generated in one degree have a linear resolution.

\begin{Proposition1}\cite{CH}\label{linquot->linres}
Let $I\subset S$ be a monomial ideal generated in degree $d$ and assume that $I$ has linear quotients. Then $I$ has a $d-$linear resolution.
\end{Proposition1}

J. Herzog, T. Hibi and X. Zheng proved that the converse also holds for monomial ideals generated in degree $2$:

\begin{Proposition1}\cite{HHZ}
Let $I$ be a monomial ideal generated in degree $2$. The following conditions are equivalent:
\begin{itemize}
	\item [(a)] I has a linear resolution.
	\item [(b)] I has linear quotients.
\end{itemize}
\end{Proposition1}

It naturally arises the generalization of the notion of linear resolution, when considering monomial ideals generated in different degrees.

\begin{Definition1}\rm\cite{HH}
A monomial ideal $I\subset S$ is \textit{componentwise linear} if, for all $j$, the ideal $I_{\langle j\rangle}$ has a linear resolution.
\end{Definition1}

The notion of componentwise linear ideals was introduced by J. Herzog and T. Hibi \cite{HH} in order to generalize the Eagon--Reiner theorem. Since in the monomial squarefree case the notion of linear resolution corresponds, by using Alexander duality, to the Cohen--Macaulay concept, the natural extension is that to componentwise linearity it corresponds the sequentially Cohen--Macaulayness, again by Alexander duality \cite{HRW}. Remarkable properties of componentwise linear ideals are described in \cite{HeHi}.

\begin{Example1}\rm
Any stable ideal is componentwise linear.
\end{Example1}

Denote, as usual, $\mm=(x_1,\ldots,x_n)$ the maximal graded ideal of $S$.

\begin{Theorem1}\cite{HRW}
If $I\subset S$ is a monomial ideal with a linear resolution, then $\mm I$ has again a linear resolution.
\end{Theorem1}

Monomial ideals with linear resolution are componentwise linear.

\begin{Theorem1}\cite{HH}
Let $I\subset S$ be a monomial ideal. If $I$ has a linear resolution, then $I$ is componentwise linear.
\end{Theorem1}

The converse of this result is not true. Stable ideals, which are not generated in one degree, are examples of componentwise linear ideals which do not have a linear resolution.


The next result comes as an extension of Proposition \ref{linquot->linres}.

\begin{Theorem1}\cite{CH}\label{lin quot comp lin}
Let $I\subset S$ be a monomial ideal which has linear quotients with respect to a minimal homogeneous system of generators of $I$. Then $I$ is componentwise linear.
\end{Theorem1}

The properties of the monomial ideals generated in one degree can be described in the next diagram:

	\[\begin{array}{ccc}
	\mbox{linear quotients} & \Longrightarrow & \mbox{linear resolution}\\
	 & & \\
	\Downarrow & & \Downarrow \\
	& & \\
	\mbox{componentwise linear}& \Longrightarrow & \mbox{componentwise linear}\\
	\mbox{quotients} & & \\
\end{array}
\]





Next, we focus on the squarefree monomial ideals. Let $I\subset S$ be a squarefree monomial ideal. For each degree, we write $I_{[j]}$ for the ideal generated by all the squarefree monomials of degree $j$ belonging to the ideal $I$. Then one has the squarefree analogue of componentwise linear notion as follow:
 
\begin{Definition1}\rm
Let $I$ be a squarefree monomial ideal in $S$. Then $I$ is called \textit{squarefree componentwise linear} if for all $j\geq 0$, the ideals $I_{[j]}$ have a linear resolution. 
\end{Definition1}


\begin{Proposition1}\cite{HeHi}
A squarefree monomial ideal $I\subset S$ is componentwise linear if and only if it is squarefree componentwise linear. 
\end{Proposition1}

	\[
\]

\subsection{Lexsegment ideals and squarefree lexsegment ideals}
The lexsegment ideals play a key role in the theory of Hilbert function and in extremal combinatorics. 

A general class of lexsegment ideals was defined by H. Hulett and H.M. Martin in \cite{HM}.

Let $k$ be a field and $S=k[x_1,\ldots,x_n]$ be the polynomial ring in $n$ variables. We assume that all the monomials of $S$ are ordered by the lexicographical order with $x_1>x_2>\ldots>x_n$. Namely, if $u=x_1^{a_1}\cdots x_n^{a_n}$ and $v=x_1^{b_1}\cdots x_n^{b_n}$ are two monomials in $S$, one has $u>_{lex}v$ if $\deg(u)>\deg(v)$ or $\deg(u)=\deg(v)$ and there exists $1\leq s\leq n$ such that $a_i=b_i$ for all $i\leq s-1$ and $a_s>b_s$. 





\begin{Definition1}\rm
Let $u$ and $v$ be monomials in $\Mon_d(S)$, $u\geq_{lex}v$. The set of monomials of degree $d$
$$\mathcal L(u,v)=\{w\in\Mon_d(S)\ :\ u\geq_{lex} w\geq_{lex} v\}$$
is called the \textit{lexsegment set} determined by the monomials $u$ and $v$.
A \textit{lexsegment ideal} is a monomial ideal generated by a lexsegment set. 
\end{Definition1}

\begin{Remark1}\rm 
In the particular case when $u=x_1^d$, one has 
$$\mathcal L(x_1^d,v)=\mathcal L^i(v)=\{w\in\Mon_d(S)\ :\ w\geq_{lex} v\}$$
and it is called an \textit{initial lexsegment}. Similarly, when the monomial $v=x_n^d$, we get the set
$$\mathcal L(u,x_n^d)=\mathcal L^f(u)=\{w\in\Mon_d(S)\ :\ u\geq_{lex} w\}$$
which is called a \textit{final lexsegment} set. An \textit{initial (final) lexsegment ideal} is a monomial ideal generated by an initial (final) lexsegment set. 
\end{Remark1}

\begin{Example1}\rm
Let $u=x_1x_2^2$, $v=x_2^2x_3$. Then the lexsegment ideal generated by the lexsegment $L(u,v)\subset k[x_1,x_2,x_3]$ is
$$I=(L(u,v))=(x_1x_2^2,x_1x_2x_3,\ x_1x_3^2,\ x_2^3,\ x_2^2x_3).$$
\end{Example1}

In \cite{EOS} it was computed the Krull dimension and the depth of lexsegment ideals in numerical terms related to the ends of the generating set. We recall here only the characterization of lexsegment ideals with depth zero since we are going to use it in the last chapter.

\begin{Proposition1}\cite{EOS}\label{depthzero} Let $I=(\mathcal L(u,v))\subset S$, where $u=x_1^{a_1}\cdots x_n^{a_n}$, $v=x_q^{b_q}\cdots x_n^{b_n}$, $q\geq2$, $a_1,b_q>0$. Then $\depth(S/I)=0$ if and only if $x_nu/x_1\geq_{lex} v$.
\end{Proposition1}

In particular, the initial and final lexsegment ideals have depth zero \cite{O1}.

The other possible values of the depth for lexsegment ideals are given in \cite{EOS}. 

After computing the Krull dimension and the depth, one may proceed at characterizing the Cohen--Macaulay lexsegment ideals. The classification can be found in \cite{EOS}.

A very useful tool in the study of sets of monomials is the notion of shadow. 
\begin{Definition1}\rm
Let $L$ be a set of monomials in $S$. The set 
$$\shad(L)=\{x_1,\ldots,x_n\}L$$
is called the \textit{shadow} of $L$. 
\end{Definition1}

One can define recursively the $i$\textit{-th shadow} as $\shad^i(L)=\shad(\shad^{i-1}(L))$.

Any initial lexsegment has the property that its shadows are again initial lexsegments, a fact which is not true for arbitrary lexsegment sets.

\begin{Definition1}\rm\cite{HM}
A \textit{completely lexsegment set} is a lexsegment whose all iterated shadows are lexsegments. 
A monomial ideal generated by a completely lexsegment set is called a \textit{completely lexsegment ideal}.
\end{Definition1}

A classical example of completely lexsegment ideal is that of initial lexsegment ideal. 

By definition, the procedure which describes the completely lexsegment sets has an infinite number of steps.
Even so, the problem of determining whether a lexsegment set is a completely lexsegment set or not can be reduced to only one step, namely to the computation of the first shadow of the lexsegment set which generates the ideal, according to the Persistence Theorem \cite{DH}.

E. De Negri and J. Herzog characterized the completely lexsegment ideals \cite{DH}. 

One may easily note that any initial lexsegment ideal is stable in the sense of Eliahou--Kervaire. Moreover, final lexsegment ideals are stable, but with respect to the order of the variables $x_n >\ldots> x_1$. Hence initial and final lexsegment ideal generated in one degree have linear quotients, therefore a linear resolution.

The lexsegment ideals with a linear resolution were classified by A. Aramova, E. De Negri and J. Herzog in \cite{ADH}. In the characterization, the authors distinguish between the completely lexsegment ideals and lexsegment ideals which are not complete. 

\begin{Theorem1}\cite{ADH}\label{completelylex}
Let $u=x_1^{a_1}\cdots x_n^{a_n}$ and $v=x_1^{b_1}\cdots x_n^{b_n}$ be monomials of degree $d$ in $S$, $u\geq_{lex}v$ such that $I=(\mathcal L(u,v))$ is a completely lexsegment ideal. Then $I$ has a linear resolution if and only if one of the following conditions holds:
\begin{itemize}
	\item[(a)]  $u=x_1^px_2^{d-p}$ and $v=x_1^px_{n}^{d-p}$, for some integer $0<p\leq d$.
	\item[(b)] $b_1<a_1-1$.
	\item[(c)] $b_1=a_1-1$ and for the greatest monomial $w<_{lex}v$, $w\in\Mon_d(S)$, one has $x_1w/x_{\max(w)}\leq_{lex} u$.
\end{itemize}
\end{Theorem1}

\begin{Theorem1}\cite{ADH}\label{noncompletelylex}
Let $u=x_1^{a_1}\cdots x_n^{a_n}$ and $v=x_2^{b_2}\cdots x_n^{b_n}$ be monomials of degree $d$ in $S$, $a_1\neq0$, $u\geq_{lex}v$. Suppose that $I=(\mathcal L(u,v))$ is not a completely lexsegment ideal. Then $I$ has a linear resolution if and only if $u$ and $v$ are of the form
	\[u=x_1x_{l+1}^{a_{l+1}}\cdots x_n^{a_n},\ v=x_lx_n^{d-1},
\]
for some $l$, $2\leq l<n$.
\end{Theorem1}

In general, a monomial ideal which has a linear resolution may not have linear quotients. For the case of lexsegment ideals these two concepts are equivalent \cite{EOS}.

\begin{Theorem1}\cite{EOS}\label{lin res lin quot lex}
Let $u=x_1^{a_1}\cdots x_n^{a_n},$ with  $a_1>0,$ and $ v=x_1^{b_1}\cdots x_n^{b_n}$ be monomials of degree $d$ with $u\geq_{lex} v,$ such that $I=(\mathcal L(u,v))$ is a lexsegment ideal. Then $I$ has a linear resolution if and only if $I$ has linear quotients.
\end{Theorem1}

M. Ishaq, in \cite{I}, contributed to the study of this class of ideals, by describing all the associated prime ideals of a lexsegment ideal. Moreover, he proved that all the lexsegment ideals are sequentially Cohen--Macaulay.   

In the following, we aim at analyzing the squarefree analogue of lexsegment ideals. 

Let $S=k[x_1,\ldots,x_n]$ be the polynomial ring in $n$ variables over a field $k$.

For an integer $d\geq 2$, let $\Mon_d^s(S)$ be the set of all squarefree monomials of degree $d$ in $S$. We consider the lexicographical order on the monomials in $S$ with $x_1>\ldots>x_n$. 

\begin{Definition1}\rm
Let $u,\ v\in\Mon_d^s(S)$ be two squarefree monomials of degree $d$, $u>_{lex}v$. 
The set $L(u,v)=\{w\in\Mon_d^s(S): u\geq_{lex}w\geq_{lex} v\}$ is called the \textit{squarefree lexsegment set} determined by $u$ and $v$.
\end{Definition1}

\begin{Remark1}\rm
As in the non-squarefree case, we have two particular squarefree lexsegment sets. One of them is the \textit{initial squarefree lexsegment set} determined by $v$, that is a set of the form $L^i(v)=\{w\in\Mon_d^s(S):w\geq_{lex} v\}$, and the other one is $L^f(u)=\{w\in\Mon_d^s(S):u\geq_{lex}w\}$, which is a \textit{final squarefree lexsegment set} determined by $u$. 
\end{Remark1}

\begin{Definition1}\rm
An \textit{(initial, final) squarefree lexsegment ideal} is the squarefree monomial ideal generated by an (initial, final) squarefree lexsegment set. 
\end{Definition1}

The concept of squarefree lexsegment ideal was introduced in \cite{AHH} by A. Aramova, J. Herzog and T. Hibi, where the notion was associated with the nowadays concept of initial squarefree lexsegment ideals. This class was also studied in \cite{ADH}, \cite{AHH}, \cite{B} and \cite{BoS}. 

\begin{Example1}\rm
If $n=6$ and $d=3$, we have 
$$\Mon_3^s(S)=\{x_1x_2x_3,x_1x_2x_4,x_1x_2x_5,x_1x_2x_6,x_1x_3x_4,x_1x_3x_5,x_1x_3x_6,x_1x_4x_5,x_1x_4x_6,$$
$$\qquad\qquad\qquad x_1x_5x_6,x_2x_3x_4,x_2x_3x_5,x_2x_3x_6,x_2x_4x_5,x_2x_4x_6,x_2x_5x_6,x_3x_4x_5,x_3x_4x_6,$$
$$x_3x_5x_6,x_4x_5x_6\}.\qquad\qquad\qquad\qquad\qquad\qquad\qquad\qquad\ \ \ $$
If we take $u=x_1x_2x_5$ and $v=x_3x_4x_5$, then the squarefree lexsegment set determined by $u$ and $v$ is  
$$L(u,v)=\{x_1x_2x_5,\ x_1x_2x_6,\ x_1x_3x_4,\ x_1x_3x_5,\ x_1x_3x_6,\ x_1x_4x_5,\ x_1x_4x_6,\ x_1x_5x_6,$$
$$\qquad\ \ x_2x_3x_4,\ x_2x_3x_5,\ x_2x_3x_6,\ x_2x_4x_5,\ x_2x_4x_6,\ x_2x_5x_6,\ x_3x_4x_5\}.\ $$
\end{Example1}

In the rest of this section, we will refer only to the squarefree lexsegment ideals, even we do not explicitly mention it.

The formal analogue of the notion of shadow is the squarefree shadow.

\begin{Definition1}\rm
Let $L$ be a subset of $\Mon_d^s(S)$. The set 
$$\shad_{s}(L)=\{x_iu\ : u\in L,\ x_i\nmid u,\ 1\leq i\leq n\}$$
is called the \textit{squarefree shadow} of $L$. 
\end{Definition1}

\begin{Definition1}\rm
A squarefree lexsegment $L$ is called a \textit{completely squarefree lexsegment} if all the iterated squarefree shadows of $L$ are again squarefree lexsegments.  
\end{Definition1}

Once again, one may define the \textit{$i$-th squarefree shadow} recursively by $\shad_{s}^i(L)=\shad_{s}(\shad_{s}^{i-1}(L))$. Initial squarefree lexsegment ideals are examples of completely squarefree lexsegment ideals. 

In \cite{B}, V. Bonanzinga proved a persistence theorem for squarefree lexsegment ideals, which allows us to determine the completely squarefree lexsegment ideals by looking only at the first squarefree shadow.

\begin{Example1}\rm
The squarefree lexsegment set $L(u,v)=\{x_1x_4x_5,x_2x_3x_4\}$ in the polynomial ring $k[x_1,\ldots,x_5]$ is not completely. Indeed, 
$$\shad_{s}(L(u,v))=\{x_1x_2x_3x_4,\ x_1x_2x_4x_5,\ x_1x_3x_4x_5,\ x_2x_3x_4x_5\}$$
is not a lexsegment, since $x_1x_2x_3x_5\notin\shad_{s}(L(u,v))$. 
\end{Example1}

\begin{Lemma1}\cite{B}
Let $u$ and $v$ be two monomials in $\Mon_d^s(S)$, $u\geq_{lex}v$. We consider the following squarefree lexsegment ideals: $I = (L(u,v))$, $J = (L^f(u))$ and $K= (L^i(v))$. If $I$ is a completely lexsegment ideal, then $I =J\cap K$.
\end{Lemma1}

A characterization of completely squarefree lexsegment ideals is given in \cite{B}, \cite{BoS}. 


Initial and final squarefree lexsegment ideals generated in degree $d$ have a $d-$linear resolution, being squarefree stable \cite{AHH}.

Moreover, the arbitrary squarefree lexsegment ideals with a linear resolution were characterized in \cite{B}, \cite{BoS}. Similarly to the non-squarefree case, the characterization is done in two steps.



A useful property of the non-squarefree lexsegment ideals is that having a linear resolution is equivalent with having linear quotients. For the squarefree lexsegment ideals, a partial similar result is proved in \cite{BEOS}.  

\begin{Proposition1}\cite{BEOS}
Let 
$$u=x_1\cdots x_hx_{i_{h+1}}\cdots x_{i_d}\geq_{lex}v=x_1\cdots x_{h-1}x_{j_{h}}x_{j_{h+1}}\cdots x_{j_d}$$
be two squarefree monomials in $S$ with $j_h\neq h$ and suppose that $I=(L(u,v))$ is not a completely squarefree lexsegment ideal. Then $I$ has a linear resolution if and only if it has linear quotients. 
\end{Proposition1}

For completely squarefree lexsegment ideals, it is not known if having a linear resolution is equivalent to having linear quotients.

	\[
\]

\subsection{Gotzmann ideals}
Important informations about standard graded $k-$ algebras can be found by computing its Hilbert function. 
	\[
\]

\subsubsection{\rm\textbf{Binomial representations}}

Let $a$ and $d$ be two positive integers and 
$$a={a_d\choose d}+{a_{d-1}\choose d-1}+\cdots+{a_j\choose j}$$
be the $d-$binomial representation of $a$, with $a_d>a_{d-1}>\cdots> a_j\geq j\geq 1$.

One can define several operators related to the binomial expansion of an integer, as follows:
\[\large{
	\begin{array}{ccccccccc}
	a^{\langle d\rangle}&=& {a_d+1\choose d+1}&+&{a_{d-1}+1\choose d}&+&\cdots &+&{a_j+1\choose j+1}.\\
	&&&&&&&&\\
	a_{\langle d\rangle}&=&{a_d\choose d-1}&+&{a_{d-1}\choose d-2}&+&\cdots &+&{a_j\choose j-1}.\\
	&&&&&&&&\\
  a^{(d)}&=&{a_d\choose d+1}&+&{a_{d-1}\choose d}&+&\cdots &+&{a_j\choose j+1}.\\
	\end{array}}
\]

We set, by convention, 
$$0^{\langle d\rangle}=0_{\langle d\rangle}=0^{(d)}=0$$
$$1^{\langle d\rangle}=1 \mbox{ and }1_{\langle d\rangle}=0.$$

\begin{Example1}\rm
Let $a=148$ and $d=5$. Then the $5-$binomial representation of $148$ is 
$$148={9\choose 5}+{6\choose 4}+{4\choose 3}+{3\choose 2}=126+15+4+3.$$
Using the above notations, one have the following:
$$148^{\langle 5\rangle}={9+1\choose 5+1}+{6+1\choose 4+1}+{4+1\choose 3+1}+{3+1\choose 2+1}=210+21+5+4=240.$$
$$148_{\langle 5\rangle}={9\choose 5-1}+{6\choose 4-1}+{4\choose 3-1}+{3\choose 2-1}=126+20+6+3=154.$$
$$148^{(5)}={9\choose 5+1}+{6\choose 4+1}+{4\choose 3+1}+{3\choose 2+1}=84+6+1+1=92.$$
\end{Example1}

\begin{Definition1}\rm
Let $a=(a_1,\ldots,a_d)$ and $b=(b_1,\ldots,b_d)$ be two elements in $\mathbb Z^d_{\geq0}$. We define the \textit{lexicographical order} on $\mathbb Z^d_{\geq 0}$ as follows: $a>_{lex}b$ if and only if the leftmost nonzero entry of the vector $a-b$ is positive. Equivalently, there exists an integer $i\geq 1$ such that $a_s=b_s$, for all $s< i$ and $a_i>b_i$.
\end{Definition1}

\begin{Lemma1}\cite{BH}\label{h<h'}
Let $a,\ a'$ and $d$ be positive integers and
$$a={a_d\choose d}+{a_{d-1}\choose d-1}+\cdots+{a_j\choose j}$$
$$a'={a'_d\choose d}+{a'_{d-1}\choose d-1}+\cdots+{a'_j\choose j}$$
be the $d-$binomial representations of $a$ and $a'$. We set $a_s=0$, for all $1\leq s<i$ and $a'_t=0$, for all $1\leq t<j$. Then $a<a'$ if and only if 
$$(a_d,a_{d-1},\ldots,a_1)<_{lex}(a'_d,a'_{d-1},\ldots,a'_1).$$
\end{Lemma1}

We list some useful properties of the above operators.
 
\begin{Lemma1}\cite{HeHi}
Let $a$ and $d$ be two positive integers. Then:
$$a^{\langle d\rangle}=a+a^{(d)}.$$
\end{Lemma1}

\begin{Lemma1}\cite{HeHi}
Let $a\geq a'$ and $d$ be positive integers. Then $a^{\langle d\rangle}\geq a'^{\langle d\rangle}$
\end{Lemma1}

\begin{Lemma1}\label{b(d)=c(d)}\cite{OOS} Let $c>b>0$ be two integers. Let
	$b={b_d\choose d}+\ldots +{b_j\choose j},$ $b_d>b_{d-1}>\ldots>b_j\geq j\geq1$, and $c={c_d\choose d}+\ldots +{c_i\choose i},$ $c_d>c_{d-1}>\ldots>c_i\geq i\geq1$, be the $d-$binomial expansions of $b$ and $c$. The following statements are equivalent:
\begin{itemize}
	\item[(i)] $b^{(d)}=c^{(d)}$.
	\item[(ii)] $j\geq2$ and $c-b\leq j-1$. 
\end{itemize}
\end{Lemma1}

\begin{proof} Let $b^{(d)}=c^{(d)}$. Since $c>b$, by Lemma \ref{h<h'}, there exists $s\leq d$ such that $c_d=b_d,\ldots,c_{s+1}=b_{s+1}$, and $c_s>b_s$. We obviously have $s+1\geq j$. Let us suppose that $s\geq j$. Since $c_s\geq b_s+1$, we get:
	\[{c_s\choose s+1}\geq{b_s+1\choose s+1}\geq{b_s\choose s+1}+{b_{s-1}\choose s}+\ldots+{b_j\choose j+1}+{b_j\choose j}>\]
	\[>{b_s\choose s+1}+{b_{s-1}\choose s}+\ldots+{b_j\choose j+1}
\]
This leads to the inequality $c^{(d)}>b^{(d)}$, which contradicts our hypothesis. Indeed, we have
	\[c^{(d)}\geq{c_d\choose d+1}+\ldots+{c_{s+1}\choose s+2}+{c_s\choose s+1}>\]
  \[>{b_d\choose d+1}+\ldots+{b_{s+1}\choose s+2}+{b_s\choose s+1}+\ldots+{b_j\choose j+1}=b^{(d)}.\]
Therefore we must have $s=j-1$. Hence $j\geq2$ and $c$ has the binomial expansion
	\[c={c_d\choose d}+\ldots+{c_j\choose j}+{c_{j-1}\choose j-1}+\ldots+{c_i\choose i}.
\]
Using the equality $c^{(d)}=b^{(d)}$ we get
	\[{c_{j-1}\choose j}+\ldots+{c_i\choose i+1}=0,
\]
which implies that $c_{j-1}=j-1,\ldots,c_i=i$. Therefore $c=b+j-i\leq b+j-1$, which proves (ii).
 
Now, let $j\geq2$ and $c\leq b+j-1$. As in the first part of the proof, let $s\leq d$ be an integer such that $c_d=b_d,\ldots,c_{s+1}=b_{s+1}$, and $c_s>b_s$. If $s\geq j$, we get the following inequalities:
	\[c={c_d\choose d}+\ldots+{c_{s+1}\choose s+1}+{c_s\choose s}+\ldots+{c_i\choose i}\geq\]
	\[\geq{b_d\choose d}+\ldots+{b_{s+1}\choose s+1}+{b_s+1\choose s}+{c_{s-1}\choose s-1}+\ldots+{c_i\choose i}\geq\]
	\[\geq{b_d\choose d}+\ldots+{b_{s+1}\choose s+1}+{b_{s}\choose s}+\ldots+{b_j\choose j}+{b_j\choose j-1}+{c_{s-1}\choose s-1}+\ldots+{c_i\choose i}=
\]
	\[=b+{b_j\choose j-1}+{c_{s-1}\choose s-1}+\ldots+{c_i\choose i}\geq b+j-i+s.
\]
Since, by hypothesis, $c-b\leq j-1$, we have $j-1\geq j-i+s$, thus $s\leq i-1$, a contradiction. Hence, $s=j-1$. Then we have

	\[c^{(d)}={c_d\choose d+1}+\ldots+{c_{s+1}\choose s+2}+{c_s\choose s+1}+\ldots+{c_{i}\choose i+1}=
\]
	\[={b_d\choose d+1}+\ldots+{b_{j}\choose j+1}+{c_{j-1}\choose s}+\ldots+{c_{i}\choose i+1}=\]\[=b^{(d)}+{c_{j-1}\choose j}+\ldots+{c_{i}\choose i+1}
\]
If we assume that $c_{j-1}\geq j$, then it follows that ${c_{j-1}\choose j-1}\geq j$. Looking at the $d$-binomial expansions of $b$ and $c$, we get $c-b\geq j$, contradiction. Hence $c_{j-1}=j-1$. This equality implies also the equalities $c_k=k$, for all $i\leq k\leq j-2$. We obtain the following binomial expansion of $c$:
	\[c={c_d\choose d}+\ldots+{c_{j}\choose j}+{j-1\choose j-1}+\ldots+{i\choose i}.
\]
Then
\[c^{(d)}={c_d\choose d+1}+\ldots+{c_{j}\choose j+1}={b_d\choose d+1}+\ldots+{b_{j}\choose j+1}=b^{(d)}.
\] 
\end{proof}

\begin{Lemma1}\label{c^d=0 dnd c<egal b}\cite{OOS}
Let $c>0$ be an integer with the binomial expansion
\[c={c_d\choose d}+\ldots+{c_i\choose i},c_d>\ldots>c_i\geq i\geq 1.
\]
The following statements are equivalent:
\begin{itemize}
	\item [(a)] $c^{(d)}=0$.
	\item [(b)] $c\leq d$. 
\end{itemize}
\end{Lemma1}

\begin{proof}
Let $c\leq d$. Then $c$ has the following binomial expansion with respect to $d$:
\[c={d\choose d}+\ldots+{i\choose i}, \mbox{ for some }i\geq 1.
\]
Hence $c^{(d)}=0$.

Now let $c^{(d)}=0$. We get
\[{c_d\choose d+1}+\ldots+{c_i\choose i+1}=0,\mbox{ which implies}\]
$$c_d=d,\ldots,\ c_i=i.$$

It follows $c=d-(i-1)\leq d$.
\end{proof}

Binomial expansions naturally appear in the context of lexsegment sets. The definition of the binomial operator $a^{\langle d\rangle}$ is justified by the following
result:

\begin{Proposition1}\cite{HeHi}\label{cardinality mon-shad}
Let $L\subset\Mon_d(S)$ be an initial lexsegment set of monomials and $a=|\Mon_d(S)\setminus L|$. Then
$$|\Mon_{d+1}(S)\setminus \shad(L)|=a^{\langle d\rangle}.$$
\end{Proposition1}

In 1927, F.S. Macaulay \cite{Mac} characterized the possible Hilbert functions of standard graded $k-$algebras.

\begin{Theorem1}[Macaulay, \cite{HeHi}]
Let $h:\mathbb Z_+\longrightarrow\mathbb Z_+$ be a numerical function. The following conditions are equivalent:
\begin{itemize}
	\item [(a)] $h$ is the Hilbert function of a standard graded $k-$algebra.
	\item [(b)] There exists an integer $n\geq 1$ and $I\subset k[x_1,\ldots,x_n]=S$ an initial lexsegment ideal such that $H(S/I,i)=h(i)$, for all $i\geq 0$. 
	\item [(c)] $h(0)=1$ and $h(j+1)\leq h(j)^{\langle j\rangle}$, for all $j>0$.
\end{itemize}
\end{Theorem1}

	\[
\]
\subsubsection{\rm\textbf{Gotzmann ideals}}

Given a graded ideal $I\subset S$, there exists a unique lexicographic ideal $I^{lex}$ such that $I$ and $I^{lex}$ have the same Hilbert function. The \textit{lexicographic ideal} $I^{lex}$ is constructed as follows. For each graded component $I_j$ of $I$, one defines $I_j^{lex}$ to be the $k-$vector space spanned by the unique degree $j$ initial lexsegment $L_j$ such that $|L_j|=\dim_k(I_j)$. Let $I^{lex}=\bigoplus\limits_{j} I_j^{lex}$. It is known that $I^{lex}$ constructed as before is indeed an ideal.

\begin{Example1}\rm
Let $I\subset k[x_1,x_2,x_3]$ be the lexsegment ideal determined by $u=x_1^2x_3$ and $v=x_2^3$, that is 
$$I=(x_1^2x_3,x_1x_2^2,x_1x_2x_3,x_1x_3^2,x_2^3).$$
We construct the lexicographic ideal $I^{lex}$. The ideal $I$ is generated in degree $3$. Since $\dim_k(I_3)=5$, we consider $I_3^{lex}$ to be the $k-$vector space spanned by the initial lexsegment set $L_3=\{x_1^3,x_1^2x_2,x_1^2x_3,x_1x_2^2,x_1x_2x_3\}$. The component of degree $4$ of the ideal $I$ is 
$$I_4=\shad(I)=\{x_1^3x_3,x_1^2x_2^2,x_1^2x_2x_3,x_1^2x_3^2,x_1x_2^3,x_1x_2^2x_3,x_1x_2x_3^2,x_2^4,x_1x_3^3,x_2^3x_3\},$$
therefore we consider $I_4^{lex}$ to be the $k-$vector space spanned by the initial lexsegment set $L_4=\{x_1^4,x_1^3x_2,x_1^3x_3,x_1^2x_2^2,x_1^2x_2x_3,x_1^2x_3^2,x_1x_2^3,x_1x_2^2x_3,x_1x_2x_3^2,x_1x_3^3\}$.
Since all the monomials belonging to $I_4^{lex}$ are divisible by the monomials from $I_3^{lex}$, we obtain $I^{lex}=(x_1^3,x_1^2x_2,x_1^2x_3,x_1x_2^2,x_1x_2x_3)$.
\end{Example1}

The notion of Gotzmann ideal is connected to the one of lexicographic ideal.

\begin{Definition1}\rm
A graded ideal $I\subset S$ generated in degree $d$ is called a \textit{Gotzmann ideal} if the number of generators of $\frak mI$ is the smallest possible, namely it is equal to the number of generators of $\frak mI^{lex}$. 
\end{Definition1}

We recall the Gotzmann persistence theorem \cite{G}, which describes the growth of the Hilbert function.

\begin{Theorem1}[Gotzmann's persistence, \cite{G}]
Let $I\subset S$ be a homogeneous ideal generated by elements of degree at most $d$. If $H(I,d+1)=H(I,d)^{\langle d\rangle}$, then $H(I,q+1)=H(I,q)^{\langle q\rangle}$, for all $q\geq d$.   
\end{Theorem1}

Therefore, by Gotzmann's persistence theorem \cite{G}, a graded ideal $I\subset S$ generated in degree $d$ is Gotzmann if and only if $I$ and $(I^{lex})_{\langle d\rangle}$ have the same Hilbert function.

\begin{Example1}\rm
$(1)$ Initial lexsegment ideals are Gotzmann ideals.

$(2)$ The ideal from the previous example is not a Gotzmann ideal. Indeed, we have to verify if $I$ and $(I^{lex})_{\langle 3\rangle}$ have the same Hilbert function. By the Gotzmann persistence theorem, we have to verify if $H(I,4)=H((I^{lex})_{\langle3\rangle},4)$. But $H(I,4)=\dim_k(I_4)=10$ and $H((I^{lex})_{\langle3\rangle},4)=|\shad(L_3)|=9$. Hence the ideal $I=(x_1^2x_3,x_1x_2^2,x_1x_2x_3,x_1x_3^2,x_2^3)$ is not a Gotzmann ideal. 
\end{Example1}



The recent work of A.H. Hoefel contributes to the characterization of Gotzmann ideals. He proved that the Gotzmann edge ideals are precisely the ones that arise from star graphs, see \cite{Hoe}. 

J. Herzog and T. Hibi generalized the notion of Gotzmann ideal as follows:

\begin{Definition1}\rm\cite{HH}
A graded ideal $I$ of $S$ is \textit{a Gotzmann ideal} if all ideals $I_{\langle j\rangle}$ are Gotzmann ideals.  
\end{Definition1}

In this sense, Gotzmann ideals are componentwise linear.

For Gotzmann ideals we have the following characterization in terms of (graded) Betti numbers \cite{HH}.

\begin{Theorem1}\cite{HH}\label{betti I and I^lex}
Let $I\subset S$ be a graded ideal. The following conditions are equivalent:
\begin{itemize}
	\item [(a)] $\beta_{ij}(S/I)=\beta_{ij}(S/I^{lex})$ for all $i,j$.
	\item [(b)] $\beta_{1j}(S/I)=\beta_{1j}(S/I^{lex})$ for all $j$.
	\item [(c)] $\beta_1(S/I)=\beta_1(S/I^{lex})$.
	\item [(d)] $I$ is a Gotzmann ideal.
\end{itemize}     
\end{Theorem1} 

Let $I$ be a Gotzmann monomial ideal generated in degree $d$. From the above results it follows that $I^{lex}$ is also generated in degree $d$ and $I$ has a linear resolution.

Using the Hilbert function, the Gotzmann ideals in $S$ which are generated by at most $n$ homogeneous polynomials are described in \cite{MH}.

A componentwise ideal $I$ is Gotzmann if and only if the Taylor resolution of $I$ is minimal, by Theorem \ref{minTay}.

\section{Simplicial complexes and Alexander duality}
Squarefree monomial ideals have interesting combinatorial properties due to their connections with simplicial complexes.

	\[
\]

\subsection{Basic notions}
We proceed at describing the general concepts which are important when dealing with simplicial complexes.
   
\begin{Definition1}\rm
A finite \textit{simplicial complex} consists of a finite set $V$ of vertices and a collection $\Delta$ of subsets of $V$ called \textit{faces} such that 
\begin{itemize}
	\item [(i)] If $v\in V$, then $\{v\}\in \Delta$;
	\item [(ii)] If $F\in\Delta$ and $G\subset F$, then $G\in\Delta$. 
\end{itemize}
\end{Definition1}

In general, we will consider the \textit{vertex set} $V=[n]=\{1,2,\ldots,n\}$. 

Let $\Delta$ be a simplicial complex on $[n]$ and $F$ a face of $\Delta$. The \textit{dimension} of $F$ is $\dim(F)=|F|-1$ and the \textit{dimension} of $\Delta$ is $\dim(\Delta)=\sup\{\dim(F):F\in\Delta\}$. 

The faces of dimension zero and $1$ are called \textit{vertices} and \textit{edges}, respectively.

The \textit{facets} of a simplicial complex are the maximal faces with respect to inclusion. It is customary to denote by $\mathcal F(\Delta)=\{F_1,\ldots,F_s\}$ the set of all the facets of a simplicial complex $\Delta$, and in this case we write $\Delta=\langle F_1,\ldots ,F_s\rangle$ and say that $\Delta$ is \textit{generated} by $F_1,\ldots,F_s$.

A simplicial complex is \textit{pure} if all its facets have the same dimension.

Let $\Delta$ be a simplicial complex on the vertex set $[n]$ of dimension $d-1$. Let $f_i$ denote the number of faces of dimension $i$ of $\Delta$. The sequence $f(\Delta)=(f_0,f_1,\ldots ,f_{d-1})$ is called the \textit{$f-$vector} of $\Delta$. By convention, we put $f_{-1}=1$. Then one can define the \textit{$h-$vector} of $\Delta$, $h(\Delta)=(h_0,h_1,\ldots,h_d,\ldots)$, by formula:
$$\sum\limits_{i=0}^d f_{i-1}(t-1)^{d-i}=\sum\limits_{i=0}^d h_it^{d-i}.$$
Equivalently,
$$\sum\limits_{i=0}^d f_{i-1}t^i(t-1)^{d-i}=\sum\limits_{i=0}^d h_it^{i}.$$

In particular, the $h-$vector has length at most $d$ and for $j=0,\ldots, d$ one has
$$h_j=\sum\limits_{i=0}^j(-1)^{j-i}{d-i\choose j-i}f_{i-1}.$$

The possible $f-$vectors of simplicial complexes have been determined by J. Kruskal and G. Katona in \cite{Ka}, \cite{Kr}. The $(f_0,\ldots,f_{d-1})\in \mathbb Z^d$ is the $f-$vector of a simplicial complex of dimension $d-1$ if and only if $0<f_{i+1}\leq f_i^{(i+1)}$, $0\leq i\leq d-2$.

\begin{Example1}\rm
We consider the simplicial complex $\Delta$: 
\newpage
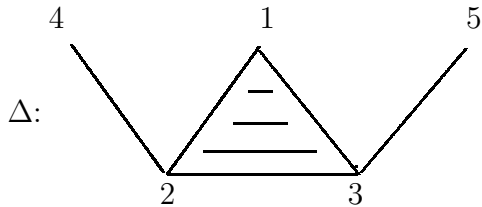
\begin{figure}[h]
\begin{center}
\unitlength 1mm 
\linethickness{0.8pt}
\ifx\plotpoint\undefined\newsavebox{\plotpoint}\fi 
\begin{picture}(70.75,27.5)(0,0)
\multiput(26,6)(.0336538462,.0467032967){364}{\line(0,1){.0467032967}}
\put(51,7){\line(1,0){.25}}
\put(26,6){\line(1,0){25.25}}
\put(30.75,9){\line(1,0){15}}
\put(34.75,12.75){\line(1,0){7.25}}
\put(36.75,17){\line(1,0){3.25}}
\put(38.25,25.5){$1$}
\put(25,2){$2$}
\put(50,2){$3$}
\put(10.25,25.5){$4$}
\put(65.75,25.5){$5$}
\put(4.75,13){$\Delta$:}
\multiput(38.25,22.5)(.0337150127,-.0419847328){393}{\line(0,-1){.0419847328}}
\multiput(25.5,6.25)(-.0336538462,.0467032967){364}{\line(0,1){.0467032967}}
\multiput(51.75,6.25)(.0337349398,.0397590361){415}{\line(0,1){.0397590361}}
\end{picture}
\caption{The geometric realization of a simplicial complex}
\label{figure1}
\end{center}
\end{figure}

In this case, the facets of $\Delta$ are:
$$\mathcal F(\Delta)=\{\{1,2,3\},\{2,4\},\{3,5\}\},$$
one of them, $\{1,2,3\}$, being of dimension $2$ and the other two of dimension $1$. Hence the simplicial complex is of dimension $\dim(\Delta)=2$. Moreover, $\Delta$ is not a pure simplicial complex.

The $f-$vector of $\Delta$ is $f(\Delta)=(5,5,1)$ and the $h-$vector is $h(\Delta)=(1,2,-2,0)$.
\end{Example1}

\begin{Definition1}\rm
Let $\Delta$ be a simplicial complex of dimension $d-1$. For each $0\leq i\leq d-1$, one can define the \textit{$i-$th skeleton} to be the simplicial complex
$$\Delta^{(i)}=\{F\in \Delta:\dim(F)\leq i\}.$$
\end{Definition1}

In particular, $\Delta^{(d-1)}=\Delta$.

Let $\Delta$ be a simplicial complex and $F$ a face of $\Delta$. One may consider the \textit{link} of $F$ in $\Delta$ to be 
$$\link_{\Delta}(F)=\{G\in \Delta:F\cup G\in\Delta,\ F\cap G=\emptyset\}.$$
In particular, $\link_{\Delta}(\emptyset)=\Delta$.

\begin{Definition1}\rm
A simplicial complex $\Delta$ on the vertex set $V$ is called \textit{disconnected} if there exist two non-empty sets $V_1$ and $V_2$, with $V=V_1\cup V_2$ and $V_1\cap V_2=\emptyset$, such that no face of $\Delta$ has vertices in both $V_1$ and $V_2$. A simplicial complex which is not disconnected is called \textit{connected}.
\end{Definition1}   

Obviously, a simplicial complex $\Delta$ is disconnected if and only if $\Delta^{(1)}$ is disconnected.

	\[
\]

\subsection{Ideals associated to a simplicial complex}

Let $S=k[x_1,\ldots,x_n]$ be the polynomial ring in $n$ variables over a field $k$. To any simplicial complex on $[n]$ one can associate a squarefree monomial ideal, as follows:

\begin{Definition1}\rm
Let $\Delta$ be a simplicial complex on the vertex set $[n]$. The \textit{Stanley--Reisner ideal} of $\Delta$ is 
$$I_{\Delta}=(x_F:F\notin \Delta),$$
where for $F\subseteq[n]$, $x_F=\prod\limits_{i\in F}x_i$.
\end{Definition1}
 
The minimal monomial generating set of the Stanley--Reisner ideal $I_{\Delta}$ consists of the monomials $x_F$, with $F$ a minimal nonface of $\Delta$.

The quotient ring $S/I_{\Delta}$, denoted by $k[\Delta]$, is called the \textit{Stanley--Reiner ring} of $\Delta$.

\begin{Example1}\rm
The Stanley--Reisner ideal of the simplicial complex from Figure \ref{figure1} is 
$$I_{\Delta}=(x_1x_4,x_1x_5,x_2x_5,x_3x_4,x_4x_5).$$
\end{Example1}

The standard primary decomposition of the Stanley--Reisner ideal is determined by the facets of the simplicial complex.

\begin{Proposition1}\cite{St}\label{St Reisner ideal decomposition}
Let $\Delta$ be the simplicial complex on the vertex set $[n]$ and $k$ be a field. Then
$$I_{\Delta}=\bigcap\limits_{F\in\mathcal F(\Delta)}P_{F^c},$$
where $P_{F^c}$ are the prime ideals generated by all the variables $x_i$ such that $i\notin F$. 
\end{Proposition1}

\begin{Corollary1}\cite{St}
If $\Delta$ is a simplicial complex of dimension $d-1$, then the Krull dimension of $k[\Delta]$ is
$$\dim(k[\Delta])=\dim(\Delta)+1=d.$$
\end{Corollary1}

\begin{Example1}\rm Consider the following simplicial complex:
	
\begin{figure}[h]
\begin{center}
\unitlength 1mm 
\linethickness{0.8pt}
\ifx\plotpoint\undefined\newsavebox{\plotpoint}\fi 
\begin{picture}(72.25,30)(0,0)
\multiput(12.75,28)(.0336826347,-.0636227545){334}{\line(0,-1){.0636227545}}
\multiput(24,6.75)(.0336676218,.0393982808){349}{\line(0,1){.0393982808}}
\put(24.25,6.75){\line(0,1){0}}
\put(24.25,6.5){\line(1,0){24.75}}
\multiput(49,6.5)(.0337389381,.046460177){452}{\line(0,1){.046460177}}
\multiput(64.25,27.5)(-.13255814,-.03372093){215}{\line(-1,0){.13255814}}
\multiput(35.75,20.25)(-.095588235,.033613445){238}{\line(-1,0){.095588235}}
\multiput(49,6.5)(-.05,.0333333){15}{\line(-1,0){.05}}
\multiput(48.25,7)(-.0337301587,.0357142857){378}{\line(0,1){.0357142857}}
\multiput(18.5,23.75)(.033602151,-.059139785){186}{\line(0,-1){.059139785}}
\multiput(24.5,22)(.033505155,-.06185567){97}{\line(0,-1){.06185567}}
\multiput(40.75,19.5)(.033613445,-.035714286){238}{\line(0,-1){.035714286}}
\multiput(46.25,20.25)(.033557047,-.035234899){149}{\line(0,-1){.035234899}}
\multiput(51.75,22)(.03370787,-.03932584){89}{\line(0,-1){.03932584}}
\multiput(57.25,24)(.03289474,-.03947368){38}{\line(0,-1){.03947368}}
\put(11.25,30){$1$}
\put(35,23.5){$2$}
\put(65.25,30){$3$}
\put(51.75,4.25){$4$}
\put(21,4.25){$5$}
\put(5,15){$\Delta$:}
\end{picture}
\end{center}
\caption{The simplicial complex $\Delta$}
\label{figure2}
\end{figure}
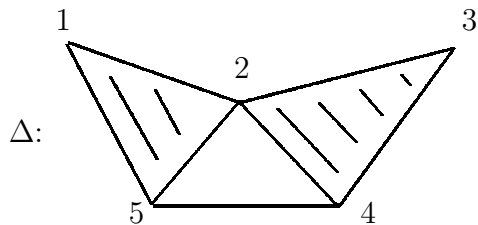

\noindent of dimension $\dim(\Delta)=2$, where the facets are $\mathcal F(\Delta)=\{\{1,2,5\},\{2,3,4\},\{4,5\}\}$. The Stanley--Reisner ideal is generated by the monomials $x_1x_3,x_1x_4,x_3x_5,x_2x_4x_5$. The standard decomposition of $I_{\Delta}$ is 
$$I_{\Delta}=(x_3,x_4)\cap(x_1,x_2,x_3)\cap(x_1,x_5).$$
\end{Example1}

\begin{Example1}\rm
For every integer $q\in[n]$, we denote by $I_{n,q}\subset S$ the squarefree monomial ideal generated by all the squarefree monomials of degree $q$. Then $I_{n,q}$ is the Stanley--Reisner ideal of $\Delta$, where $\Delta$ is generated by all the subsets of $[n]$ of cardinality $q-1$. Applying Proposition \ref{St Reisner ideal decomposition}, we obtain 
$$I_{n,q}=\bigcap\limits_{F\in\mathcal F(\Delta)}P_{F^c}=\bigcap\limits_{|F|=q-1}P_{F^c}.$$
\end{Example1}



For a simplicial complex $\Delta$, the multiplicity of the Stanley--Reisner ring $k[\Delta]$ can be expressed in terms of the $f-$vector.

\begin{Lemma1}\cite{BH}\label{multiplicity}
Let $\Delta$ be a simplicial complex of dimension $d-1$ and $k$ be a field. Then the multiplicity of the Stanley--Reisner ring of $\Delta$ is $e(k[\Delta])=f_{d-1}$. 
\end{Lemma1}
\begin{Definition1}\rm
A simplicial complex is called \textit{Cohen--Macaulay} if $k[\Delta]$ is a Cohen--Macaulay ring.
\end{Definition1}

The property of a simplicial complex of being Cohen--Macaulay depends on the characteristic of the base field. See \cite{BH} for some examples.

\begin{Proposition1}\cite{BH}
Any Cohen--Macaulay simplicial complex is pure.
\end{Proposition1}

An important class of Cohen--Macaulay simplicial complexes is given by the shellable simplicial complexes. We firstly recall their definition.

\begin{Definition1}\rm
A pure simplicial complex $\Delta$ is \textit{shellable} if its facets can be ordered $F_1,F_2,\ldots,F_m$ such that for all $1\leq j < i\leq m$, there exist some $v\in F_i\setminus F_j$ and some $s\in\{1,2,\ldots,i-1\}$ with $F_i\setminus F_s =\{v\}$.
\end{Definition1}

\begin{Theorem1}\cite{HeHi}
A pure shellable simplicial complex is Cohen–-Macaulay over an arbitrary field.
\end{Theorem1}

\begin{Proposition1}\cite{BH}\label{CM connected}
Let $\Delta$ be a simplicial complex.
\begin{itemize}
	\item [(a)] If $\dim(\Delta)=0$, then $\Delta$ is Cohen--Macaulay.
	\item [(b)] If $\Delta$ is disconnected, then $\depth(k[\Delta])=1$. In particular, a Cohen--Macaulay simplicial complex of positive dimension is connected.
	\item [(c)] Suppose $\dim(\Delta)=1$. The following conditions are equivalent:
\begin{itemize}
	\item [(i)] $\Delta$ is connected.
	\item [(ii)] $\Delta$ is shellable.
	\item [(iii)] $\Delta$ is Cohen--Macaulay.
\end{itemize}
\end{itemize}
\end{Proposition1}

An important Cohen--Macaulay criterion for simplicial complexes is due to G.A. Reisner \cite{Re}.

\begin{Theorem1}[Reisner, \cite{Re}]
\label{Reisner}
Let $\Delta$ be a simplicial complex on the vertex set $[n]$ of dimension $d-1$. Then $\Delta$ is Cohen-Macaulay over $k$ if and only if 
$\tilde{H}_i(\Delta;k)=0$ for all $i<d-1$ and the links of all vertices of $\Delta$ are Cohen-Macaulay.
\end{Theorem1}

For a squarefree monomial ideal $I\subset S$, one may use the following criterion due to T. Hibi \cite{H} to compute $\depth(S/I)$. 

\begin{Lemma1}\label{CM skeleton}\cite{H}
Let $\Delta$ be a simplicial complex of dimension $d-1$ and $\Delta^{(i)}$ be its $i-$th skeleton, for all $0\leq i\leq d-1$. Then:
\begin{enumerate}
	\item [(a)] If $\Delta$ is Cohen--Macaulay, then $\Delta^{(i)}$ is Cohen--Macaulay, for all $0\leq i\leq d-1$;
	\item [(b)] $\depth(k[\Delta])=\max\{i+1:k[\Delta^{(i)}]\mbox{ is Cohen--Macaulay }\}$.
\end{enumerate}
\end{Lemma1}

\begin{Example1}\rm
We consider the Stanley--Reisner ideal 
$$I_{\Delta}=(x_1x_3,x_1x_4,x_3x_5,x_2x_4x_5)$$
associated to the simplicial complex from Figure \ref{figure2}. The ring $k[\Delta]$ is not Cohen--Macaulay, since it is not pure. Moreover, $\dim(k[\Delta])=3$ and $\depth(k[\Delta])=2$, by Lemma \ref{CM skeleton}. 
\end{Example1}

Another important ideal which can be associated to a simplicial complex is the \textit{facet ideal}. R. S. Villarreal in \cite{Vi1} defined the notion of facet ideal of a graph, which is called \textit{edge ideal}. Later, S. Faridi \cite{Fa} generalized the study of facet ideal, namely to the case of simplicial complexes. 

\begin{Definition1}
Let $\Delta$ be a simplicial complex on the vertex set $[n]$. The \textit{facet ideal} of $\Delta$ is the squarefree monomial ideal
$$I(\Delta)=(x_F: F\in\mathcal F(\Delta)).$$
\end{Definition1}

Thus if $\Delta=\langle F_1,\ldots,F_s\rangle$, then the facet ideal of $\Delta$ is $I(\Delta)=(x_{F_1},\ldots,x_{F_s})$.

\begin{Example1}\rm
The facet ideal of the simplicial complex from Figure \ref{figure2} is
$$I(\Delta)=(x_1x_2x_5,x_2x_3x_4,x_4x_5)\subset k[x_1,\ldots,x_5].$$
\end{Example1}

	\[
\]


\subsection{Edge ideals}

Let $G$ be a finite simple graph on the vertex set $V(G)$, with the edge set $E(G)$. Without loosing the generality, we may assume that $V(G)=[n]$. A finite graph can be viewed as a $1-$dimensional simplicial complex. 

Let $S=k[x_1,\ldots,x_n]$ be the polynomial ring in $n$ variables over a field $k$. An important class of squarefree monomial ideals consists of edge ideals.

\begin{Definition1}\rm
The \textit{edge ideal} $I(G)$ associated to the graph $G=(V(G),E(G))$ is the squarefree monomial ideal in $S$ generated by all the monomials $x_ix_j$, with $\{i,j\}\in E(G)$.
\end{Definition1}

Let $G=(V(G),E(G))$ be a graph.
\begin{Definition1}\rm
A subset $A\subset V(G)$ is called a \textit{minimal vertex cover} of $G$ if every edge of $G$ is incident to one vertex in $A$ and there is no proper subset of $A$ with this property.   
\end{Definition1}

There is a strong relation between the minimal vertex covers of a graph and the minimal prime ideals of the edge ideal.

\begin{Proposition1}\cite{Vi}\label{min vertex cover}
Let $G$ be a graph on the vertex set $[n]$. If $P\subset S$ is the ideal generated by $A=\{x_{i_1},\ldots,x_{i_r}\}$, then $P$ is a minimal prime ideal of $I(G)$ if and only if $A$ is a minimal vertex cover of $G$.
\end{Proposition1}

\begin{Example1}\rm
We consider the graph $G$ on the vertex set $V(G)=\{1,\ldots,6\}$:

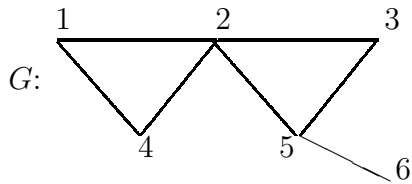
\begin{figure}[h]
\begin{center}
\unitlength 1mm 
\linethickness{0.8pt}
\ifx\plotpoint\undefined\newsavebox{\plotpoint}\fi 
\begin{picture}(70,25.75)(0,0)
\put(17.25,21.5){\line(1,0){21.25}}
\multiput(38.5,21.5)(-.0337171053,-.0419407895){304}{\line(0,-1){.0419407895}}
\multiput(28.25,8.75)(-.0336391437,.0382262997){327}{\line(0,1){.0382262997}}
\put(38.75,21.5){\line(1,0){21.25}}
\multiput(59.75,21.5)(-.0337171053,-.0419407895){304}{\line(0,-1){.0419407895}}
\multiput(49,8.75)(-.0336391437,.0382262997){327}{\line(0,1){.0382262997}}
\put(49.5,8.75){\line(2,-1){12}}
\put(17,23){$1$}
\put(38.25,23){$2$}
\put(60.75,23){$3$}
\put(28,6){$4$}
\put(46.75,6){$5$}
\put(62.25,3){$6$}
\put(10.75,15){$G$:}
\end{picture}

\end{center}
\caption{The graph G}
\label{}
\end{figure}
	
The sets $A_1=\{1,2,5\}$, $A_2=\{1,3,4,5\}$, $A_3=\{1,2,3,6\}$, $A_4=\{2,3,4,6\}$ and $A_5=\{2,4,5\}$ are all the minimal vertex covers of $G$.
\end{Example1}

	\[
\]
\subsection{Alexander duality} The Alexander duals of simplicial complexes represents a new tool in the study of algebraic and combinatorial
properties of squarefree monomial ideals. N. Terai and T. Hibi \cite{TH} used the Alexander duality in the study of Stanley--Reisner rings. Later, J.A. Eagon and V. Reiner \cite{ER} proved a fundamental result which relates homological data of the Stanley--Reisner ideal of a simplicial complex to combinatorial properties of its Alexander dual.

\begin{Definition1}\rm
Let $\Delta$ be a simplicial complex on the vertex set $[n]$. The simplicial complex
$$\Delta^\vee=\{[n]\setminus F: F\notin\Delta\}$$ 
is called the \textit{Alexander dual} of $\Delta$.
\end{Definition1}

One may note that $G$ is a facet of $\Delta^\vee$ if $G=[n]\setminus F$, for some minimal nonface $F$ of $\Delta$.

It is easy to check that $(\Delta^\vee)^\vee=\Delta$.

\begin{Example1}\rm
Let $\Delta$ be the simplicial complex on the vertex set $\{1,2,3,4,5\}$:

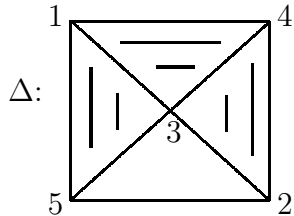
\begin{figure}[h]
\begin{center}
\unitlength 1mm 
\linethickness{0.8pt}
\ifx\plotpoint\undefined\newsavebox{\plotpoint}\fi 
\begin{picture}(40.75,30.5)(0,0)
\put(13.25,28.75){\line(0,-1){24}}
\put(13.25,4.75){\line(1,0){26.5}}
\put(39.75,4.75){\line(0,1){23.75}}
\put(39.75,28.5){\line(-1,0){26.25}}
\multiput(13.5,28.5)(.03728693182,-.03373579545){704}{\line(1,0){.03728693182}}
\multiput(13.5,5)(.03766140603,.03371592539){697}{\line(1,0){.03766140603}}
\put(16,22.5){\line(0,-1){10.75}}
\put(19.5,19){\line(0,-1){4.75}}
\put(37.5,23){\line(0,-1){12.25}}
\put(34,18.75){\line(0,-1){4.75}}
\put(20,25.75){\line(1,0){13.25}}
\put(24.75,22.5){\line(1,0){5}}
\put(10.25,28){$1$}
\put(40.75,28){$4$}
\put(40.75,3){$2$}
\put(10.25,3){$5$}
\put(26,12.5){$3$}
\put(5,18){$\Delta$:}
\end{picture}

\end{center}
\caption{The simplicial complex $\Delta$}
\label{figure3}
\end{figure}

The Alexander dual of $\Delta$ is $\Delta^\vee=\langle\{1,2,3\},\{1,4\},\{3,4,5\}\rangle$ and has the following geometric realization:
\begin{figure}[h]
\begin{center}
\unitlength 1mm 
\linethickness{0.8pt}
\ifx\plotpoint\undefined\newsavebox{\plotpoint}\fi 
\begin{picture}(74,24.25)(0,0)
\put(6.5,5){\line(1,0){65.75}}
\put(28.75,21.75){\line(1,0){18.5}}
\multiput(47.25,21.75)(.0496031746,-.0337301587){504}{\line(1,0){.0496031746}}
\multiput(6.75,5)(.0452716298,.0337022133){497}{\line(1,0){.0452716298}}
\multiput(29,21.5)(.0337301587,-.0634920635){252}{\line(0,-1){.0634920635}}
\multiput(37.5,5.5)(.0337370242,.0553633218){289}{\line(0,1){.0553633218}}
\put(14.75,7.75){\line(1,0){16}}
\put(22,12.25){\line(1,0){7}}
\put(46,7.5){\line(1,0){16}}
\put(47.75,12.75){\line(1,0){6.5}}
\put(28.25,23.75){$1$}
\put(48,23.75){$4$}
\put(5.5,1.5){$2$}
\put(37.25,1.5){$3$}
\put(74,1.5){$5$}
\put(4,15.75){$\Delta^\vee$:}
\end{picture}
\end{center}
\caption{The Alexander dual of $\Delta$}
\label{figure4}
\end{figure}
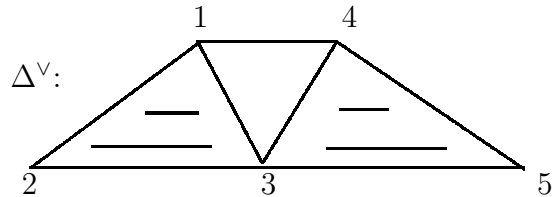
\end{Example1}

Next, we denote by $\Delta^c$ the simplicial complex with the facets $\mathcal F(\Delta^c)=\{[n]\setminus F: F\in\mathcal F(\Delta)\}$.

The Alexander duality allows us to connect the Stanley--Reisner ideal with the facet ideal.

\begin{Proposition1}\cite{HeHiZh} 
If $\Delta$ is a simplicial complex, then
$$I_{\Delta^\vee}=I(\Delta^c).$$
\end{Proposition1}

This last result, connected with the standard decomposition described before, gives us an effective method to determine the minimal monomial set of generators for the Stanley--Reisner ideal associated to the Alexander dual of a simplicial complex.

\begin{Corollary1}\cite{HeHi}
Let $\Delta$ be a simplicial complex and $I_{\Delta}=P_{F_1}\cap\ldots\cap P_{F_s}$ be the standard primary decomposition of $I_{\Delta}$, where each $F_j\subset [n]$. Then $G(I_{\Delta^\vee})=\{x_{F_1},\ldots, x_{F_s}\}$.
\end{Corollary1}

There are known some useful identities which describes the relations between subcomplex operations and operations with ideals:  
\begin{Lemma1}\cite{O1}\label{dual subcomplex}
Let $\Delta$ be a simplicial complex on the vertex set $[n]$ and $\Delta_1,\Delta_2$ subcomplexes of $\Delta$. Then:
\begin{itemize}
	\item[(a)] $I_{\Delta_1^{\vee}}\cap I_{\Delta_2^{\vee}}=I_{(\Delta_1\cap\Delta_2)^{\vee}}$;
	\item[(b)] $\Delta=\Delta_1\cup\Delta_2$ if and only if $I_{\Delta^{\vee}}=I_{\Delta_1^{\vee}}+I_{\Delta_2^{\vee}}$;
	\item[(c)] $\Delta=\Delta_1\cap\Delta_2$ if and only if $I_{\Delta}=I_{\Delta_1}+I_{\Delta_2}$;
	\item[(d)] $\Delta=\Delta_1\cap\Delta_2$ if and only if $\Delta^{\vee}=\Delta_1^{\vee}\cup\Delta_2^{\vee}$.
\end{itemize}
\end{Lemma1}

An application of the Alexander duality is given by the Eagon--Reiner theorem \cite{ER} which provides a powerful tool in the study of Cohen--Macaulay simplicial complexes.

\begin{Theorem1}[Eagon--Reiner, \cite{ER}]\label{Eagon--Reiner}
Let $\Delta$ be a simplicial complex on $[n]$ and let $k$ be a field. Then the Stanley--Reisner ideal $I_{\Delta}\subset k[x_1,\cdots, x_n]$ has a linear resolution if and only if $k[\Delta^\vee]$ is Cohen--Macaulay. More precisely, $I_{\Delta}$ has a $q-$linear resolution if and only if $k[\Delta^\vee]$ is Cohen--Macaulay of dimension $n-q$.
\end{Theorem1}

Another interesting result, relating the projective dimension of a Stanley--Reisner ideal to the regularity of the Stanley--Reisner ring of the Alexander dual is due to N. Terai \cite{T2}.

\begin{Proposition1}\cite{T2}
Let $\Delta$ be a simplicial complex. Then
$$\projdim(I_\Delta)=\reg(k[\Delta^\vee]).$$
\end{Proposition1}

The concept of sequentially Cohen--Macaulay simplicial complexes was introduced by R. Stanley \cite{St}.

\begin{Definition1}\rm
A simplicial complex $\Delta$ is \textit{sequentially Cohen--Macaulay} if the Stanley--Reisner ring $k[\Delta]$ is sequentially Cohen--Macaulay.
\end{Definition1}

A pure simplicial complex is sequentially Cohen--Macaulay if and only if it is Cohen--Macaulay, \cite{St}.

A.M. Duval \cite{Du} proved that if the Stanley--Reisner ideal of a simplicial complex is sequentially Cohen--Macaulay, then the simplicial complex has Cohen--Macaulay subcomplexes.

A generalization of the Eagon--Reiner theorem (Theorem \ref{Eagon--Reiner}) is due to J. Herzog, V. Reiner and V. Welker.

\begin{Theorem1}\cite{HRW}\label{comp lin-seq}
Let $\Delta$ be a simplicial complex on $[n]$. Then $I_{\Delta}\subset S$ is componentwise linear if and only if $\Delta^\vee$ is sequentially Cohen--Macaulay.
\end{Theorem1}



\section{Arithmetical rank of squarefree monomial ideals}

Let $S=k[x_1,\ldots,x_n]$ be the polynomial ring in $n$ variables over a field $k$. Let $I\subset S$ be a homogeneous ideal and 
$\sqrt{I}$ its radical.
\begin{Definition1}\rm The \textit{arithmetical rank} of $I$ is defined as
\[
\ara(I)=\min\{r\in \NN\colon \text{ there exist }a_1,\ldots,a_r\in I \text{ such that }\sqrt{I}=\sqrt{(a_1,\ldots,a_r)}\}.
\]
\end{Definition1}

For the study of the arithmetical rank of an ideal we refer to \cite{SV}, \cite{ScV} and \cite{L}. The problem of determining the minimum number of equations required to generate a monomial ideal up to radical was firstly studied by P. Schenzel and W. Vogel \cite{SV}, Th. Schmitt and W. Vogel \cite{ScV} and G. Lyubeznik \cite{L}.


Geometrically, $\ara(I)$ is the smallest number of hypersurfaces whose intersection is set-theoretically equal to the algebraic set defined 
by $I$, if $k$ is algebraically closed.

It is known that, for an ideal $I\subset S$, we always have $\ara(I)\geq\height(I)$, by the Krull's principal ideal theorem \cite{Vi}. 

\begin{Definition1}\rm
 An ideal $I\subset S$ which satisfies $\ara(I)=\height(I)$ is called \textit{set-theoretic complete intersection}.
\end{Definition1}

An upper bound for the arithmetical rank of an ideal is given by the Eisenbud--Evans theorem.

\begin{Theorem1}[Eisenbud--Evans, \cite{Ei-Ev}] Let $R$ be a Noetherian ring of dimension $n$, and suppose that $R$ is a polynomial ring of the form $R=S[x]$, for some ring $S$. Let $I\subset R$ be an ideal. Then there exists $n$ elements $g_1,\ldots,g_n\in I$ such that $\sqrt{I}=\sqrt{(g_1,\ldots,g_n)}$.
\end{Theorem1}

Hence the arithmetical rank of an ideal $I\subset S=k[x_1,\ldots, x_n]$ is less than or equal to $n$, where $n$ is the number of variables. Moreover, if $I$ is homogeneous with the property that $\sqrt I\subsetneq\frak m=(x_1,\ldots,x_n)$, then $\ara(I)\leq n-1$.  



For a squarefree monomial ideal $I\subset S$ an upper bound of $\ara(I)$ is given by  H.G. Gr\"abe \cite{Gr}. 

\begin{Theorem1}\cite{Gr}
Let $I$ be a squarefree monomial ideal. Then
$$\ara(I)\leq n-\indeg(I)+1,$$
where $\indeg(I)$ is the \textit{initial degree} of $I,$ that is, $\indeg(I)=\min\{q\colon I_q\neq 0\}.$
\end{Theorem1}

Let $\cd(I)=\max\{i\in \ZZ\colon H_I^i(S)\neq 0\},$ where $H_I^i(S)$ denotes the $i$-th local cohomology module of $S$ with support at $V(I)$. The number $\cd(I)$ is called the \textit{cohomological dimension} of $I.$ By expressing the local cohomology modules in terms of Cech complex, one can see that $\ara(I)$ is bounded below by $\cd(I)$.

In his paper, G. Lyubeznik \cite{L} obtained equality between cohomological dimension and projective dimension.   

\begin{Theorem1}\cite{L} 
Let $I\subset S$ be a squarefree monomial ideal. Then $$\cd(I)=\projdim_S(S/I).$$ 
\end{Theorem1}

As a consequence, we have:
\begin{Corollary1}\cite{L}
Let $I\subset S$ be a squarefree monomial ideal. Then:
\[
\projdim_S(S/I)=\cd (I)\leq \ara(I)\leq n-\indeg(I)+1.
\] 
In particular, if $I$ is set--theoretic complete intersection, then $I$ is Cohen--Macaulay.
\end{Corollary1}

An example of ideal for which $\projdim_S(S/I)<\ara(I)$ is due to Z. Yan \cite{Ya}:

\begin{Example1}\rm\cite{Ya}
Let $I=(uvw,uvy,vwx,uwz,uxy,uxz,vxz,vyz,wxy,wyz)$ be the ideal in $S=k[u,v,w,x,y,z]$. Then $I$ is the Stanley-Reisner ideal of a
triangulation of $\mathrm{P}^{2}(\mathrm{S})$ with six vertices. Using the \'etale cohomology, Z. Yan proved that $\ara(I)=4$. Moreover, if $\chara(k)\neq 2$, then $\projdim(S/I)=3$.
\end{Example1}

It naturally arises the problem of determining conditions under which we have $\projdim_S(S/I)=\ara(I)$. 

There are many instances when the equality $\projdim_S(S/I)=\ara(I)$ holds. For example, if $\Delta$ is a disconnected simplicial complex, then the Stanley--Reisner ideal of $\Delta$ satisfies this equality.

Many authors such as M. Barile, N. Terai, K. Yoshida, M. Morales, K. Kimura, M. Kummini were involved in finding classes of squarefree monomial ideals $I\subset S$ which satisfy the equality $\projdim_S(S/I)=\ara(I)$, see for instance \cite{B2}, \cite{BT}, \cite{BT1}, \cite{Ku}.

We briefly recall some well-known cases.

W. Vogel, Th. Schmitt and P. Schenzel  \cite{SV}, \cite{ScV} studied the case of squarefree monomial ideals $I$ with $|\Ass_S(S/I)|=\indeg(I)$ 

\begin{Theorem1}\cite{SV}, \cite{ScV}
If $I\subset S$ is a squarefree monomial ideal with the property that $|\Ass_S(S/I)|=\indeg(I)$, then $\projdim_S(S/I)=\ara(I)$.
\end{Theorem1}

K. Kimura, N. Terai and K. Yoshida extended this result:
 
\begin{Theorem1}\cite{KTY}
Let $I$ be a squarefree monomial ideal of the polynomial ring $S=k[x_1,\ldots,x_n]$. If one of the following holds: 
\begin{enumerate}
	\item [(a)] $\mu(I)\leq \projdim(S/I)+1$
	\item [(b)] $|\Ass_S(S/I)|=\reg(I)$
	\item [(c)] $|\Ass_S(S/I)|=\indeg(I)+1,$
\end{enumerate}
then $\projdim_S(S/I)=\ara(I)$, where $\mu(I)$ is the minimal number of monomial generators of $I$. 
\end{Theorem1}

Other cases when the equality $\projdim_S(S/I)=\ara(I)$ holds are given in \cite{KTY1}.

\begin{Proposition1}\cite{KTY1}
If $I$ is a squarefree monomial ideal which has one of the properties:
\begin{itemize}
	\item [(a)] $\mu(I)\leq \height(I)+2$
	\item [(b)] $I$ is a squarefree monomial ideal with a linear resolution and $|\Ass_S(S/I)|=\indeg(I)+2$,
\end{itemize}
then the equality $\projdim_S(S/I)=\ara(I)$ holds.
\end{Proposition1}

M. Morales in \cite{M} studied the case of squarefree monomial ideals generated in degree $2$ with a linear resolution.

\begin{Theorem1}\cite{M}
Let $I$ be a squarefree monomial ideal which has a $2-$linear resolution. Then $\projdim_S(S/I)=\ara(I)$.
\end{Theorem1}

A useful tool which allows us to compute the arithmetical rank of a squarefree monomial ideal is the Schmitt--Vogel Lemma (see for instance \cite{SV}).

\begin{Lemma1}[Schmitt--Vogel, \cite{SV}]\label{svlemma}
Let $I\subset S$ be a squarefree monomial and $A_1,\ldots, A_r$ be some subsets of the set of monomials of $I.$ Suppose that the 
following conditions hold:
\begin{itemize}
	\item [(SV1)] $|A_1|=1$ and $A_i$ is a  finite set for any $2\leq i\leq r$;
	\item [(SV2)] The union of all the sets $A_i, i=\overline{1,r},$ contains the set of the minimal monomial generators of $I;$
	\item [(SV3)] For any $i\geq 2$ and for any two different monomials $m_1,m_2\in A_i$ there exist $j<i$ and a monomial 
	$m^{\prime}\in A_j$ such that $m^{\prime}|m_1m_2.$
\end{itemize}
Let  $g_i=\sum_{m_i\in A_i}m_i$ for $1\leq i\leq r.$ Then $\sqrt{(g_1,\ldots,g_r)}=I.$ In particular, $\ara(I)\leq r.$  
\end{Lemma1}

\begin{Example1}\rm
Let $I=(x_1x_2x_3,x_1x_4,x_1x_5,x_2x_4,x_2x_5,x_3x_4x_5)$ be a squarefree monomial ideal of $S=k[x_1,\ldots,x_5]$. The height of the ideal is $\height(I)=3$, since the standard primary decomposition is:
$$I=(x_1,x_2,x_3)\cap(x_1,x_2,x_4)\cap(x_1,x_2,x_5)\cap(x_1,x_4,x_5)\cap(x_2,x_4,x_5)\cap(x_3,x_4,x_5).$$

One has that $3=\height(I)\leq \ara(I)\leq 4$, by the Eisenbud--Evans theorem.

We construct the sets
	
$$ A_1=\{x_1x_5\}$$
$$A_2= \{x_1x_4,x_2x_5\}$$
$$A_3=\{x_1x_2x_3,x_2x_4,x_3x_4x_5\}.$$

It is easy to see that the sets satisfy the conditions of the Schmitt--Vogel lemma, hence $\ara(I)\leq 3$.

We obtain that $\ara(I)=3$ and the ideal is set--theoretic complete intersection.
\end{Example1}

\chapter{Classes of sequentially Cohen--Macaulay squarefree monomial ideals}

In this chapter we aim at determining the minimal primary decomposition of completely squarefree lexsegment ideals. As an application, we will characterize the sequentially Cohen--Macaulay completely squarefree lexsegment ideals. The results from this chapter are contained in \cite{OO}. 

\section[Primary decomposition]{Primary decomposition for completely squarefree lexsegment ideals}

The goal of this section is to determine the minimal primary decomposition of a completely squarefree lexsegment ideal. 
Initial squarefree lexsegment ideals are an example of completely squarefree lexsegment ideals. It was proved \cite{B} that any completely squarefree lexsegment ideal $I=(L(u,v))$ can be written as the intersection of the initial squarefree lexsegment ideal $(L^i(v))$ with the final squarefree lexsegment ideal $(L^f(u))$. Thus we will firstly determine the minimal primary decomposition for initial and final squarefree lexsegment ideals. 

Let $I=(L^i(v))\subset k[x_1,\ldots,x_n]$ be an initial squarefree lexsegment ideal, where $v=x_{j_1}\cdots x_{j_q}$. We may assume that $j_1\geq 2$. Otherwise, if $j_1=1$, then $I=(x_1)\cap (L(x_2\cdots x_q,v/x_1))$ and the problem reduces to compute the minimal primary decomposition of an initial squarefree lexsegment ideal in a fewer number of variables. 

\begin{Theorem1}\label{prim dec initial}
Let $I\subset k[x_1,\ldots,x_n]$ be the initial squarefree lexsegment ideal generated in degree $q$, determined by the monomial $v=x_{j_1}\cdots x_{j_q}$, with $2\leq j_1<\ldots<j_q\leq n$. Consider the sets $A_t=[j_t]\setminus\{j_1,\ldots,j_{t-1}\}$, for $1\leq t\leq q$. Then $I$ has the minimal primary decomposition of the form:
$$I=\left(\bigcap\limits_{t=1}^{q}(x_i\ :\ i\in A_t)\right)\cap\left(\bigcap\limits_\twoline{F\subset[n],\ |F|=q-1}{F\cap A_t\neq\emptyset,\ \forall t}P_{F^c}\right).$$  
\end{Theorem1}
  
\begin{proof}
Let $\Delta$ be the simplicial complex on the vertex set $[n]$ such that $I=I_{\Delta}$. We show that the facets of $\Delta$ are exactly the sets $[n]\setminus A_t$, $1\leq t\leq q$, together with the sets $F\subset [n]$ with $|F|=q-1$ and such that $F\cap A_t\neq\emptyset$ for all $t$. By applying Proposition \ref{St Reisner ideal decomposition} we get the desired formula.

In the first place, we observe that all the sets $G\subset [n]$ with $|G|=q-1$ are faces of $\Delta$ since $x_G\notin I$. Therefore, the facets of $\Delta$ have the cardinality at least $q-1$. Secondly, we note that $G$ is a facet of $\Delta$ if and only if $x_G\notin I$ and $x_{G\cup\{i\}}\in I$, for all $i\in [n]\setminus G$.

Let $G$ be a subset of the set $[n]$ with $|G|=q-1$. Then $G$ is a facet of $\Delta$ if and only if $x_{G\cup\{i\}}\in I$, for all $i\in[n]\setminus G$, which is equivalent to $x_Gx_{\max([n]\setminus G)}\geq_{lex} v$. We show that this last condition is equivalent to $G\cap A_t\neq\emptyset$, for all $t$. It is obvious that if $G\cap A_t=\emptyset$, for some $1\leq t\leq q$, then $G\subset \{j_1,\ldots ,j_{t-1}\}\cup\{j_t+1,\ldots,n\}$. Hence $x_Gx_{\max([n]\setminus G)}<_{lex} v$, which is a contradiction. Conversely, let $G\cap A_t\neq\emptyset$, for all $1\leq t\leq q$ and assume that $x_Gx_{\max([n]\setminus G)}<_{lex} v$. Then either $x_G\leq_{lex}v/x_{\max(v)}=x_{j_1}\cdots x_{j_{q-1}}$ or $x_{G\setminus\max(G)}x_{\max([n]\setminus G)}\leq_{lex} x_{j_1}\cdots x_{j_{q-1}}$. In each case, we get a contradiction with the condition $G\cap A_t\neq\emptyset$, for all $1\leq t\leq q$.

Next, we look at the facets of cardinality greater than or equal to $q$. Let $G\subset[n]$ be such a facet and $x_G=x_{l_1}\cdots x_{l_r}$, with $r\geq q$. It is clear that if $l_1<j_1$, then $x_G\in I$, thus $l_1\geq j_1$. The only facet with $l_1>j_1$ is $G=\{j_1+1,\ldots, n\}=[n]\setminus A_1$. Assume now that $l_1=j_1$. Then we obtain $l_2\geq j_2$. The only facet with $l_2>j_2$ is $G=\{j_1,j_2+1,\ldots, n\}=[n]\setminus A_2$. By applying this argument step by step, it follows that the facets of cardinality greater than or equal to $q$ are exactly $[n]\setminus A_t$, with $t=1,\ldots, q$.
\end{proof}

\begin{Example1}\rm
Let $I=(L^i(x_2x_4x_5))\subset k[x_1,\ldots ,x_6]$ be an initial squarefree lexsegment ideal. We consider the sets $A_1=\{1,2\}$, $A_2=\{1,3,4\}$ and $A_3=\{1,3,5\}$. Then the minimal primary decomposition of $I$ is of the form:
$$I=(x_1,x_2)\cap (x_1,x_3,x_4)\cap (x_1,x_3,x_5)\cap(x_3,x_4,x_5,x_6)\cap(x_2,x_4,x_5,x_6)\cap(x_2,x_3,x_5,x_6)\cap$$
$$\cap(x_2,x_3,x_4,x_6)\cap(x_2,x_3,x_4,x_5)\cap(x_1,x_4,x_5,x_6).$$
\end{Example1}

\begin{Corollary1}\label{dim}
Let $I=(L^i(v))\subset S$ be an initial squarefree lexsegment ideal, where $v=x_{j_1}\cdots x_{j_q}$, and $\Delta$ be the simplicial complex on the vertex set $[n]$ such that $I=I_{\Delta}$. Then:
\begin{itemize}
	\item [(a)] The dimension of the Stanley--Reisner ring of $\Delta$ is $n-j_1$.
	\item [(b)] The depth of the Stanley--Reisner ring of $\Delta$ is $q-1$.
	\item [(c)] The simplicial complex $\Delta$ is pure if and only if $v=x_{n-q+1}\cdots x_n$. Moreover, $\Delta$ is a pure simplicial complex if and only if $\Delta$ is Cohen--Macaulay.
\end{itemize}
\end{Corollary1}

\begin{proof}
(a) The height of the Stanley--Reisner ideal associated to $\Delta$ is 
$$\height I=\min\{j_1,j_2-1,\ldots,j_q-(q-1),n-q+1\}=j_1.$$
 Hence $\dim k[\Delta]=\dim(S/I)=n-j_1$.  

(b) The claim is obvious if $I= I_{n,q}$. We assume that $I\neq I_{n,q}$. Since $I$ is generated by squarefree monomials of degree $q$, we have that any subset $F$ of $[n]$, with $|F|\leq q-1$, is a face of $\Delta$. This implies that $I_{\Delta^{(q-2)}}$ is Cohen--Macaulay, because it is generated by all the monomials of degree $q$. By Lemma \ref{CM skeleton} we have that $\depth (S/I)\geq q-1$.  

The facets of $\Delta$ are of cardinality $n-j_1,n-j_2+1,\ldots,n-j_q+q-1,q-1$, with $n-j_1>q-1$. It results that $\Delta^{(q-1)}$ is not a pure simplicial complex, therefore it is not Cohen--Macaulay. By Lemma \ref{CM skeleton}, we obtain $\depth (S/I)\leq q-1$. 

(c) It is known that any Cohen--Macaulay simplicial complex is pure. 

The simplicial complex $\Delta$ is pure if all its facets have the same dimension. Using Proposition \ref{St Reisner ideal decomposition}, we obtain that $n-j_1=n-j_2+1=\ldots=n-j_q+q-1=q-1$. Therefore the monomial $v$ must have the form $v=x_{n-q+1}\cdots x_n$. 

Since the ideal generated by all the squarefree monomials of degree $q$ is Cohen--Macaulay, it results that $\Delta$ is a pure simplicial complex if and only if $\Delta$ is Cohen--Macaulay.
\end{proof}

As a consequence of the Proposition \ref{min vertex cover}, we obtain: 

\begin{Corollary1}\label{min vertex cover initial}
Let $G$ be the graph with the edge ideal $I(G)$, where $I(G)$ is the initial squarefree lexsegment ideal generated in degree $2$, determined by the monomial $v=x_{j_1}x_{j_2}$. Then 
$$I(G)=(x_i\ :\ i\in A_1)\cap(x_i\ :\ i\in A_2)\cap\left(\bigcap\limits_{i=1}^{j_1-1}P_{[n]\setminus\{i\}}\right),$$
where $A_1=[j_1],\ A_2=[j_2]\setminus\{j_1\}$. The sets $A_1,\ A_2$ and $[n]\setminus\{i\}$, for $i=1,\ldots ,j_1-1$, are the minimal vertex covers of $G$.
\end{Corollary1}

\begin{Corollary1}\label{multiplicity initial}
Let $I\subset S$ be the initial squarefree lexsegment ideal generated in degree $q$, determined by the monomial $v=x_{j_1}\cdots x_{j_q}$, with $2\leq j_1<\ldots<j_q\leq n$, $v\neq x_{n-q+1}\cdots x_n$, and $\Delta$ the simplicial complex with the Stanley--Reisner ideal $I$. If $s$ is the smallest integer such that $j_i=j_1+(i-1)$, for all $1\leq i< s$ and $j_s\geq j_1+s$, then the multiplicity of $k[\Delta]$ is $e(k[\Delta])=s-1$.
\end{Corollary1}

\begin{proof}
By Lemma \ref{multiplicity} and Corollary \ref{dim} (a), the multiplicity of $k[\Delta]$ is $e(k[\Delta])=f_{n-j_1-1}$. Therefore, we need to determine the number of facets of $\Delta$ of cardinality $n-j_1$. This is equivalent, by Theorem \ref{prim dec initial}, to determine the number of minimal prime ideals which contains $I$, with $\height(p)=j_1$ and $j_1\neq n-q+1$. In the notations of Theorem \ref{prim dec initial}, if $s$ is the smallest integer such that $j_i=j_1+(i-1)$, for all $1\leq i< s$ and $j_s\geq j_1+s$, then $|A_i|=j_i-(i-1)=j_1$, for all $1\leq i\leq s-1$ and $|A_i|=j_i-(i-1)>j_1$, for $i\geq s$. Therefore, $e(k[\Delta])=s-1$.   
\end{proof}


Next, we discuss the case of a final squarefree lexsegment ideal $(L^f(u))$, where $u$ is a monomial in $S$. Note that we may reduce to the hypothesis $x_1\mid u$. Indeed, otherwise, $x_1,\ldots,x_{\min(u)-1}$ are regular on $S/I$, hence they do not belong to any associated prime of $I$. Therefore, computing the primary decomposition of $I$ in $S$ is equivalent to compute the primary decomposition of $I\cap k[x_{\min(u)},\ldots,x_n]$.  

\begin{Lemma1}\label{deg-dual}
Let $I=(L^f(u))\subset S$ be the final squarefree lexsegment ideal determined by the monomial $u=x_1x_{i_2}\cdots x_{i_q}$, $2\leq i_2<\ldots<i_q\leq n$. We denote by $\Delta$ the simplicial complex with $I=I_{\Delta}$. Then the Stanley--Reisner ideal of the Alexander dual of $\Delta$, $I^\vee$, is generated in degree $n-q$ and $n-q+1$.
\end{Lemma1}

\begin{proof}
Since $u=x_1x_{i_2}\cdots x_{i_q}$, one may easily check that $\shad(I)=\Mon^s_{q+1}(S)$. Moreover, $\shad^s(I)=\Mon^s_{q+s}(S)$, for all $s\geq 1$. This implies that all the squarefree monomials of degree greater than or equal to $q+1$ belong to $I$. Hence, $I^{\vee}$ is generated in degree greater than or equal to $n-q$. 

On the other hand, all the monomials $x_G$, with $|G|\leq q-1$ do not belong to $I$, thus all the monomials $x_{[n]\setminus G}$ of degree greater than or equal to $n-q+1$ belong to $I^\vee$.

Therefore, $I^\vee$ is generated in degree $n-q$ and $n-q+1$.     
\end{proof}

The following result gives us the minimal primary decomposition of a final squarefree lexsegment ideal. 

\begin{Theorem1}\label{prim dec final}
Let $I\subset S$, $I\neq I_{n,q}$, be the final squarefree lexsegment ideal generated in degree $q$, determined by the monomial $u=x_1x_{i_2}\cdots x_{i_q}$, $2\leq i_2<\ldots<i_q\leq n$. Let $F=\{i_2,\ldots,i_q\}$ and $x_F=\prod\limits_{i\in F}x_i$. Then $I$ has the minimal primary decomposition of the form:
$$I=\left(\bigcap\limits_\twoline{G\subset[n],\ |G|=n-q+1}{x_G\geq_{lex}x_{F^c}}P_G\right)\cap\left(\bigcap\limits_\twoline{G\subset[n]\setminus\{1\},\ |G|=n-q+1}{x_{G\setminus\min(G)}\geq_{lex}x_{F^c\setminus\{1\}}}P_G\right)\cap\left(\bigcap\limits_\twoline{G\subset[n],\ |G|=n-q}{x_{F^c\setminus\{1\}}>_{lex}x_G}P_G\right).$$
\end{Theorem1}

\begin{proof}
Let $\Delta$ be the simplicial complex on the vertex set $[n]$ such that $I=I_{\Delta}$. By Lemma \ref{deg-dual}, we know that all the facets of $\Delta$ have cardinality $q-1$ or $q$. A facet $H$ of $\Delta$ is characterized by the condition $x_H\notin I$ and $x_{H\cup\{i\}}\in I$, for all $i\in [n]\setminus H$. 

Let $H$ be a subset of $[n]$ of cardinality $q$. Note that, if $|H|=q$, then $x_Hx_i\in \Mon_{q+1}^s(S)=\shad(I)\subset I$, for all $i\in [n]\setminus H$, hence $H$ is a facet of $\Delta$ if and only if $x_H\notin I$. We have $x_H\notin I$ if and only if $x_H>_{lex} u=x_1x_F$, that is $x_{H^c}<_{lex}x_{F^c\setminus\{1\}}$. Therefore, by denoting $G=H^c$, we get the last family of minimal prime ideals in the formula of the theorem.

Now, we look at the facets of cardinality $q-1$. Let $H$ be a subset of $[n]$, with $|H|=q-1$. Then $x_H\notin I$, thus it remains to characterize the sets $H$ with $x_{H\cup\{i\}}\in I$, for all $i\in [n]\setminus H$. This is equivalent to $x_Hx_{\min([n]\setminus H)}\in I$, that is $u=x_1x_F\geq_{lex}x_Hx_{\min([n]\setminus H)}$. We distinguish two cases:

\textit{Case 1:} If $\min([n]\setminus H)=1$, then $x_1x_F\geq_{lex}x_Hx_{\min([n]\setminus H)}$ is equivalent to $x_F\geq _{lex} x_H$, that is $x_{H^c}\geq_{lex} x_{F^c}$. Therefore, in this case, we get the first family of minimal prime ideals of $I$.

\textit{Case 2:} If $\min([n]\setminus H)>1$, then $1\in H$. By taking the complements in $x_1x_F\geq_{lex}x_Hx_{\min([n]\setminus H)}$, we obtain $x_{([n]\setminus H)\setminus\{\min([n]\setminus H)\}}\geq _{lex} x_{F^c\setminus\{1\}}$. By setting $G=[n]\setminus H$, we get the second family of minimal prime ideals.
\end{proof}

\begin{Example1}\rm
Let $I=(L^f(x_1x_2x_5))\subset k[x_1,\ldots ,x_6]$ be a final squarefree lexsegment ideal. In this case, the monomial $x_{F^c}$ is $x_1x_3x_4x_6$. The minimal primary decomposition of $I$ is of the form:
$$I=(x_1,x_2,x_3,x_4)\cap (x_1,x_2,x_3,x_5)\cap (x_1,x_2,x_3,x_6)\cap (x_1,x_2,x_4,x_5)\cap(x_1,x_2,x_4,x_6)\cap $$
$$\cap(x_1,x_2,x_5,x_6)\cap(x_1,x_3,x_4,x_5)\cap (x_1,x_3,x_4,x_6)\cap(x_3,x_5,x_6)\cap(x_4,x_5,x_6)\cap$$
$$\cap(x_2,x_3,x_4,x_5)\cap(x_2,x_3,x_4,x_6).$$
\end{Example1}

The minimal primary decomposition allows us to compute the Krull dimension of the quotient ring $S/I$.
 
\begin{Corollary1}\label{dim final}
Let $I=(L^f(u))\subset S$ be the final squarefree lexsegment ideal determined by the monomial $u$, where $u=x_1x_{i_2}\cdots x_{i_q}$, $2\leq i_2<\ldots<i_q\leq n$, $u\neq x_1x_2\cdots x_q$. Let $\Delta$ be the simplicial complex with the Stanley--Reisner ideal $I$. Then:
\begin{itemize}
	\item [(a)] $\dim(S/I)=q$.
	\item [(b)] $\depth(S/I)=q-1$.
	\item [(c)] The multiplicity of $k[\Delta]$ is $e(k[\Delta])=|\{x_G: x_{F^c\setminus\{1\}}>_{lex}x_G, |G|=n-q\}|$.
\end{itemize}
\end{Corollary1}

\begin{proof}
(b) As in the first part of the proof of Corollary \ref{dim} (b), we have that $\depth (S/I)\geq q-1$.

We prove that $\Delta ^{(q-1)}$ is not pure. Indeed, because $u\neq x_1\cdots x_q$, we have that $\{x_1,\ldots,x_q\}$ is a facet of $\Delta^{(q-1)}$. We consider $\tau=\{x_{n-q+2},\ldots,x_n\}$. Then $|\tau|=q-1$ and $\tau$ is a maximal face of $\Delta^{(q-1)}$ because all the monomials $x_ix_{n-q+2}\cdots x_n\in I$, for all $1\leq i\leq n-q+1$. Therefore $\Delta ^{(q-1)}$ is not pure. This implies that $k[\Delta ^{(q-1)}]$ is not Cohen--Macaulay and we get $\depth(S/I)\leq q-1$.

(c) The multiplicity of $k[\Delta]$ is $e(k[\Delta])=f_{q-1}$, by (a) and Lemma \ref{multiplicity}. The number of facets of $\Delta$ of cardinality $q$, that is $f_{q-1}$, is the same with the number of minimal prime ideals of $I$ of $\height(p)=n-q$. Thus, by Theorem \ref{prim dec final}, the desired conclusion follows.
\end{proof}

\begin{Corollary1}\label{min vertex cover final}
Let $G$ be the graph with the edge ideal $I(G)$, where $I(G)$ is the final squarefree lexsegment ideal determined by the monomial $u=x_1x_{i_2}$, with $i_2>2$. Then 
$$I(G)=\left(\bigcap\limits_{s\geq i_2}P_{[n]\setminus\{s\}}\right)\cap\left(\bigcap\limits_{s<i_2}P_{[n]\setminus\{1,s\}}\right).$$
The minimal vertex covers of $G$ are the sets $[n]\setminus\{s\}$, with $s\geq i_2$ together with the sets $[n]\setminus\{1,s\}$, with $s<i_2$. 
\end{Corollary1}

\medskip

Let $I=(L(u,v))$ be a completely squarefree lexsegment ideal. By \cite{B}, one may write $I=(L^i(v))\cap(L^f(u))$. This allows us to determine the standard primary decomposition of an arbitrary completely squarefree lexsegment ideal.

\begin{Theorem1}\label{prim dec completely}
Let $I\subset S$ be a completely squarefree lexsegment ideal generated in degree $q$, determined by the monomials $u=x_1x_F$ and $v=x_{j_1}\cdots x_{j_q}$, $2\leq j_1<\ldots<j_q\leq n$. Then the minimal primary decomposition of $I$ is the following:
$$I=\bigcap\limits_{|A_t|\leq n-q}P_{A_t}\cap\bigcap\limits_\twoline{|A_t|=n-q+1}{u/x_1\geq_{lex} v/x_{j_t}}P_{A_t}\cap\bigcap\limits_\twoline{G\subset[n],\ |G|=q-1,\ G\cap A_t\neq\emptyset,\ \forall t}{u/x_1\geq_{lex} x_G}P_{G^c}\cap\bigcap\limits_\twoline{P\in \Min(L^f_{S}(u))}{\height P=n-q} P,$$
if $x_2\nmid u$, and
$$I=\bigcap\limits_{|A_t|\leq n-q}P_{A_t}\cap\bigcap\limits_\twoline{|A_t|=n-q+1}{u/x_1\geq_{lex} v/x_{j_t}}P_{A_t}\cap\bigcap\limits_\twoline{G\subset[n]\setminus\{1\},\ |G|=q-1,\ G\cap A_t\neq\emptyset,\ \forall t}{u/x_1\geq_{lex} x_G}P_{G^c}\cap$$
$$\cap\bigcap\limits_\twoline{G\subset[n]\setminus\{1\},\ |G|=n-q+1}{x_{G\setminus\min(G)}\geq_{lex}x_{F^c\setminus\{1\}}}P_{G}\cap\bigcap\limits_\twoline{P\in \Min(L^f_{S}(u))}{\height P=n-q} P,\mbox{ otherwise.}$$
\end{Theorem1}

\begin{proof}
To begin with, we describe the facets of the simplicial complex $\Delta$ associated with $I$ which have cardinality greater than $q$. 

Let $G$ be a facet of $\Delta$ with $|G|>q$. Then $x_G\notin I$ and for all $i\in[n]\setminus G$, $x_{G\cup\{i\}}\in \shad^e(I)$, for some $e>0$, where $\shad^e(I)$ is the $e-$th shadow of $I$. Since $I$ is a completely squarefree lexsegment ideal, $$\shad(I)=\shad(L^i(v))\cap\shad(L^f(u))=\shad(L^i(v))\cap\Mon_{q+1}^s(S)=\shad(L^i(v)).$$ 

We also note that since $|G|>q$, then $x_G\in\shad^e(L^f(u))$, thus $x_G\in (L^f(u))$. Therefore $G$ is a facet of $\Delta$ if and only if $x_G\notin(L^i(v))$ and $x_{G\cup\{i\}}\in \shad^e(L^i(v))$, for some $e>0$. This is equivalent to the fact that $G$ is a facet of $\Gamma_1$, where $\Gamma_1$ is the simplicial complex associated with $(L^i(v))$. By Theorem \ref{prim dec initial}, it follows that $G=[n]\setminus A_t$, for some $1\leq t\leq q$, where $A_t=[j_t]\setminus\{j_1,\ldots,j_{t-1}\}$ and $|A_t|<n-q$.

In the second step of the proof, we describe the facets of $\Delta$ of cardinality $q$. Let $G$ be a facet of $\Delta$ with $|G|=q$. Then $x_G\notin I$ and $x_{G\cup\{i\}}\in I$, for all $i\in[n]\setminus G$, since, as we noticed above, $\shad(L^f(u))=\Mon_{q+1}^s(S)$. Therefore $G$ is a facet of $\Delta$ with $|G|=q$ if and only if $x_G\notin G(I)=L(u,v)$ and $x_{G\cup \max([n]\setminus G)}\geq_{lex}vx_{\max([n]\setminus\supp(v))}$. In other words, $G$ is a facet of $\Delta$ if and only if $G$ is a facet of $\Gamma_2$, where $\Gamma_2$ is the simplicial complex associated with $(L^f(u))$. This is equivalent to $P_{G^c}\in\Min(L^f(u))$, or $x_G<_{lex}v$ and $x_{G\cup \max([n]\setminus G)}\geq_{lex}vx_{\max([n]\setminus\supp(v))}$. The later condition says that $G$ is a facet of $\Gamma_1$, which means that $G=[n]\setminus A_t$, for some $t$ such that $|A_t|=n-q$. 

Finally, let us describe the facets $G$ of $\Delta$ with $|G|=q-1$. We have that $G$ is a facet of $\Delta$ if and only if $x_{G\cup\{i\}}\in I$, for all $i\in[n]\setminus G$, hence if and only if $G$ is a facet of $\Gamma_1$ and $x_{G\cup \min([n]\setminus G)}\leq_{lex}u$. 

We have to consider the following cases:

\textit{Case 1:} There is an integer $1\leq t\leq q$ such that $G=[n]\setminus A_t$, that is $G=\{j_1,\ldots,j_{t-1}\}\cup\{j_t+1,\ldots,n\}$. Then we obtain $\min([n]\setminus G)=1$ and the condition $x_{G\cup \min([n]\setminus G)}\leq_{lex}u$ is equivalent to $x_G\leq_{lex}u/x_1$.

Since $|G|=q-1$, it follows that $j_t-(t-1)=n-q+1$, that is $j_t=n-q+t$. In this case, $v=x_{j_1}\cdots x_{j_{t-1}}x_{j_t}x_{j_t+1}\cdots x_n$, thus $x_G=v/x_{j_t}$. Therefore, $G=[n]\setminus A_t$ is a facet of $\Delta$ of cardinality $q-1$ if and only if $u/x_1\geq_{lex} v/x_{j_t}$. 

\textit{Case 2:} Let $G$ be a facet of $\Gamma_1$ with $|G|=q-1$ and $G\cap A_t\neq\emptyset$ for all $1\leq t\leq q$, such that $x_{G\cup \min([n]\setminus G)}\leq_{lex}u$. If $1\notin G$, then $\min([n]\setminus G)=1$, hence $x_G\leq_{lex}u/x_1$. 

Let now consider $1\in G$ and assume that $i_2\geq 3$. In this case, the condition $x_{G\cup \min([n]\setminus G)}\leq_{lex}u$ is equivalent to $x_{G\setminus\{1\}}x_{\min([n]\setminus G)}\leq_{lex}u/x_1=x_{i_2}\cdots x_{i_q}$, which is imposible since $x_2\mid x_{G\setminus\{1\}}x_{\min([n]\setminus G)}$. Therefore, in this case, the proof is completed.

What is left is to consider $1\in G$ and $i_2=2$. Then $G$ must satisfy the following conditions: $G\cap A_t\neq \emptyset$, for all $1\leq t\leq q$  and $x_{([n]\setminus G)\setminus \{\min([n]\setminus G)\}}\geq_{lex}x_{F^c\setminus \{1\}}$, the later one being equivalent to $x_{G\setminus\{1\}}x_{\min([n]\setminus G)}\leq_{lex}u/x_1$.  
\end{proof}

\begin{Example1}\rm
Let $I=(L(x_1x_3x_4x_5, x_3x_4x_6x_7))\subset k[x_1,\ldots ,x_7]$ be a completely squarefree lexsegment ideal. Then the minimal primary decomposition of $I$ is of the form:
$$I=(x_1,x_2,x_3)\cap (x_1,x_2,x_4)\cap (x_1,x_2,x_5,x_6)\cap (x_1,x_2,x_5,x_7)\cap(x_1,x_2,x_6,x_7)\cap$$
$$\cap(x_3,x_4,x_5)\cap(x_3,x_4,x_6)\cap(x_3,x_4,x_7)\cap(x_3,x_5,x_6)\cap(x_3,x_5,x_7)\cap$$
$$\cap(x_3,x_6,x_7)\cap(x_4,x_5,x_6)\cap(x_4,x_5,x_7)\cap(x_4,x_6,x_7)\cap(x_5,x_6,x_7).$$
\end{Example1}

\begin{Example1}\rm
Let $I=(L(x_1x_2x_4x_5, x_3x_4x_5x_7))\subset k[x_1,\ldots ,x_7]$ be a completely squarefree lexsegment ideal. In this case, the minimal primary decomposition of $I$ is of the form:
$$I=(x_1,x_2,x_3)\cap (x_1,x_2,x_4)\cap(x_1,x_2,x_5)\cap(x_1,x_2,x_6,x_7)\cap (x_1,x_3,x_6,x_7)\cap$$
$$\cap(x_1,x_3,x_5,x_7)\cap(x_1,x_3,x_5,x_6)\cap (x_1,x_3,x_4,x_7)\cap(x_1,x_3,x_4,x_6)\cap$$
$$\cap(x_1,x_3,x_4,x_5)\cap(x_2,x_3,x_4,x_5)\cap(x_2,x_3,x_4,x_6)\cap(x_2,x_3,x_4,x_7)\cap$$
$$\cap(x_2,x_3,x_5,x_6)\cap(x_2,x_3,x_5,x_7)\cap(x_2,x_3,x_6,x_7)\cap(x_4,x_5,x_6)\cap$$
$$\cap(x_4,x_5,x_7)\cap(x_4,x_6,x_7)\cap(x_5,x_6,x_7).$$
\end{Example1}

\begin{Corollary1}\label{dim+multiplicity completely}
Let $I\subset S$ be a completely squarefree lexsegment ideal determined by the monomials $u=x_1x_{i_2}\cdots x_{i_q}=x_1x_F$, $2\leq i_2<\ldots<i_q\leq n$ and $v=x_{j_1}\cdots x_{j_q}$, $2\leq j_1<\ldots<j_q\leq n$, and $\Delta$ be the simplicial complex with the Stanley--Reisner ideal $I$. Then:
\begin{itemize}
	\item [(a)] $\dim(S/I)=n-j_1$.
	\item [(b)] Let $s$ be the smallest integer such that $j_i=j_1+(i-1)$, for all $1\leq i< s$ and $j_s\geq j_1+s$, and $t=|\{x_G: x_{F^c\setminus\{1\}}>_{lex}x_G, |G|=n-q\}|$. Then the multiplicity of $k[\Delta]$ is 
	\[e(k[\Delta])=\left\{\begin{array}{cc}
			s-1&,\mbox{ if }j_1<n-q \\
			s+t-1 &,\mbox{ if }j_1=n-q.   \\
	\end{array}\right.\]
	\item [(c)] If $I$ is generated in degree $2$ and $u\neq x_1x_2$, $v\neq x_{n-1}x_{n}$, then
	$$I=P_{A_1}\cap P_{A_2}\cap\left(\bigcap\limits_{i_2\leq s<j_1}P_{[n]\setminus\{s\}}\right)\cap\left(\bigcap\limits_{s<i_2}P_{[n]\setminus\{1,s\}}\right),\mbox{ if }j_2\leq n-1,$$
and 
$$I=P_{A_1}\cap\left(\bigcap\limits_{i_2\leq s\leq j_1}P_{[n]\setminus\{s\}}\right)\cap\left(\bigcap\limits_{s<i_2}P_{[n]\setminus\{1,s\}}\right),\mbox{ if }j_2=n,$$
where $A_1=[j_1]$ and $A_2=[j_2]\setminus\{j_1\}$.
\end{itemize}
\end{Corollary1}

\begin{proof} The statements (a) and (c) follow by Theorem \ref{prim dec completely}.

(b) The multiplicity of $k[\Delta]$ is $e(k[\Delta])=f_{n-j_1-1}$. To determine the number of facets of $\Delta$ of cardinality $n-j_1$ is equivalent to determine the number of minimal prime ideals which contains $I$ of $\height(p)=j_1$. Let $s$ be the smallest integer such that $j_i=j_1+(i-1)$, for all $1\leq i< s$ and $j_s\geq j_1+s$. Then $|A_i|=j_i-(i-1)=j_1$, for all $1\leq i\leq s-1$ and $|A_i|=j_i-(i-1)>j_1$, for $i\geq s$. By Theorem \ref{prim dec completely}, if $j_1<n-q$, then $e(k[\Delta])=s-1$, and if $j_1=n-q$, then $e(k[\Delta])=s-1+t$, where $t$ is the cardinality of the set $\{x_G: x_{F^c\setminus\{1\}}>_{lex}x_G, |G|=n-q\}$.
\end{proof}

In the hypothesis of the Corollary \ref{dim+multiplicity completely} (c), from the above formula one derives the minimal vertex covers of the graph $G$ with $I(G)=I$.

\section[Sequentially Cohen--Macaulay ideals]{Completely squarefree lexsegment ideals which are sequentially Cohen--Macaulay}

This section is devoted to determining the completely squarefree lexsegment ideals which are sequentially Cohen--Macaulay. 

For an ideal $I$, we will denote by $I^\vee$ its Alexander dual. By Theorem \ref{comp lin-seq}, to prove that the ideal $I$ is sequentially Cohen--Macaulay is equivalent to prove that its Alexander dual is componentwise linear. If $I\subset S$ is a graded ideal, then we denote by $I_{\langle j\rangle}$ the ideal generated by all homogeneous polynomials of degree $j$ belonging to $I$. A graded ideal $I\subset S$ is componentwise linear if $I_{\langle j\rangle}$ has a linear resolution for all $j$.

As we did for the primary decomposition, we will start the study of the sequentially Cohen--Macaulay property with the case of initial squarefree lexsegment ideals. Moreover, we will prove a more general result, which provides a new class of componentwise squarefree ideals. 

In \cite{MH}, the so-called canonical critical ideals have been studied. We recall the definition. A homogeneous ideal $I\subset S$ is called \textit{canonical critical} if it is of the form $I=(f_1x_1,f_1f_2x_2,\ldots,f_1\cdots f_{s-1}x_{s-1},f_1f_2\cdots f_s)$, for some homogeneous polynomials $f_1,\ldots,f_s$, with $f_i\in k[x_i,\ldots ,x_n]$, for each $1\leq i\leq s$ and with $\deg(f_s)>0$, where $1\leq s\leq n$. 

We consider now the class of squarefree monomial (canonical) critical ideals. Let $S=k[x_1,\ldots,x_n]$ be the polynomial ring in $n$ variables over a field $k$ and $I$ be a squarefree monomial ideal in $S$. We denote by $G(I)$ the minimal system of monomial generators of $I$.

\begin{Definition1}\rm
A squarefree monomial ideal $I\subset S$ is called \textit{critical} if it is obtained by the following recursive procedure:
\begin{enumerate}
	\item [(a)] The ideal $I$ is of the form $I=(x_i,m)$, for some squarefree monomial $m$ with $x_i\nmid m$,
	\item [(b)] There is a variable $x_j$, a squarefree monomial $m'$ and $J$ a critical ideal, with $x_j\nmid mm'$ and $\supp(m)\cap\supp(m')=\emptyset$, for all the monomials $m\in G(J)$, such that $I=(x_j)+m'J$.
\end{enumerate}
\end{Definition1}

\begin{Definition1}\rm
A squarefree monomial ideal $I\subset S$ is called \textit{canonical critical} if there exists $w$ a monomial in $S$ such that $I=wJ$, where $J$ is a critical squarefree monomial ideal.
\end{Definition1}

\begin{Proposition1}\label{critical comp lin}
Let $I\subset S$ be a critical squarefree monomial ideal. Then $I$ has linear quotients. In particular, any canonical critical squarefree monomial ideal has linear quotients. 

Moreover, if $I$ is a canonical critical ideal, then $I$ is componentwise linear. 
\end{Proposition1}

\begin{proof}
We use induction on $|G(I)|$.

If $|G(I)|=2$, then $I=(x_i,m)$, for some squarefree monomial $m$ and $x_i\nmid m$. Since $(x_i):m=(x_i)$, it follows that $I$ has linear quotients. 
 
Assume that the assertion holds for all critical ideals $J$, with $|G(J)|=s$, with $s\geq 2$. By definition, there is a variable $x_j$, a squarefree monomial $m'$ and $J$ a critical ideal, with $x_j\nmid mm'$ and $\supp(m)\cap\supp(m')=\emptyset$, for all the monomials $m\in G(J)$, such that $I=(x_j)+m'J$ and $|G(I)|=s+1$. 

Let us assume that $G(J)=\{x_i,m_2,\ldots, m_s\}$ such that $J$ has linear quotients with respect to $x_i,m_2,\ldots, m_s$. We order the minimal system of generators of $I$, $G(I)=\{x_j,m'x_i,m'm_2,\ldots ,m'm_s\}$. For every $2\leq t\leq s$, we have 
$$(x_j,m'x_i,\ldots ,m'm_{t-1}):(m'm_t)=(x_j):(m'm_t)+(m'x_i,\ldots ,m'm_{t-1}):(m'm_t)=$$
$$=(x_j)+(x_i,\ldots ,m_{t-1}):(m_t)$$
is generated by variables, by the induction hypothesis. Hence $I$ has linear quotients.

The last statement in the proposition is a consequence of \cite{JZ}.
\end{proof}

This result gives us a new class of componentwise linear ideals. We will use the properties of these ideals for the case of initial squarefree lexsegment ideals.

\begin{Proposition1}\label{sqCM initial}
Let $I\subset k[x_1,\ldots,x_n]$ be the initial squarefree lexsegment ideal, generated in degree $q$, determined by the monomial $v=x_{j_1}\cdots x_{j_q}$, with $2\leq j_1<\cdots<j_q\leq n$. Then $I$ is sequentially Cohen--Macaulay.
\end{Proposition1}

\begin{proof}
We need to prove that the Alexander dual $I^\vee$ is componentwise linear. By the primary decomposition, Theorem \ref{prim dec initial}, one has $I^\vee=J+K$, where $J=(x_{A_t}:1\leq t\leq q)$ and $K=(x_{G^c}:|G|=q-1,G\cap A_t\neq\emptyset,1\leq t\leq q)$. Since $J$ is generated in degree at most $n-q$ and $K$ is generated in degree $n-q+1$, we obtain that $I^\vee_{\langle j \rangle}=J_{\langle j \rangle}$, for all $j<n-q+1$.

For all $j<n-q+1$, one has that $I^\vee_{\langle j \rangle}=J_{\langle j \rangle}$. Therefore $I^\vee_{\langle j \rangle}$ has a linear resolution, since $J$ is a canonical critical squarefree monomial ideal. Indeed, we have 
$$J=x_1\cdots x_{j_1-1}(x_{j_1}+x_{j_1+1}\cdots x_{j_2-1}(x_{j_2}+x_{j_2+1}\cdots x_{j_3-1}(\ldots))).$$

We prove that $I^\vee_{\langle n-q+1 \rangle}=I_{n,n-q+1}$, that is, $I^\vee_{\langle n-q+1 \rangle}$ is the ideal generated by all the squarefree monomials of degree $n-q+1$ in $S$. This will end our proof since $I_{n,n-q+1}$ has a linear resolution.

Let $m=x_{F^c}\in G(I_{n,n-q+1})$, with $|F|=q-1$, be a squarefree monomial. If $F\cap A_t\neq\emptyset$, for all $1\leq t\leq q$, then $m\in G(I^\vee_{\langle n-q+1 \rangle})$. Assume that there is an integer $1\leq t\leq q$ such that $F\cap A_t=\emptyset$. It results that $F^c\supset A_t$, hence $x_{A_t}\mid m$, thus $m\in G(I^\vee_{\langle n-q+1 \rangle})$. 

The other inclusion is trivial, namely $I^\vee_{\langle n-q+1 \rangle}\subseteq I_{n,n-q+1}$.
\end{proof}

The final squarefree lexsegment ideals are sequentially Cohen--Macaulay, as it follows from the next result.

\begin{Proposition1}\label{sqCM final}
Let $I\subset S$ be the final squarefree lexsegment ideal generated in degree $q$, determined by the monomial $u=x_1x_{i_2}\cdots x_{i_q}$, $2\leq i_2<\ldots<i_q\leq n$. Then $I$ is sequentially Cohen--Macaulay.
\end{Proposition1}

\begin{proof}
We may assume that $I\neq I_{n,q}$. We will prove that the Alexander dual of $I$, $I^\vee$, is componentwise linear. 

It is easy to see, by Theorem \ref{prim dec final}, that $I^\vee_{\langle n-q\rangle}$ has a linear resolution, since it is a final squarefree lexsegment ideal. 

We show that $I^\vee_{\langle n-q+1\rangle}$ is the ideal generated by all the squarefree monomials of degree $n-q+1$. Indeed, we have the inclusion $I^\vee_{\langle n-q+1\rangle}\subset I_{n,n-q+1}$. For the other one, let $m=x_{G^c}$ be a squarefree monomial of degree $n-q+1$. It is clear that, using the notation $u=x_1x_F$, if $m=x_{G^c}\geq_{lex}x_{F^c}$, then $m\in G(I^\vee_{\langle n-q+1\rangle})$. Otherwise, assume that $m=x_{G^c}<_{lex}x_{F^c}$. We have to analyze two cases:

\textit{Case 1:} If $x_{F^c\setminus\{1\}}>_{lex}x_{G^c\setminus\min(G^c)}$, then $x_{G^c\setminus\min(G^c)}\in G(I^\vee_{\langle n-q\rangle})$ and $m=x_{G^c}\in I^\vee_{\langle n-q+1\rangle}$.

\textit{Case 2:} If $x_{F^c\setminus\{1\}}\leq_{lex}x_{G^c\setminus\min(G^c)}$, then $m=x_{G^c}\in G(I^\vee_{\langle n-q+1\rangle})$, which ends the proof. 
\end{proof}

We now focus on arbitrary completely squarefree lexsegment ideals. In order to characterize the completely squarefree lexsegment ideals which are sequentially Cohen--Macaulay, we have to establish when $I^\vee_{\langle j\rangle}$ has a linear resolution, for every $j\leq n-q+1$. Analyzing the minimal primary decompositions obtained if $x_2\mid u$ or $x_2\nmid u$, one can see that

\begin{eqnarray}
I^\vee_{\langle j\rangle}=(x_{A_t}:|A_t|\leq n-q)_{\langle j\rangle}\mbox{ for all }j<n-q\mbox{ and } \eqname{$1$}\label{1}
\end{eqnarray}

\begin{eqnarray}
I^\vee_{\langle n-q\rangle}=(x_{A_t}:|A_t|\leq n-q)_{\langle n-q\rangle}+(x_G:x_{F^c\setminus\{1\}}>_{lex}x_G,|G|=n-q), \eqname{$2$}\label{20}
\end{eqnarray}
\noindent
where $u=x_1x_F$ and $A_t=[j_t]\setminus\{j_1,\ldots,j_{t-1}\},\mbox{ for }1\leq t\leq q.$ Note that only for $I^\vee_{\langle n-q+1\rangle}$ we have to treat separate cases, given by the conditions $x_2\mid u$ or $x_2\nmid u$. 

\begin{Lemma1}\label{dual lower deg}
For every $j<n-q$, the ideal $I^\vee_{\langle j\rangle}$ has a linear resolution.
\end{Lemma1}

\begin{proof}
The ideal $I^\vee_{\langle j\rangle}$ has a linear resolution, for every $j<n-q$, since it is a canonical critical squarefree monomial ideal. 
\end{proof}

In order to determine when
$$I^\vee_{\langle n-q\rangle}=(x_{A_t}:|A_t|\leq n-q)_{\langle n-q\rangle}+(x_G:x_{F^c\setminus\{1\}}>_{lex}x_G,|G|=n-q)$$
has a linear resolution, we need some preparatory results.

Firstly, let $1\leq s\leq q$ be the unique index with the property that $j_s\leq n-q+s-1$ and $j_{s+1}>n-q+s$. Then we have $(x_{A_t}:|A_t|\leq n-q)_{\langle n-q\rangle}=(x_{A_1},x_{A_2},\ldots ,x_{A_s})_{\langle n-q\rangle}$.

\begin{Proposition1}\label{dual deg n-q}
The ideal $(x_{A_t}:|A_t|\leq n-q)_{\langle n-q\rangle}$ is an initial squarefree lexsegment ideal determined by the monomial $x_{A_s}x_{q+j_s-s+2}\cdots x_n$. 
\end{Proposition1}

\begin{proof}We start by noticing that if the cardinality of $A_s$ equals $n-q$, then $x_{A_s}=x_{A_s}x_{q+j_s-s+2}\cdots x_n$.

For the inclusion "$\subseteq$", let $m$ be a minimal monomial generator of $(x_{A_t}:|A_t|\leq n-q)_{\langle n-q\rangle}=(x_{A_1},x_{A_2},\ldots ,x_{A_s})_{\langle n-q\rangle}$, that is $m=x_{A_t}m'$, for some $1\leq t\leq s$ with $A_t\cap \supp(m')=\emptyset$ and we assume that $t$ is the smallest with this property. We have to prove that $m\geq_{lex} x_{A_s}x_{q+j_s-s+2}\cdots x_n$. 

Assume by contradiction that $x_{A_s}x_{q+j_s-s+2}\cdots x_n>_{lex}m=x_{A_t}m'$. We observe that $x_{l}\mid x_{A_t}m'$, whence $x_{l}\mid x_{A_s}x_{q+j_s-s+2}\cdots x_n$, for all $l\in A_t$, $l\neq j_t$. Thus if $s\neq t$, then we reach a contradiction. Therefore, we must have $t=s$ and by the relation $x_{A_s}x_{q+j_s-s+2}\cdots x_n>_{lex}m=x_{A_t}m'$ we get $x_{q+j_s-s+2}\cdots x_n>_{lex}m'$. This is a contradiction, hence $m\geq_{lex} x_{A_s}x_{q+j_s-s+2}\cdots x_n$.

Conversely, for the inclusion "$\supseteq$", let $m\geq_{lex} x_{A_s}x_{q+j_s-s+2}\cdots x_n$ be a squarefree monomial of degree $n-q$. We claim that $\supp(m)\cap\supp(v)\neq\emptyset$. Indeed, if we assume that $\supp(m)\cap\supp(v)=\emptyset$, by degree considerations, we obtain that $m=x_{[n]\setminus\supp(v)}=x_{\{1,\ldots n\}\setminus\{j_1,\ldots,j_q\}}$. Hence  
$$x_{A_s}x_{q+j_s-s+2}\cdots x_n=x_{\{1,\ldots,j_s\}\setminus\{j_1,\ldots,j_{s-1}\}}x_{q+j_s-s+2}\cdots x_n>_{lex}x_{\{1,\ldots n\}\setminus\{j_1,\ldots,j_q\}}=m,$$
which is a contradiction. 

Therefore $\supp(m)\cap\supp(v)\neq\emptyset$. Let $j_t=\min\{i\in[n]:i\in \supp(m)\cap\supp(v)\}$. We claim that 
$$\{r\in \supp(m): r\leq j_t\}=[j_t]\setminus\{j_1,\ldots,j_{t-1}\}=A_t.$$

Indeed, by hypothesis, we have 
$$m=x_{\{r\in \supp(m): r\leq j_t\}}x_{\{r\in \supp(m): r> j_t\}}\geq_{lex} x_{A_s}x_{q+j_s-s+2}\cdots x_n=$$
$$=x_{\{1,\ldots,j_t\}\setminus\{j_1,\ldots,j_{t-1}\}}x_{\{j_t,\ldots,j_s\}\setminus\{j_t,\ldots,j_{s-1}\}}x_{q+j_s-s+2}\cdots x_n.$$
Since $j_t$ is minimal, we get $\{r\in \supp(m): r\leq j_t\}\cap\{j_1,\ldots,j_{t-1}\}=\emptyset$. Hence $x_{\{r\in \supp(m): r\leq j_t\}}=x_{[j_t]\setminus\{j_1,\ldots,j_{t-1}\}}$. Otherwise $x_{\{r\in \supp(m): r\leq j_t\}}<_{lex}x_{[j_t]\setminus\{j_1,\ldots,j_{t-1}\}}$, thus $m<_{lex}x_{A_s}x_{q+j_s-s+2}\cdots x_n$, a contradiction.

It results that $\{r\in \supp(m): r\leq j_t\}=[j_t]\setminus\{j_1,\ldots,j_{t-1}\}=A_t$ and $x_{A_t}\mid m$, thus $m\in (x_{A_t}:|A_t|\leq n-q)_{\langle n-q\rangle}$, which ends the proof.
\end{proof}

Returning to the ideal $I^\vee_{\langle n-q\rangle}$, by (\ref{20}), it follows that $I^\vee_{\langle n-q\rangle}$ is the sum of an initial with a final squarefree lexsegment ideal, both of them generated in the same degree. In the following, in order to determine the ideals $I^\vee_{\langle n-q\rangle}$ with a linear resolution, we will analyze a more general problem.

\medskip
\medskip
\medskip
We consider $J=(L^i(w))$ and $K=(L^f(m))$ be initial and final squarefree lexsegment ideals, generated in degree $d$, such that $x_1\mid w$ and $x_1\nmid m$. The ideals $J$ and $K$ have a $d-$linear resolution. We are interested in determining when the sum $J+K$ has a $d-$linear resolution.

The statement of the next proposition is probably known. We give a short proof for the convenience of the reader.

\begin{Proposition1}\label{intersection}
The ideal $J+K$ has a $d-$linear resolution if and only if the ideal $J\cap K$ has a $(d+1)-$linear resolution.
\end{Proposition1}

\begin{proof}
We consider the following exact sequence of $S-$modules
\begin{eqnarray}
0\rightarrow J\cap K\rightarrow J\oplus K\rightarrow J+K\rightarrow 0.\eqname{$3$}\label{30}
\end{eqnarray}

We obviously have $\reg(J\cap K)\geq d+1$ since the initial degree of $J\cap K$ is $d+1$. If $J+K$ has a linear resolution, equivalently, $\reg(J+K)=d$, then, by using the exact sequence (\ref{30}), we get $\reg(J\cap K)\leq\max\{d,d+1\}=d+1$. Therefore, $\reg(J\cap K)=d+1$, which means that $J\cap K$ has a $(d+1)-$linear resolution. 

The proof of the converse works similarly. 
\end{proof}

\begin{Lemma1}\label{intersection d+1}
The ideal $J\cap K$ is generated in degree $d+1$ if and only if $J\cap K=(L(x_1m,wx_{\max([n]\setminus\supp(w))}))$.
\end{Lemma1}

\begin{proof}
If $J\cap K$ is generated in degree $d+1$, then it is easy to see that $G(J\cap K)=L(x_1m,wx_{\max([n]\setminus\supp(w))})$. Indeed, 
$$G(J\cap K)\subseteq L^i(wx_{\max([n]\setminus\supp(w))})\cap L^f(x_1m)=L(x_1m,wx_{\max([n]\setminus\supp(w))}).$$
For the other inclusion, let $\gamma\in L(x_1m,wx_{\max([n]\setminus\supp(w))})$. It is clear that $x_1\mid \gamma$. We obtain that $m\geq \gamma/x_1$, that is $\gamma\in K$. Moreover, $\gamma\in J$ because $\gamma/x_{\max(\gamma)}\geq_{lex}w$. Hence $\gamma\in J\cap K$, thus $\gamma\in G(J\cap K)$, as desired.

The converse is clear.
\end{proof}

\medskip
\medskip

We return to the problem of characterizing the completely squarefree lexsegment ideals which are sequentially Cohen--Macaulay.

Let $I\subset S$ be a completely squarefree lexsegment ideal generated in degree $q$, determined by the monomials $u=x_1x_{i_2}\cdots x_{i_q}=x_1x_F$, $2\leq i_2<\ldots<i_q\leq n$ and $v=x_{j_1}\cdots x_{j_q}$, with $2\leq j_1<\ldots<j_q\leq n$, which is not an initial or a final squarefree lexsegment ideal. Denote by $A_t=[j_t]\setminus\{j_1,\ldots,j_{t-1}\}$ for $1\leq t\leq q$. Let $s$ be the unique index such that $|A_s|\leq n-q$ and $|A_{s+1}|=n-q+1$.

Recall that if $x_2\nmid u$, then 

$$I^\vee_{\langle n-q+1\rangle}=(x_{A_t}:|A_t|\leq n-q)_{\langle n-q+1\rangle}+(x_{A_t}:|A_t|=n-q+1, u/x_1\geq_{lex} v/x_{j_t})+$$
$$+(x_{G^c}:G\subset[n],\ |G|=q-1,\ G\cap A_t\neq\emptyset,1\leq t\leq q,u/x_1\geq_{lex} x_G)+$$
$$+(x_G: x_{F^c\setminus\{1\}}>_{lex}x_G,|G|=n-q)_{\langle n-q+1\rangle}$$
\noindent
and if $x_2\mid u$, then 

$$I^\vee_{\langle n-q+1\rangle}=(x_{A_t}:|A_t|\leq n-q)_{\langle n-q+1\rangle}+(x_{A_t}:|A_t|=n-q+1, u/x_1\geq_{lex} v/x_{j_t})+$$
$$+(x_{G^c}:G\subset[n]\setminus\{1\},\ |G|=q-1,\ G\cap A_t\neq\emptyset,1\leq t\leq q,u/x_1\geq_{lex} x_G)+$$
$$+(x_G:x_{F^c\setminus\{1\}}>_{lex}x_G,|G|=n-q)_{\langle n-q+1\rangle}+(x_G:x_{G\setminus\min(G)}\geq_{lex}x_{F^c\setminus\{1\}},|G|=n-q+1).$$

\begin{Theorem1}\label{sqCM completely}
With the above notations, $I$ is sequentially Cohen--Macaulay if and only if $(L^i(x_{A_s}x_{q+j_s-s+2}\cdots x_n))\cap(L^f(\succ(x_{F^c\setminus\{1\}})))$ has a $(n-q+1)-$linear resolution.
\end{Theorem1}

\begin{proof}
We have to establish when the Alexander dual $I^\vee$ is componentwise linear. By Lemma \ref{dual lower deg}, $I^\vee_{\langle j\rangle}$ has a linear resolution, for all $j<n-q$.

By Proposition \ref{dual deg n-q}, we have 
$$I^\vee_{\langle n-q\rangle}=(L^i(x_{A_s}x_{q+j_s-s+2}\cdots x_n))+(L^f(\succ(x_{F^c\setminus\{1\}}))).$$
Hence $I^\vee_{\langle n-q\rangle}$ has a linear resolution if and only if, using Proposition \ref{intersection}, 
$$(L^i(x_{A_s}x_{q+j_s-s+2}\cdots x_n))\cap(L^f(\succ(x_{F^c\setminus\{1\}})))$$
has a $(n-q+1)-$linear resolution. 

We end the proof by showing that $I^\vee_{\langle n-q+1\rangle}=I_{n,n-q+1}$. 

\textit{Case 1:} We assume that $x_2\nmid u$. The inclusion $I^\vee_{\langle n-q+1\rangle}\subseteq I_{n,n-q+1}$ is clear. For the other one, let $m=x_{H^c}$ be a squarefree monomial of degree $n-q+1$. If $x_{F^c}>_{lex}m$, then $m\in (L^f(\succ(x_{F^c})))=(L^f(\succ(x_{F^c\setminus\{1\}})))_{\langle n-q+1\rangle}$.

Assume that $m=x_{H^c}\geq_{lex}x_{F^c}$, equivalently $x_H\leq_{lex} x_{F}$. If $H\cap A_t\neq\emptyset$ for all $1\leq t\leq q$, then $m\in (x_{G^c}: |G|=q-1, G\subset[n],G\cap A_t\neq\emptyset,1\leq t\leq q,u/x_1\geq_{lex} x_G)\subseteq I^\vee_{\langle n-q+1\rangle}$. Otherwise, if there is some $t$ such that $H\cap A_t=\emptyset$, then we must have $j_t-t+1+q-1\leq n$, that is $j_t-t+1\leq n-q+1$. If $j_t-t+1\leq n-q$, then $x_{A_t}\in (x_{A_t}:|A_t|\leq n-q)$. Moreover, by $H\cap A_t=\emptyset$ we obtain $A_t\subseteq H^c$, thus $x_{A_t}\mid m$ and $m\in (x_{A_t}:|A_t|\leq n-q)_{\langle n-q+1\rangle}$. If $j_t-t+1=n-q+1$, then $v=x_{j_1}\cdots x_{j_{t-1}}x_{n-q+t}\cdots x_n$. In particular, $H^c=A_t$ and by the relation $m=x_{H^c}=x_{A_t}\geq_{lex}x_{F^c}$ we get $u/x_1=x_F\geq_{lex} v/x_{j_t}$. Hence $m\in (x_{A_t}: |A_t|=n-q+1,\ u/x_1\geq_{lex} v/x_{j_t})\subseteq I^\vee_{\langle n-q+1\rangle}$, which ends the proof.

\textit{Case 2:} For the case when $x_2\mid u$, we have the following. Let $m=x_{H^c}$ be a squarefree monomial of degree $n-q+1$. 

We analyze separately the cases when $x_1\mid m$ and $x_1\nmid m$.

Firstly, if $x_1\nmid m$, then either $x_{H^c\setminus\min(H^c)}\geq_{lex}x_{F^c\setminus\{1\}}$, thus $m=x_{H^c}\in G(I^\vee_{\langle n-q+1\rangle})$, or $x_{H^c\setminus\min(H^c)}<_{lex}x_{F^c\setminus\{1\}}$. This later inequality implies $x_{H^c\setminus\min(H^c)}\in G(I^\vee_{\langle n-q\rangle})$, thus $x_{H^c}\in I^\vee_{\langle n-q+1\rangle}$. 

Assume that $m$ is divisible by $x_1$. By the minimal primary decomposition, it remains to study the following cases:
\begin{itemize}
	\item $H\cap A_t\neq\emptyset$ for all $1\leq t\leq q$ and $u/x_1<_{lex} x_H$;
	\item there exists an integer $1\leq t\leq q$ such that $H\cap A_t=\emptyset$,
\end{itemize}

If we assume that $H\cap A_t\neq\emptyset$ for all $1\leq t\leq q$ and $u/x_1=x_F<_{lex} x_H$, then $m=x_{H^c}<_{lex}x_{F^c}$. Since $x_1\mid x_{H^c}$, we obtain $x_{H^c\setminus\{1\}}\in (L^f(\succ(x_{F^c\setminus\{1\}})))$, thus $x_{H^c}\in I^\vee_{\langle n-q+1\rangle}$.

For the case when there is some $1\leq t\leq q$ such that $H\cap A_t=\emptyset$, we obtain that $|H|+|A_t|\leq n$. This implies that $j_t-t+1\leq n-q+1$. If $|A_t|=j_t-t+1\leq n-q$, then $x_{A_t}\in(x_{A_t}:|A_t|\leq n-q)$ and $x_{A_t}\mid x_{H^c}$, because $A_{t}\subset H^c$. Thus $m=x_{H^c}\in (x_{A_t}:|A_t|\leq n-q)_{\langle n-q+1\rangle}$. Finally, if $|A_t|=j_t-t+1=n-q+1$, then $v=x_{j_1}\cdots x_{j_{t-1}}x_{n-q+t}\cdots x_n$ and $H^c=A_t$. Moreover, if $u/x_1=x_F\geq_{lex} v/x_{j_t}$, then $m=x_{A_t}\in(x_{A_t}:|A_t|=n-q+1, u/x_1\geq_{lex} v/x_{j_t})$. Otherwise, if $u/x_1=x_F<v/x_{j_t}$, we obtain $x_{F^c}>_{lex} x_{A_t}=x_{H^c}$ which implies $x_{H^c\setminus\{1\}}\in G(I^\vee_{\langle n-q\rangle})$, thus $m=x_{H^c}\in I^\vee_{\langle n-q+1\rangle}$.
\end{proof}

\chapter{Homological invariants of squarefree lexsegment ideals}

In this chapter we are interested in characterizing the homological invariants of an arbitrary squarefree lexsegment ideal. In the previous chapter we were able to completely describe the initial and final squarefree lexsegment ideals. For these, we have computed the Krull dimension, the depth, the multiplicity and the Castelnuovo--Mumford regularity. We saw that these invariants are described in terms of the ends of the squarefree lexsegment or in terms of the degree in which the ideal is generated. 

Passing to the completely squarefree lexsegment ideals, we completed the computation of the Krull dimension and of the multiplicity. The full description of the invariants for arbitrary squarefree lexsegment ideals seems to be very difficult. We are able to give complete answers for squarefree lexsegment ideals generated in small degrees.

Based on these results in small degrees, we formulate some conjectures and questions in arbitrary degrees.
 
The results appearing in this chapter can be found in our paper \cite{EOT} and in \cite{EKOT}. There are also some results from \cite{OO} presented in the first section of this chapter.

\section{Bounds for $\depth(S/I)$}
In this section we will give lower and upper bounds for $\depth(S/I)$, where $I\subset S$ is an arbitrary squarefree lexsegment ideal. 

We begin with a very useful lemma.

\begin{Lemma1}\label{lex ideal complementary}
Let $I=(L(x_G,x_H))\subset S$ be a squarefree lexsegment ideal, with $|G|=|H|$. 
Then $I$ has a linear resolution if and only if the ideal $(L(x_{[n]\setminus H},x_{[n]\setminus G}))$ has a linear resolution.
\end{Lemma1}

\begin{proof}
It is enough to demonstrate that if we have a relation between the minimal monomial generators of $(L(x_G,x_H))$, then we obtain the same relation in $(L(x_{[n]\setminus H},x_{[n]\setminus G}))$. 

Let $x_G\geq_{lex} m_1>_{lex}m_2\geq_{lex} x_H$ be two minimal monomial generators of $I$ and let $x_{F_1}m_1-x_{F_2}m_2=0$ be a relation. Denote $m_1=x_{T_1}$ and $m_2=x_{T_2}$. Then the relation $x_{F_1}m_1=x_{F_2}m_2$ is equivalent to $x_{[n]\setminus (T_{1}\cup F_1)}=x_{[n]\setminus (T_2\cup F_2)}$. Equivalently, $x_{F_2}x_{[n]\setminus T_1}=x_{F_1}x_{[n]\setminus T_2}$, which give us the same relation between the minimal monomial generators $x_{[n]\setminus T_1}$ and $x_{[n]\setminus T_2}$ of $(L(x_{[n]\setminus H},x_{[n]\setminus G}))$. In particular, they have the same linear relations. 
\end{proof}

Next, we will give upper and lower bounds for the depth of an arbitrary squarefree lexsegment ideal.

Let $I$ be an arbitrary squarefree lexsegment ideal, not necessarily complete. As it follows from the proof of Corollary \ref{dim} (b), one can see that $\depth(S/I)\geq q-1$. The following result will characterize the ideals with $\depth(S/I)>q-1$.

Denote $\succ(m)=\max\{w: m>_{lex}w\}$ and $\pred(m)=\min\{w: w>_{lex}m\}$.

\begin{Proposition1}
Let $I=(L(u,v))$ be an arbitrary squarefree lexsegment ideal, with $x_1\mid u$ and $x_1\nmid v$. Then $\depth(S/I)>q-1$ if and only if $$\succ(v)/x_{\max(\succ(v))}\geq_{lex}\pred(u)/x_1$$
and $(L(\succ(v)/x_{\max(\succ(v))},\pred(u)/x_1))$ has a linear resolution. 
\end{Proposition1}

\begin{proof}
Let $\Delta$ be the simplicial complex such that $I=I_{\Delta}$. By Lemma \ref{CM skeleton}, one has $\depth(S/I)>q-1$ if and only if $\Delta^{(q-1)}$ is Cohen--Macaulay. Equivalently, by Eagon--Reiner theorem \cite [Theorem 8.1.9]{HeHi}, $I_{(\Delta^{(q-1)})^\vee}$ has a linear resolution. We denote by $\Gamma=\Delta^{(q-1)}$ and by $\Gamma_1=\langle\{1,s_2,\ldots ,s_q\}: x_1x_{s_2}\cdots x_{s_q}>_{lex}u\rangle$ and $\Gamma_2=\langle\{t_1,\ldots ,t_q\}: v>_{lex}x_{t_1}\cdots x_{t_q}\rangle$. It is clear that $\Gamma=\Gamma_1\cup\Gamma_2$, equivalent, by \cite{O1}, to $I_{\Gamma^\vee}=I_{\Gamma_1^\vee}+I_{\Gamma_2^\vee}$. Since $I_{\Gamma_1^\vee}=I(\Gamma_1^c)$, we obtain that $I_{\Gamma_1^\vee}=(L^f(\succ(x_{[n]\setminus\supp(u)}))$. In a similar way, one has $I_{\Gamma_2^\vee}=(L^i(\pred(x_{[n]\setminus\supp(v)}))$. 

Then $I_{\Gamma^\vee}=I_{\Gamma_1^\vee}+I_{\Gamma_2^\vee}$ has a linear resolution if and only if, by Proposition \ref{intersection}, the ideal $I_{\Gamma_1^\vee}\cap I_{\Gamma_2^\vee}$ has an $(n-q+2)-$linear resolution. Equivalently, by Lemma \ref{intersection d+1}, $(L(x_1\succ(x_{[n]\setminus\supp(u)}),\pred(x_{[n]\setminus\supp(v)})x_{\max([n]\setminus([n]\setminus\supp(v)))})$ has an $(n-q+2)-$linear resolution. By Lemma \ref{lex ideal complementary}, this is equivalent to $(L(\succ(v)/x_{\max(\succ(v))},\pred(u)/x_1))$ has a linear resolution. 
\end{proof}

\begin{Proposition1}
Let $I=(L(u,v))$ be an arbitrary squarefree lexsegment ideal with $u\neq v$ generated in degree $q$. Then $\depth(S/I)\leq n-2$.
\end{Proposition1}

\begin{proof}
Firstly, assume by contradiction that $\depth(S/I)=n-1$. Hence $S/I$ is Cohen--Macaulay, or equivalently, $I^\vee$ has a linear resolution. Moreover, $I^\vee$ is generated in degree $1$, since there is a minimal prime ideal of $I$ of height $1$. Then any minimal prime ideal of $I$ has height $1$, thus $I$ is a principal ideal, contradiction.

Hence $\depth(S/I)\leq n-2$. 
\end{proof}

\section{Lexsegment edge ideals}
The results of this section are contained in \cite{EOT}, joint work with V. Ene and N. Terai. Throughout this section, a squarefree lexsegment ideal generated in degree $2$ will be called a \textit{lexsegment edge ideal}.

Let $u=x_1x_i, v=x_jx_r$ be two squarefree monomials of degree $2$ such that $u>_{lex} v$ and $I=(L(u,v))$ the lexsegment edge 
ideal generated by the set $L(u,v)$. We begin the study of invariants with the computation of the Krull dimension. 

The Krull dimension of the initial and final lexsegment edge ideals can be easily determined by the minimal primary decomposition, given in the previous chapter. For arbitrary lexsegment edge ideals, we will present a different proof than in the original paper \cite{EOT}, which involves the completely lexsegment edge ideals. 

\begin{Corollary1}\label{lex edge ideals compl}
Let $u=x_1x_i$, $v=x_jx_r$ and $I=(L(u,v))$ be the lexsegment edge ideal. Then $I$ is a completely lexsegment edge ideal if and only if $j\geq i-2.$
\end{Corollary1} 

\begin{proof}
The statement can be checked easily by showing directly that the squarefree shadow of $L(u,v)$ is a squarefree lexsegment if and only if $j\geq i-2.$
\end{proof} 

\begin{Proposition1}
Let $I=(L(u,v))$ be a lexsegment edge ideal which is neither initial nor final and is determined by $u=x_1x_i$ and $v=x_jx_r.$ Then 
$\dim(S/I)=n-j.$
\end{Proposition1}

\begin{proof}
If we consider that $I$ is a completely lexsegment edge ideal, then the conclusion follows from the Corollary \ref{dim+multiplicity completely}.

Assume that $I$ is a lexsegment edge ideal which is not complete, equivalently, by Corollary \ref{lex edge ideals compl}, to $j<i-2$. In this case, we write
$$I=(x_1)\cap(x_i,\ldots,x_n)+(L(x_2x_3,x_jx_r)).$$ 
Using the minimal primary decomposition of an initial squarefree lexsegment ideal (Theorem \ref{prim dec initial}), we get
$$I=(x_1)\cap(x_i,\ldots,x_n)+(x_t:t\in\{2,\ldots,j\})\cap(x_t:t\in\{2,\ldots,r\}\setminus\{j\})\cap$$
$$\cap\bigcap\limits_{2\leq s\leq j-1}(x_t:t\in \{2,\ldots,n\}\setminus\{s\}).$$
After some computations, we have
$$I=(x_t:t\in\{1,\ldots,j\})\cap(x_t:t\in\{1,\ldots,r\}\setminus\{j\})\cap(x_t:t\in\{2,\ldots,j\}\cup\{i,\ldots,n\})\cap$$
$$\cap(x_t:t\in(\{2,\ldots,r\}\setminus\{j\})\cup\{i,\ldots,n\})\cap\bigcap\limits_{2\leq s\leq j-1}(x_t:t\in \{2,\ldots,n\}\setminus\{s\})$$
which is an irredundant decomposition of $I$. In particular, we obtain that $\height(I)=j$, thus $\dim(S/I)=n-j$, as desired.  
\end{proof}

We are also interested in computing the depth of $S/I$ for an arbitrary lexsegment edge ideal $I.$

\begin{Proposition1}\label{depth1}
Let $u=x_1x_i, v=x_jx_r$ with $j\geq 2,$ and $I=(L(u,v)).$ Then $\depth(S/I)=1$ if and only if $x_{i-1}x_n\geq_{lex} v.$
\end{Proposition1}

\begin{proof}
Let $\Delta$ be the simplicial complex on the vertex set $[n]$ whose Stanley--Reisner ideal is $I.$ It is known that $\depth(S/I)=1$ 
if and only if $\Delta$ is disconnected, which, in turn, is equivalent to the fact that the first skeleton 
$\Delta^{(1)}=\{F\in \Delta\colon \dim F\leq 1\}$ of $\Delta$ is disconnected. 

In the first place we consider $\Delta^{(1)}$ disconnected. Let $V_1,V_2\neq \emptyset,$ $V_1\cup V_2=[n],$ $V_1\cap V_2=\emptyset,$ 
and such that no face of $\Delta^{(1)}$ has vertices in both $V_1$ and $V_2.$ One may assume that $1\in V_1.$ Then, since 
$\{1,2\},\ldots,\{1,i-1\}\in \Delta^{(1)},$ we must have $2,\ldots,i-1\in V_1.$ Let us assume that $v>_{lex}x_{i-1}x_n.$ Then 
$\{\ell, n\}\in \Delta^{(1)}$ for all $\ell\geq i-1$ which implies that $i,\ldots,n\in V_1$ as well. This leads to $V_1=[n]$ which is a 
contradiction to our hypothesis.

For the converse, let $x_{i-1}x_n\geq_{lex} v.$ We claim that $\Delta^{(1)}$ is disconnected. Indeed, one may choose 
$V_1=\{1,\ldots,i-1\}$ and $V_2=\{i,\ldots,n\}$ and observe that for any $1\leq r\leq i-1$ and $i\leq s\leq n$ we have $x_rx_s\in I$,
hence $\{r,s\}\not\in \Delta^{(1)}.$
\end{proof}

\begin{Corollary1}\label{pd1}
Let $u$ and $v$ as in the above proposition. Then $\projdim_S(S/I)=n-1$ if and only if $x_{i-1}x_n\geq_{lex}v.$
\end{Corollary1}

Next we compute the depth of $S/I$ in the case when $v=x_jx_r$ with $j\geq 2$ and $v>_{lex}x_{i-1}x_n.$ This is made in two steps. Firstly, in the next lemma we investigate the case $j\geq 3.$

\begin{Lemma1}\label{depth2}
Let $I=(L(u,v))$ where $u=x_1x_i, v=x_jx_r$, $j\geq 3,$ and $v>_{lex}x_{i-1}x_n.$ Then $\depth(S/I)=2.$
\end{Lemma1}

\begin{proof}
By the hypothesis on $v$ we have $\depth(S/I)\geq 2.$ Let $\Delta$ be the simplicial complex on $[n]$ such that $I=I_{\Delta}.$  We 
claim that $\{1,2\}$ is a facet of $\Delta.$ Indeed, if $3\leq 
p\leq n,$ then $\{1,2,p\}\not\in \Delta$ since $x_2x_p\in I_{\Delta}.$ Therefore, by Proposition \ref{depth-bound}, we get $\depth(S/I)\leq 2$. We therefore 
have 
\[
n-2\leq \projdim_S(S/I)\leq n-2.
\]
It follows that $\projdim_S(S/I)=n-2$ and $\depth(S/I)=2.$
\end{proof}

It remains to consider the case $v=x_2x_r$ for some $r\geq 3.$

\begin{Lemma1}\label{depth3}
Let $u=x_1x_i, v=x_2x_r>_{lex}x_{i-1}x_n$ and $I=(L(u,v)).$ Then
\[
\depth(S/I)=\left\{
\begin{array}{ll}
	2, & \text{ if }\ r\geq i,\\
	i+1-r, & \text{ if }\ i>r.
\end{array}\right.
\]
\end{Lemma1}

\begin{proof}
Let us first consider $r\geq i.$ One may easily see that $I$ has the following primary decomposition
\[
I=(x_1,x_2)\cap(x_1,x_3,\ldots,x_r)\cap(x_2,x_i,\ldots,x_n)\cap(x_3,\ldots,x_n).
\]
Hence $n-2\geq \projdim_S(S/I)\geq n-2$, using Proposition \ref{depth-bound}, which implies $\projdim_S(S/I)=n-2$ and, thus, $\depth(S/I)=2.$

For $i>r$ one checks that the minimal monomial generators of $I,$ let us say, $m_1,\ldots,m_r,$ satisfy the following condition: 
for any $1\leq i\leq r,$ there exists $1\leq j\leq n$ such that $x_j| m_i$ and $x_j\not|m_{\ell}$ for all $\ell\neq i.$ This implies that the Taylor resolution of $S/I$ is minimal and, therefore, $\projdim_S(S/I)$ is equal to the number of the minimal monomial generators of $I$, that is, $\projdim_S(S/I)=n+r-i-1.$ Consequently, $\depth(S/I)=i+1-r.$
\end{proof}

Based on the above formulae for dimension and depth, we can easily recover the characterization of the Cohen--Macaulay lexsegment edge ideals given in \cite{BS2}.

\begin{Corollary1}
Let $I=(L(u,v))$ be a lexsegment edge ideal with $x_1|u$ and $u \ne v$. Then $I$ is Cohen--Macaulay if and only if one of the following conditions holds:
\begin{itemize}
	\item [(i)] $I=I_{n,2}$.
	\item [(ii)] $u=x_1x_n$ and $v\in\{x_2x_3,x_{n-2}x_{n-1},x_{n-2}x_n\}$ for $n\geq 4$.
	\item [(iii)] $u=x_1x_{n-1}$, $v=x_{n-2}x_n$ for $n\geq 3$.
\end{itemize}
\end{Corollary1}

It remains to consider the regularity of the lexsegment edge ideals.

We first notice that if $I$ is an initial or final lexsegment edge ideal, then $\reg(I)=2$ since the ideal has a linear resolution.
Therefore we may consider that $u\neq x_1x_2$, that is, $i\geq 3,$ and $v\neq x_{n-1}x_n,$ in other words, $2\leq j\leq n-2.$

\begin{Lemma1}
Let $I=(L(u,v))$ be a lexsegment edge ideal. Then $\reg(I)\in \{2,3\}.$
\end{Lemma1}

\begin{proof}
The ideal $I$ can be decomposed as $I=J+J^{\prime}$ where $J$ is generated by the lexsegment $L(u,x_1x_n)$ and $J^{\prime}$ by 
$L(x_2x_3,v).$ Both ideals $J$ and $J^{\prime}$ have a linear resolution, hence $\reg(J)=\reg(J^{\prime})=2$. By \cite{KM} (see also \cite{H} and \cite{T}), it follows that $\reg(I)\leq \reg(J)+\reg(J^{\prime})-1=3.$
\end{proof}

This easy lemma shows that we have to distinguish only between two possible values of the regularity of $I.$

The original proof of this result uses the characterization of squarefree lexsegment ideals with a linear resolution given in \cite{B} and \cite{BoS}. We present here a different proof.

\begin{Proposition1}
Let $I=(L(u,v))$ be a lexsegment edge ideal determined by $u=x_1x_i, v=x_jx_r, j\geq 2.$ Then
\[
\reg(I)=\left\{ 
\begin{array}{ll}
	3, & \text{ if } i\geq j+2 \text{ and } x_n\not| v \\
	2, & { otherwise.}
\end{array}\right.
\]
\end{Proposition1}

\begin{proof}

We will analyze the following cases:

\textit{Case 1:} We firstly consider $i<j+2$. In this case, the ideal $I$ has linear quotients with respect to the order 
$$x_2x_3,\ldots,x_jx_r,x_1x_i,\ldots,x_1x_n$$
of the minimal monomial generators. Thus we have $\reg(I)=2$.

\textit{Case 2:} Assume that $i\geq j+2$ and $x_n\mid v$. In this case, we write the ideal $I=J+K$, where $J=(L^f(x_1x_i))=(x_1)\cap (x_i,\ldots,x_n)$ and $K=(L^i(x_jx_r))\subset k[x_2,\ldots,x_n]$. 

The ideal $I=J+K$ has a $2-$linear resolution if and only if $J\cap K$ has a $3-$linear resolution, by Proposition \ref{intersection}. Moreover, $J\cap K$ has a $3-$linear resolution if and only if, by Eagon--Reiner's theorem, $(J\cap K)^\vee$ is Cohen--Macaulay of dimension $n-3$. 

A primary decomposition of $J\cap K$ is of the form
$$J\cap K=(x_1)\cap (x_i,\ldots,x_n)\cap (x_2,\ldots,x_j)\cap(x_2,\ldots,x_{j-1},x_{j+1},\ldots ,x_r)\cap$$
$$\cap\bigcap\limits_{2\leq r\leq j-1}(x_s:s\in\{2,\ldots ,n\}\setminus\{r\}).$$

Since $i\geq j+2$ and $r=n$, we have that the minimal primary decomposition is  
$$J\cap K=(x_1)\cap (x_i,\ldots,x_n)\cap (x_2,\ldots,x_j).$$

This implies that $(J\cap K)^\vee=(x_1,x_i\cdots x_n,x_2\cdots x_j)$ is a complete intersection, thus it is Cohen--Macaulay of dimension $n-3$.

\textit{Case 3:} In the hypothesis $i\geq j+2$ and $x_n\nmid v$, using the same notations as before, we have that the standard primary decomposition of $J\cap K$ is given by
$$J\cap K=(x_1)\cap(x_i,\ldots,x_n)\cap(x_2,\ldots,x_j)\cap(x_2,\ldots,x_{j-1},x_{j+1},\ldots ,x_r).$$
This implies that the ideal $(J\cap K)^\vee$ is not unmixed, since both ideals $(x_1,x_2,x_j)$ and $(x_1,x_j,x_{j+1},x_n)$ are minimal prime containing $(J\cap K)^\vee$. Hence $(J\cap K)^\vee$ is not Cohen--Macaulay. By Eagon--Reiner's theorem, we obtain that $J\cap K$ does not have a linear resolution, thus $\reg(J+K)=3$.
\end{proof}

The arithmetical rank of a squarefree lexsegment edge ideals is computed in the following result.  
\begin{Theorem1}\label{aralei}
Let $I=(L(u,v))$ be a lexsegment edge ideal. Then 
\[
\ara(I)=\projdim_S(S/I).
\]
\end{Theorem1}

\begin{proof} Let $u=x_1x_i$ and $v=x_jx_r$ such that $u\geq_{lex}v.$
In the first place we observe that the statement is obviously true if $j=1$ since, for instance, $I$ is isomorphic as an $S$-module 
to the ideal generated by the variables $x_i,\ldots,x_r.$  Hence we may assume that $j\geq 2.$ We will consider separately the case 
$j=2$. 

Let $j\geq 3.$ By Corollary \ref{pd1}, we have $\projdim_S(S/I)=n-1$ if and only if $x_{i-1}x_n\geq_{lex}v.$ If this is the case, then,
 since $\projdim_S(S/I)\leq \ara(I)\leq n-1,$ we get the required equality. 
 
 Now let $v>_{lex}x_{i-1}x_n$ in the same hypothesis on $j,$ namely $j\geq 3.$ We have $\projdim_S(S/I)=n-2.$ We are going to 
 distinguish two cases to study. In both cases we show that $\ara(I)=n-2=\projdim_S(S/I)$ by using Schmitt--Vogel Lemma.
 
 \textit{Case 1:} Let $i=4$ or $x_{i-1}x_i\geq_{lex}v>_{lex}x_{i-1}x_n.$ In particular, by our assumption $j\geq 3,$  we 
 have  $i\geq 4.$ We display the minimal monomial generators of $I$ in an upper triangular tableau as follows. In the first row we put the generators 
 divisible by $x_2$ ordered decreasingly with respect to the lexicographic order except the monomial $x_2x_n$ which is intercalated
between the monomials $x_2x_{i-1}$ and $x_2x_i.$ In the same way we order on the second row the monomials divisible by $x_3,$ 
intercalating the monomial $x_3x_n$ between $x_3x_{i-1}$ and $x_3x_i.$ We continue in this way up to the row containing the monomials divisible by $x_{i-2}.$ On the next row we put the monomials $x_1x_n,x_1x_i,x_1x_{i+1},\ldots,x_1x_{n-1},$ and, finally, on the last row, we put the remaining generators, namely $x_{i-1}x_i,\ldots,v.$ Then our tableau looks as follows.
\[
\begin{array}{lllllllll}
	x_2x_3 & x_2x_4 & \ldots & x_2x_{i-1} & \underline{x_2x_n} & x_2x_i & \ldots & x_2x_{n-2} & x_2x_{n-1}\\
	       & x_3x_4 & \ldots & x_3x_{i-1} & \underline{x_3x_n} & x_3x_i & \ldots & x_3x_{n-2} & x_3x_{n-1}\\   
	       &        &        & \vdots     &     \vdots         &\vdots  &        &   \vdots   & \vdots     \\
	       &        &        & x_{i-2}x_{i-1}&\underline{x_{i-2}x_n} & x_{i-2}x_i & \ldots & x_{i-2}x_{n-2}& x_{i-2}x_{n-1}\\
	       &        &        &               &\underline{x_1x_n} &  \underline{x_1x_i} & \ldots & \underline{x_1x_{n-2}}&
	       \underline{x_1x_{n-1}}\\
	       &        &        &               &                   & x_{i-1}x_i & \ldots \ v
\end{array}
\]

Next we define the sets $A_1,A_2,\ldots,A_{n-2}$ in the following way. In the first set we put the monomial from the left-up corner 
of the tableau. In the second set we put the two monomials from the left up parallel to the diagonal of the triangular tableau. In the third set we collect the three monomials from the next parallel to the diagonal, and so on. Explicitly, the sets are the following ones.
\[
\begin{array}{ll}
	A_1 = & \{x_2x_{n-1}\},\\
	A_2 = & \{x_2x_{n-2}, x_3x_{n-1}\},\\
	A_3 = & \{x_2x_{n-3},x_{3}x_{n-2},x_4x_{n-1}\}, \\
	\vdots &  \\
	A_{n-i+1}= & \{x_2x_n, x_3x_i,x_4x_{i+1},\ldots \},\\
	\vdots &  \\
	A_{n-2}= & \{x_2x_3,x_3x_4,\ldots, x_{i-2}x_{i-1},x_1x_n,x_{i-1}x_i\}.
\end{array}.
\]
One may easy check that the sets $A_1,\ldots,A_{n-2}$ verify all the conditions from Lemma \ref{svlemma}. We give only a brief 
explanation concerning the third condition. Indeed if one picks up two different monomials in the set $A_j$ for some $j\geq 2,$ let us 
say $m_1$ from the $r-$th row and $m_2$ from the $s-$th row of the tableau with $r<s,$ then the monomial $m^{\prime}$ at the 
intersection of the 
$r-$th row and the column of $m_2$ divides the product $m_1m_2$ and $m^{\prime}\in A_{\ell}$ for some $\ell < j.$

\textit{Case 2:} Let $x_3x_4\geq_{lex} v=x_jx_r>_{lex}x_{i-1}x_i$ and $i\geq 5.$ Then we construct a similar triangular tableau to that one 
from the previous case, but we preserve the decreasing lexicographic order in each row. In this tableau we will add the underlined 
monomials in the $(j-1)-$th row.
\[
\begin{array}{llllllllll}
x_2x_3 & x_2x_4 & \ldots & x_2x_j & x_2x_{j+1} & \ldots & x_2x_r & x_2x_{r+1} & \ldots & x_2x_n\\
       &x_3x_4 & \ldots & x_3x_j & x_3x_{j+1} & \ldots & x_3x_r & x_3x_{r+1} & \ldots & x_3x_n\\
       &       &        & \vdots & \vdots     &        & \vdots &  \vdots    &        &  \vdots\\
       &       &        & x_{j-1}x_j & x_{j-1}x_{j+1} & \ldots & x_{j-1}x_r & x_{j-1}x_{r+1} & \ldots & x_{j-1}x_n\\   
       &       &        &            & x_jx_{j+1}& \ldots & x_jx_r=v & \underline{x_1x_jx_{r+1}}& \ldots & \underline{x_1x_jx_n}\\
       &       &        &            &           & x_1x_i &  &  \ldots  &  & x_1x_n     

\end{array}
\]
Note that in this case it is impossible to have $i=j+1$. Indeed, if $i=j+1,$ then, by our hypothesis we have 
$x_jx_r>_{lex}x_jx_{j+1}$, which is impossible. 

One may easy check that the sets $A_1=\{x_2x_n\}, A_2=\{x_2x_{n-1}, x_3x_n\}, A_3=\{x_2x_{n-2},x_3x_{n-1},x_4x_n\},\ldots, 
A_{n-2}=\{x_2x_3,x_3x_4,\ldots,x_jx_{j+1}\}$ verify the conditions from Lemma \ref{svlemma}, thus $\ara(I)\leq n-2.$ Since we also have 
$\projdim_S(S/I)=n-2,$ we get that $\ara(I)=\projdim_S(S/I).$

To finish the proof, we only have to consider the case $j=2,$ that is, $u=x_1x_i$ and $v=x_2x_r$ for some $i$ and $r$ such that
$v>_{lex}x_{i-1}x_n.$ Note that, in particular, we have $i-1\geq 2,$ that is, $i\geq 3.$

If $i>r$, then it is easily seen that the Taylor resolution of $I$ is minimal. This implies that $\projdim_S(S/I)=\mu(I)$, where 
$\mu(I)$ denotes the number of the minimal monomial generators of $I.$ Therefore, $\ara(I)=\mu(I)=\projdim_S(S/I).$

If $r\geq i,$ we show that $\ara(I)=\projdim_S(S/I)=n-2$ by using again Lemma \ref{svlemma}. In this case we put the generators of 
$I$ in a $2-$row tableau. 
\[
\begin{array}{llllllll}
	       &        & x_1x_i & \ldots & x_1x_r & x_1x_{r+1} & \ldots & x_1x_n\\
x_2x_3   & \ldots & x_2x_i & \ldots & x_2x_r &            &        &
\end{array}
\]
 If $i> 3,$ we add to the second row the monomials $x_1x_2x_{r+1},\ldots,x_1x_2x_n.$
 
 We get the tableau
\[
\begin{array}{llllllll}
	       &        & x_1x_i & \ldots & x_1x_r & x_1x_{r+1} & \ldots & x_1x_n\\
x_2x_3   & \ldots & x_2x_i & \ldots & x_2x_r &  \underline{x_1x_2x_{r+1}}& \ldots & \underline{x_1x_2x_n}\\
\end{array}
\]
and set
\[
A_1=\{x_1x_2x_n\}, A_2=\{x_1x_n,x_1x_2x_{n-1}\},A_3=\{x_1x_{n-1},x_1x_2x_{n-2}\}, \ldots
\]
\[
\ldots, A_{n-r}=\{x_1x_{r+2},x_1x_2x_{r+1}\},
A_{n-r+1}=\{x_1x_{r+1},x_2x_r\},\ldots, A_{n-2}=\{x_2x_3\}.
\]

If $i=3,$ then we add the monomials $x_1x_2x_{r+1},\ldots, x_1x_2x_{n-1}$ to the initial tableau and get
\[
\begin{array}{llllllll}
x_1x_3	 &    x_1x_4     & \ldots & x_1x_r & x_1x_{r+1} & \ldots & x_1x_{n-1} & x_1x_n\\
x_2x_3   &    x_2x_4     & \ldots  & x_2x_r & \underline{x_1x_2x_{r+1}}& \ldots & \underline{x_1x_2x_{n-1}} & \\
\end{array}
\]
We set 
\[
A_1=\{x_1x_3\}, A_2=\{x_1x_4,x_2x_3\}, A_3=\{x_1x_5,x_2x_4\},\ldots, A_{r-2}=\{x_1x_r,x_2x_{r-1}\},
\]
\[
A_{r-1}=\{x_1x_{r+1},x_2x_r\},A_r=\{x_1x_{r+2},x_1x_2x_{r+1}\},\ldots,A_{n-2}=\{x_1x_n,x_1x_2x_{n-1}\}.
\]
In both cases, by using Lemma \ref{svlemma}, we get $\projdim_S(S/I)=n-2\leq \ara(I)\leq n-2,$ hence $\ara(I)=n-2=\projdim_S(S/I).$
\end{proof}

The Cohen--Macaulay lexsegment edge ideals are precisely the ones which are set--theoretic complete intersection, by Theorem \ref{aralei}.

\begin{Corollary1}
Let $I$ be a lexsegment edge ideal. Then the following statements are equivalent:
\begin{itemize}
	\item [(a)] $I$ is Cohen--Macaulay.
	\item [(b)] $I$ is a set--theoretic complete intersection.
\end{itemize}
\end{Corollary1}

Let $I=(L(u,v))$ be the lexsegment edge ideal generated by the set $L(u,v)$ and $I^{\vee}$ be its Alexander dual ideal. Then we have 

\begin{eqnarray*}
I^{\vee} &= &(x_1,x_{i})\cap(x_1,x_{i+1})\dots \cap (x_1,x_{n})
\cap(x_2,x_3)\cap(x_2,x_{4})\cap \dots \cap (x_2,x_{n})\\
&&\cap(x_3,x_4)\cap(x_3,x_{5})\cap \dots \cap (x_3,x_{n})\cap \dots
\cap(x_j,x_{j+1})\cap \dots \cap (x_j,x_{r}),
\end{eqnarray*}
which is an unmixed ideal of height two (see, e.g., \cite[Proposition 1.1]{T1}).
In this section we show the equality 
$\ara(I^{\vee})=\projdim_S(S/I^{\vee})$. 
Since $\projdim_S(S/I^{\vee})=\reg(I)$, we have the following:

\begin{Proposition1}
Let $I=(L(u,v))$ be a lexsegment edge ideal determined by $u=x_1x_i, v=x_jx_r, j\geq 2.$ Then
\[
\projdim_S(S/I^{\vee})=\left\{ 
\begin{array}{ll}
	3, & \text{ if } i\geq j+2 \text{ and } x_n\not| v \\
	2, & { otherwise.}
\end{array}\right.
\]
\end{Proposition1}

Now we determine the arithmetical rank of the Alexander dual of a lexsegment edge ideal. The lexsegment edge ideals enlarge the class of squarefree monomial ideals satisfying the equality between the arithmetical rank and the projective dimension.

\begin{Theorem1}\label{aradual}
Let $I=(L(u,v))$ be a lexsegment edge ideal. Then 
\[
\ara(I^{\vee})=\projdim_S(S/I^{\vee}).
\]
\end{Theorem1}

\begin{proof}
We may assume that $u=x_1x_i, v=x_jx_r$.
If $j=1$, then $I^{\vee}=(x_1,x_{i})\cap(x_1,x_{i+1}) \cap \dots \cap (x_1,x_{r})
 =(x_1,x_{i}x_{i+1}\dots x_{r})$ is a (set--theoretic) complete intersection.
Hence we may assume that $j \geq 2$.

Now we assume that  $i\leq j+1 \text{ or } r=n$.
Then we have  $\projdim_S(S/I^{\vee})=\height I^{\vee}=2$, and  $S/I^{\vee}$ is Cohen--Macaulay.
In this case $I^{\vee}$ is a set--theoretic complete intersection by K. Kimura \cite{K}.
Hence $\ara(I^{\vee})=\projdim_S(S/I^{\vee})=2$.

Next we assume that  $i\geq j+2 \text{ and } r\ne n$. 
Let $J^{\vee}$ be the Alexander dual ideal of $J=(L(x_1x_i, x_{j-1}x_n))$.
Then we have $\ara(J^{\vee})=\projdim_S(S/J^{\vee})=2$.
Hence there exist $f_1,f_2 \in S$ such that
$\sqrt{(f_1,f_2)}=J^{\vee}.$
Then we have
\begin{eqnarray*}
I^{\vee} &= &J^{\vee} \cap(x_j,x_{j+1})\cap(x_j,x_{j+2})\cap \dots \cap (x_j,x_{r})\\
&= &\sqrt{(f_1,f_2)} \cap(x_j,x_{j+1}x_{j+2}\dots x_{r})\\
&= &\sqrt{(x_jf_1, x_jf_2, x_{j+1}x_{j+2}\dots x_{r}f_1, x_{j+1}x_{j+2}\dots x_{r}f_2)}\\
&= &\sqrt{(x_jf_1, x_jf_2+x_{j+1}x_{j+2}\dots x_{r}f_1, x_{j+1}x_{j+2}\dots x_{r}f_2)}.
\end{eqnarray*}
We have $3 =\projdim_S(S/I^{\vee}) \le \ara(I^{\vee}) \le 3$.
Hence $\ara(I^{\vee})=\projdim_S(S/I^{\vee})=3$, as desired.
\end{proof}

\section{Squarefree lexsegment ideals generated in degree $3$}

We begin by computing the depth of a squarefree lexsegment ideal generated in degree $3$. The results of this section are contained in \cite{EKOT}, joint work with V. Ene, K. Kimura and N. Terai. The main results are listed in the following theorem, which gives a complete description of $S/I$ when $I$ is a squarefree lexsegment ideal generated in degree $3$.

\begin{Theorem1}\label{main}
Let $u=x_1x_{i_2}x_{i_3}$ and $v=x_{j_1}x_{j_2}x_{j_3}$ where $j_1\geq 2,$ be two squarefree monomials of degree $3$ and $I=(L(u,v))$ the squarefree 
lexsegment ideal determined by them. Then:
\begin{itemize}
	\item [(a)] $\depth(S/I)=2$ if $x_{i_2-1}x_{i_3-1}x_n\geq_{lex} v.$ 
	\item [(b)] $\depth(S/I)=4$ if $i_2=4$, $i_3\geq 6$ and $j_1=2,\ j_2=3,\  j_3<i_3-1$ or $i_2\geq 5$ and $j_1=2,\ j_2=3,\ i_2-1\leq j_3\leq n-1.$
	\item [(c)] $\depth(S/I)=i_2-j_3+3$ if  $i_2 > 4, j_2=2,$ and $j_3\leq i_2-1.$
	\item [(d)] $\depth(S/I)=3$ in all the other cases. 
\end{itemize}
\end{Theorem1}

To begin with, we  characterize those ideals $I=(L(u,v))$ for which $\depth(S/I)=2$, which is equivalent to $\projdim(S/I)=n-2$. In the sequel we 
denote by $\Delta$ the simplicial complex on the vertex set $[n]$ whose Stanley--Reisner ideal $I_{\Delta}$ is equal to $I.$

\begin{Proposition1}
\label{depth2}
Let $I=(L(u,v))$ be the squarefree lexsegment ideal determined by $u=x_1x_{i_2}x_{i_3}$ and $v=x_{j_1}x_{j_2}x_{j_3}$. Then $\depth(S/I)=2$ if and only if $x_{i_2-1}x_{i_3-1}x_n\geq v.$ 
\end{Proposition1}

\begin{proof} 
By Lemma \ref{CM skeleton}, we get $\depth(S/I)\geq 2$.

Let $x_{i_2-1}x_{i_3-1}x_n\geq_{lex} v.$ We show that the simplicial complex $\Gamma=\Delta^{(2)}$ is not Cohen--Macaulay. Then Lemma \ref{CM skeleton} implies that 
$\depth(S/I)\leq 2,$ which yields the desired equality.

Let us assume that $\Gamma$ is Cohen--Macaulay. Then, by Reisner's criterion, it follows, in particular, that the links of all vertices of $\Gamma$ are Cohen--Macaulay. As all the links are of dimension $1$, it results that all the links are connected. 

We have  $F\in \Gamma$ if and only if $x_F\not\in I_{\Delta}$ and $\deg(x_F)=3.$ Therefore $F\in\Gamma$ if and only if it is of one of 
the following forms:
\begin{itemize}
	\item $F=\{1,r,s\}$ with $1< r< s$ such that $x_r x_s>_{lex} x_{i_2}x_{i_3}.$ We denote by $\mathcal{F}_1$ the set of all facets of  
	$\Gamma$ of this form.
	\item $F=\{p,r,s\}$ with $2\leq p < r < s$ and such that $x_p x_r x_s < v.$ Let $\mathcal{F}_2$ be the set of all facets of $\Gamma$
	of this form.
\end{itemize} 

Since $v\leq_{lex} x_{i_2-1}x_{i_3-1}x_n,$ it follows that $j_1\geq i_2-1.$ We investigate in turn the  cases $j_1 > i_2,$ $j_1=i_2$, 
$j_1=i_2-1$, and show in each of them that our assumption on $\Gamma$ to be Cohen--Macaulay leads to a contradiction.\\

\textit{Case 1:} Let $j_1 > i_2$ and take $\link_{\Gamma}\{j_1\}.$ If $j_1\in F=\{1,r,s\}$ for some $F\in \mathcal{F}_1,$ then $j_1=s$ since 
$r\leq i_2.$ If $j_1\in F=\{p,r,s\}$ for some $F\in \mathcal{F}_2,$ then $j_1=p$ since $x_p x_r x_s < v.$ Therefore, 
$\link_{\Gamma}\{j_1\}$ is generated by the following two disjoint sets of facets: $\{1,r\}$ with $r\leq i_2$ and $\{r,s\}$ with 
$s > r > j_1 >i_2,$ thus $\link_{\Gamma}\{j_1\}$ is disconnected, which is in contradiction to the assumption on $\Gamma$.\\

Now we are going to consider $j_1=i_2$ or $j_1=i_2-1$. Note that in this last situation, we have $j_2 > i_3-1$ by the assumption on 
$v.$\\

\textit{Case 2:} $j_1=i_2$ or $j_1=i_2-1$ and $j_2 > i_3$. In these cases we look at the link of $i_3$ in $\Gamma.$ In the first place we 
observe that if $i_3\in F=\{1,r,s\}$ for some $F\in \mathcal{F}_1,$ we must have $s=i_3$ and $r < i_2.$ In order to show that 
$\link_{\Gamma}\{i_3\}$ is disconnected  and, thus, derive a contradiction, we have to prove that there is no facet in 
$\link_{\Gamma}\{i_3\}$ contained in some facet $F\in \mathcal{F}_2$ which contains $i_2-1.$ Let us assume that there exists some 
$F=\{p,r,s\}\in\mathcal{F}_2$ such that $i_2-1\in F.$ Then we have $p=i_2-1$ and the inequality $x_p x_r x_s < v$ yields 
$r\geq j_2 > i_3,$ that is $i_3\not\in F.$ Therefore, we proved in this case that $\link_{\Gamma}\{i_3\}$ is disconnected. \\

\textit{Case 3:} Let $j_1=i_2-1$, $j_2 = i_3$, and $i_3-1 > i_2.$, that is $v=x_{i_2-1} x_{i_3} x_{j_3}.$ In this case we are going to derive a contradiction 
with the Cohen--Macaulayness of $\Gamma$ in the following way. Since $\Gamma$ is assumed to be Cohen--Macaulay and $\dim \Gamma=2,$ by Reisner's criterion, we have $\tilde{H_1}(\Gamma; k)=0.$ For $i=1,2,$ let $\Gamma_i$ be the simplicial complex generated by 
$\mathcal{F}_i.$ We consider the following part of the Mayer--Vietoris sequence:
\[
\dots\to \tilde{H}_1(\Gamma; k)\to \tilde{H}_0(\Gamma_1\cap \Gamma_2; k) \to \tilde{H}_0(\Gamma_1; k)\oplus 
\tilde{H}_0(\Gamma_2; k)\to \dots.
\]

Since $\Gamma_1$ and $\Gamma_2$ are obviously connected, our assumption on $\Gamma$ implies that $\Gamma_1\cap \Gamma_2$ should 
be connected as well. Let us take $G\in \Gamma_1\cap \Gamma_2.$ Since $1\not\in F$ for all $F\in \mathcal{F}_2,$ we have $G=\{r,s\}$ 
for some $1 < r < s$ such that $x_r x_s > x_{i_2} x_{i_3}$, in particular $r\leq i_2,$ and $\{r,s\}\subset F$ for some $F\in \mathcal{F}_2,$ that is, 
$x_r x_s x_t < x_{i_2-1}x_{i_3} x_{j_3}$ for some $t,$ or  $x_p x_r x_s < x_{i_2-1}x_{i_3} x_{j_3}$ for some $p,$ or 
$x_r x_t x_s < x_{i_2-1}x_{i_3} x_{j_3}$ for some $t.$ Note also that $r\leq i_2$ and $\{r,s\}\subset F$ for some $F\in \mathcal{F}_2$
implies that $r\in \{i_2-1,i_2\}.$ 

Let us take first $r=i_2-1.$ Then, since $\{r,s\}\subset F$ for some $F\in\mathcal{F}_2,$ we have $s\geq i_3,$ hence $G\in \{\{i_2-1,s\}\colon s\geq i_3\}.$ On the other hand, if $r=i_2$, since $x_r x_s > x_{i_2} x_{i_3}$, we get that 
$G=\{r,s\}\in \{\{i_2,i_2+1\},\ldots, \{i_2,i_3-1\}\}$. Consequently, if $i_3-1 > i_2,$ it follows that $\Gamma_1\cap \Gamma_2$ is 
disconnected, which is a contradiction.\\

Finally, we have to consider\\

 \textit{Case 4:} Let $j_1=i_2-1,$ $j_2=i_3$, and $i_3-1=i_2.$ Under these last conditions we get 
$\link_{\Gamma}\{i_2\}$ disconnected. Indeed, if $i_2\in \{1,r,s\}\in \mathcal{F}_1$, we must have $s=i_2,$ which implies $r\leq i_2-1.$
If $i_2\in\{p,r,s\}\in\mathcal{F}_2,$ we must have $i_2=p.$ This implies that the facets of the link of $i_2$ in $\Gamma$ are 
of the form $\{1,r\}$ for some $r\leq i_2-1$ or $\{r,s\}$ for some $s > r > i_2.$ Therefore, $\link_{\Gamma}\{i_2\}$ is disconnected.

For the converse, we show that if 
\begin{equation}\label{eq1}
v>_{lex} x_{i_2-1}x_{i_3-1}x_n,
\end{equation}
then $\Gamma=\Delta^{(2)}$ is a Cohen--Macaulay simplicial complex, whence $\depth(S/I)\geq 3,$ a contradiction. 

Let $\Gamma_i$ be the complex generated by the facets $F\in \mathcal{F}_i$ for $i=1,2.$ We first show that $\Gamma_1\cap \Gamma_2$ is a connected graph. Note that for any facet $\{1,r,s\}\in \mathcal{F}_1$ we have $x_rx_s>_{lex}x_{i_2}x_{i_3},$ hence $r<i_2$ or $r=i_2$ and $s<i_3.$ On the other hand, for $\{p,q,t\}\in \mathcal{F}_2$, we have $t, q> p \geq j_1.$ We will determine all the edges of $\Gamma_1\cap \Gamma_2$. The connectedness of 
$\Gamma_1\cap \Gamma_2$ will become then obvious.

We first look at the edges $\{r,s\}$  with $r=j_1.$ If $v > x_{j_1}x_{j_2}x_n,$ then 
$\{j_1,j_2\},\{j_1,j_2+1\},\ldots,\{j_1,n\}$ are edges in $\Gamma_1\cap \Gamma_2$. Indeed, for all $0\leq t\leq n-1-j_2,$ we have $x_{j_1}x_{j_2+t}>x_{i_2}x_{i_3}$, since, by (\ref{eq1}), $j_1<i_2.$ Therefore, we get $\{1,j_1,j_2+t\}\in \Gamma_1$ for $0\leq t\leq n-1-j_2.$ It is obvious that, in our hypothesis on $v,$ we also have $\{j_1,j_2+t,n\}\in \Gamma_2$. Let us take now $j_1< r \leq i_2-1.$ Then any $2-$set $\{r,s\}$ such that $\{1,r,s\}\in \mathcal{F}_1$ is an edge of $\Gamma_1\cap \Gamma_2$. Indeed, in this case, $x_rx_s>_{lex} x_{i_2}x_{i_3}$ and $v>_{lex} x_rx_sx_t$ for $t>s.$ Therefore, $\{r,r+1\}, \{r,r+2\},\ldots, \{r,n\}$ are edges of $\Gamma_1\cap \Gamma_2$ for all $j_1< r\leq i_2-1.$ Finally, let $r=i_2.$ Then 
$s\leq i_3-1$ in order to have $x_rx_s >_{lex} x_{i_2}x_{i_3}.$ Moreover, $v>_{lex} x_{i_2}x_{s}x_t$ for any $t>s>i_2.$ Consequently, $\{i_2,i_2+1\},\ldots, \{i_2,i_3-1\}$ are edges of $\Gamma_1\cap \Gamma_2.$ Now, one may easily see that $\Gamma_1\cap \Gamma_2$ is a connected graph, thus, by Proposition~\ref{CM connected}, it follows that $\Gamma_1\cap \Gamma_2$ is Cohen--Macaulay. On the other hand, $\Gamma_1$ and $\Gamma_2$ are shellable. Indeed, one may consider the decreasing lexicographic order on the facets of $\Gamma_1$ and the increasing reverse lexicographic order on the facets of $\Gamma_2$ and get in each case the shelling order of the facets. Hence $\Gamma_1$ and $\Gamma_2$ are Cohen--Macaulay. From the Mayer--Vietoris sequence
\[
\dots\to \tilde{H}_1(\Gamma_1; k)\oplus 
\tilde{H}_1(\Gamma_2; k) \to  \tilde{H}_1(\Gamma; k)\to \tilde{H}_0(\Gamma_1\cap \Gamma_2; k) \to  \dots,
\]
by applying Theorem \ref{Reisner}, it follows that $\tilde{H}_1(\Gamma; k)=0$. One may easily check by inspection that the links of all the vertices of $\Gamma$ are connected. By applying again Theorem~\ref{Reisner}, it follows that $\Gamma$ is Cohen--Macaulay.
\end{proof}

In the next proposition we find conditions to have $\depth(S/I)=3.$

\begin{Proposition1}\label{depth3}
If one of the following conditions holds:
\begin{itemize}
	\item [(a)] $x_1x_3x_4\geq_{lex} u=x_1x_{i_2}x_{i_3}\geq_{lex} x_1x_4x_5$ and $v>_{lex} x_{i_2-1}x_{i_3-1}x_n$,
	\item [(b)] $x_1x_4x_6\geq_{lex} u=x_1x_{i_2}x_{i_3}\geq_{lex} x_1x_4x_n$ and $x_2x_3x_{i_3-1}\geq v>_{lex}x_3 x_{i_3-1}x_n$,
	\item [(c)] $x_1x_5x_6\geq_{lex} u=x_1x_{i_2}x_{i_3}$ and $x_2x_3x_n\geq_{lex} v>_{lex}x_{i_2-1}x_{i_3-1}x_n$,
\end{itemize}
then $\depth(S/I)=3$. 
\end{Proposition1}

\begin{proof}
By the previous proposition, we have $\depth(S/I)\geq 3.$ We have to show that $\depth(S/I)\leq 3.$ For this purpose, in most cases, we find a facet of 
$\Delta$ of cardinality $3.$ This implies that there exists an associated prime ideal of $I$ of dimension $n-3.$ By Proposition \ref{depth-bound}, we get the desired conclusion. 

(a) Let $F=\{1,2,4\}.$ We claim that $F\in \mathcal{F}(\Delta).$  Indeed, it is clear that $F\in \Delta$. Next, we observe that $x_3x_F\in I$ since
$x_2x_3x_4\in I$ and $x_2x_3x_4| x_3x_F,$ and, for any $s\geq 5$ we have $x_1x_4x_s|x_sx_F,$ thus $x_sx_F\in I$ since $x_1x_4x_s\in I.$

(b) For $u=x_1x_4x_t$ for some $t\geq 6$ and $x_2x_3x_{t-1}>_{lex}v>_{lex} x_3x_{t-1}x_n$, one easily checks, as above, that $F=\{1,3,t\}\in \mathcal{F}(\Delta).$ Now, let $u=x_1x_4x_t$ for some $t\geq 6$ and $v=x_2x_3x_{t-1}.$ In this case, the simplicial complex $\Delta$ is 
	\[
	\Delta=\left\langle\begin{array}{ccc}
	\{1,2,3,t\} & \ldots & \{1,2,3,n\}\\
	\{1,3,4,5\} & \ldots &\{1,3,4,t-1\}\\
	\{2,3,t,\ldots,n\}&\{2,4,\ldots,n\}&\{3,4,\ldots,n\}\\
	\end{array}
	\right\rangle.
\]
It follows that 

	\[
\link_{\Delta}\{1,3\}=\left\langle\begin{array}{ccc}
	\{2,t\} & \ldots & \{2,n\}\\
	\{4,5\} & \ldots &\{4,t-1\}\\
	\end{array}
	\right\rangle
\]
thus it is disconnected. By Theorem~\ref{Reisner}, $\Delta$ it is not Cohen--Macaulay, thus $\depth(S/I)\leq 3$, by Lemma \ref{CM skeleton}.

(c) In this last case it is obvious that $F=\{1,2,3\}\in \mathcal{F}(\Delta).$  
\end{proof}

In order to get the characterization of the depth of $S/I$ in the next steps, we need to use the inequality $\projdim(S/I)\leq \ara(I)$. As a byproduct we obtain the arithmetical rank of $I.$

\begin{Proposition1}
Set $u=x_1x_{i_2}x_{i_3}$ and $v=x_2x_3x_{j_3}$. Let $I=(L(u,v))$ be the squarefree lexsegment ideal of $L(u,v)$ in $S$. If either of the conditions $(1)$ $i_2=4,\ i_3\geq 6,\ j_3<i_3-1$ or $(2)$ $i_2\geq 5,\ i_2-1\leq j_3\leq n-1$ holds, then $\ara(I)\leq n-4$. 
\end{Proposition1}

\begin{proof}
\textit{Case 1:} $i_2=4,\ i_3\geq 6,\ j_3<i_3-1$.

Consider the following tableau. Set $i=i_3$ and $j=j_3$.

{\tiny{$$\begin{array}{lllllllllll}
x_1x_5x_6 & x_1x_5x_7 &  \ldots & x_1x_5x_{j+1} & x_1x_5x_{j+2} & \ldots & x_1x_5x_{i-1} & \underline{x_2x_3x_5} & x_1x_5x_{i} & \ldots & x_1x_5x_{n} \\
          & x_1x_6x_7 &  \ldots & x_1x_6x_{j+1} & x_1x_6x_{j+2} & \ldots & x_1x_6x_{i-1} & \underline{x_2x_3x_6} & x_1x_6x_{i} & \ldots & x_1x_6x_{n} \\                     &           &         & \vdots        & \vdots        &        & \vdots        & \vdots                & \vdots      &        & \vdots      \\ 
          &           &         & x_1x_jx_{j+1} & x_1x_jx_{j+2} & \ldots & x_1x_jx_{i-1} & \underline{x_2x_3x_j} & x_1x_jx_{i}  & \ldots & x_1x_jx_{n} \\  
          &           &         &               & x_1x_{j+1}x_{j+2} &\ldots & x_1x_{j+1}x_{i-1} & \underline{x_2x_3x_4x_{j+1}} & x_1x_{j+1}x_{i}  & \ldots & x_1x_{j+1}x_{n} \\ 
          &           &         &               &                   &       & \vdots & \vdots & \vdots  &  & \vdots \\ 
          &           &         &               &                   &       & x_1x_{i-2}x_{i-1} & \underline{x_2x_3x_4x_{i-1}} & x_1x_{i-2}x_{i}  & \ldots & x_1x_{i-2}x_{n} \\ 
          &           &         &               &                   &       &                   & \underline{x_2x_3x_4} & \underline{x_1x_4x_{i}}  & \ldots & \underline{x_1x_4x_{n}} \\ 
          &           &         &               &                   &       &                   &                       & x_1x_{i-1}x_{i}  & \ldots & x_1x_{i-1}x_{n} \\ 
          &           &         &               &                   &       &                   &                       &                  &        & \vdots \\ 
          &           &         &               &                   &       &                   &                       &                  &        & x_1x_{n-1}x_n 
\end{array}$$}}
Next we define the sets $A_1,A_2,\ldots,A_{n-4}$ as follows: 
$$\begin{array}{ll}
A_1= & \{x_1x_5x_n\},\\
A_2= & \{x_1x_5x_{n-1},x_1x_6x_n\},\\           
A_3= & \{x_1x_5x_{n-2},x_1x_6x_{n-1},x_1x_7x_n\},\\          
\vdots &  \\ 
A_{n-4}= & \{x_1x_5x_6,x_1x_6x_7,\ldots,x_1x_{i-2}x_{i-1},x_2x_3x_4,x_1x_{i-1}x_i,\ldots,x_1x_{n-1}x_n\}. 
\end{array}$$

Then $A_1,A_2,\ldots,A_{n-4}$ satisfy the conditions of the Schmitt--Vogel Lemma and the desired inequality follows.
\medskip

\textit{Case 2:} $i_2\geq 5,\ i_2-1\leq j_3\leq n-1$.

Consider the following tableau. Set $h=i_2$, $i=i_3$ and $j=j_3$.

$$\begin{array}{lllllll}
                  & x_1x_hx_i & \ldots\ldots & \ldots\ldots &  \ldots & x_1x_{h}x_{n} & \underline{x_2x_3x_{h}} \\
x_1x_{h+1}x_{h+2} & \ldots & x_1x_{h+1}x_{j+1} & x_1x_{h+1}x_{j+2} &  \ldots & x_1x_{h+1}x_{n} & \underline{x_2x_3x_{h+1}} \\
                  &        & \vdots            & \vdots            &         & \vdots          & \vdots  \\
                  &        & x_1x_{j}x_{j+1} & x_1x_{j}x_{j+2} &  \ldots & x_1x_{j}x_{n} & \underline{x_2x_3x_{j}} \\
                  &        &                 & x_1x_{j+1}x_{j+2} &  \ldots & x_1x_{j+1}x_{n} & \underline{x_2x_3x_4x_{j+1}} \\
                  &        &                 &                   &         & \vdots & \vdots \\       
                  &        &                 &                   &         & x_1x_{n-1}x_{n} & \underline{x_2x_3x_4x_{n-1}} \\
                  &        &                 &                   &         &   & \underline{x_2x_3x_4} \\
                  &        &                 &                   &         &   & \underline{x_2x_3x_5} \\
                  &        &                 &                   &         &   & \vdots \\
                  &        &                 &                   &         &   & \underline{x_2x_3x_{h-1}}
\end{array}	$$

Next we define the sets $A_1,A_2,\ldots,A_{n-4}$ as follows: 

$$\begin{array}{ll}
A_1= & \{x_2x_3x_h\},\\
A_2= & \{x_1x_hx_{n},x_2x_3x_{h+1}\},\\           
A_3= & \{x_1x_hx_{n-1},x_1x_{h+1}x_{n},x_2x_3x_{h+2}\},\\          
\vdots &  \\ 
A_{n-4}= & \{x_2x_3x_{h-1}\}. 
\end{array}$$

Then $A_1,A_2,\ldots,A_{n-4}$ satisfy the conditions of the Schmitt--Vogel Lemma and the desired inequality follows.
\end{proof}

\begin{Proposition1}\label{depth4}
If one of the following conditions holds:
\begin{itemize}
	\item [(a)] $x_1x_4x_6\geq u=x_1x_4x_{i_3}\geq x_1x_4x_n$ and $v> x_2x_3x_{i_3}$,
	\item [(b)] $x_1x_5x_6\geq u=x_1x_{i_2}x_{i_3}\geq x_1x_{n-1}x_n$ and $x_2x_3x_{i_2-1}\geq v\geq x_2x_3x_{n-1},$
\end{itemize}
then $\depth(S/I)=4$.
\end{Proposition1}

\begin{proof}
By the above proposition, it remains to prove the inequality $\depth(S/I)\leq 4$. We are going to show this by similar arguments to those used in the proof of Proposition~\ref{depth3}. Assume that (a) holds and take $F=\{1,2,3,i_3\}.$ One easily observes that $F\in \mathcal{F}(\Delta)$. In case
(b), we take $F=\{1,2,3,n\}.$ Then, if $4\leq s\leq j_3,$ we get $u> x_2x_3x_s \geq v,$ hence $F\cup \{s\}\not\in \Delta.$ If $s> j_3,$ thus $s\geq i_2, $ we get $u\geq x_1x_sx_n>v,$ whence $F\cup\{s\}\not\in\Delta.$
\end{proof}

\begin{Corollary1}
If one of the following conditions holds:
\begin{itemize}
	\item [(a)] $x_1x_4x_6\geq u=x_1x_4x_{i_3}\geq x_1x_4x_n$ and $v> x_2x_3x_{i_3}$,
	\item [(b)] $x_1x_5x_6\geq u=x_1x_{i_2}x_{i_3}\geq x_1x_{n-1}x_n$ and $x_2x_3x_{i_2-1}\geq v\geq x_2x_3x_{n-1},$
\end{itemize}
then $\ara(I)=n-4$.
\end{Corollary1}

Finally, to get a full picture of the values of $\depth(S/I)$, we have to consider $u=x_1x_{i_2}x_{i_3}, v=x_2x_3x_{j_3}$, and $j_3\leq i_2-1.$ To begin with, we give an upper bound for $\ara(I)$.

\begin{Proposition1}\label{aralast}
Let $u=x_1x_{i_2}x_{i_3}, v=x_2x_3x_{j_3}$, with $j_3\leq i_2-1.$ Then 
\[
\ara(I)\leq n-(i_2-j_3+3).
\]
\end{Proposition1}

\begin{proof}
We consider the following tableau:

$\footnotesize{
\begin{array}{llllllll}
\underline{x_1x_{i_2}x_{i_2+1}x_n} & \underline{x_1x_{i_2}x_{i_2+2}x_n} & \ldots & \underline{x_1x_{i_2}x_{i_3-1}x_n }& x_1x_{i_2}x_{i_3} &  \ldots & x_1x_{i_2}x_{n-1} &  x_1x_{i_2}x_{n}\\
	                   & x_1x_{i_2+1}x_{i_2+2} & \ldots & x_1x_{i_2+1}x_{i_3-1} & x_1x_{i_2+1}x_{i_3} & \ldots & x_1x_{i_2+1}x_{n-1}  & x_1x_{i_2+1}x_{n}\\
            	       &        &        & \vdots     & \vdots                 &  &\vdots      &   \vdots      \\
            	       &        &        &  x_1x_{i_3-2}x_{i_3-1}  & x_1x_{i_3-2}x_{i_3} &  \ldots & x_1x_{i_3-1}x_{n-1}  & x_1x_{i_3-1}x_{n}\\
	                   &        &        &   & x_1x_{i_3-1}x_{i_3} &  \ldots & x_1x_{i_3-1}x_{n-1}  & x_1x_{i_3-1}x_{n}\\
            	       &        &        &            &                  &  &        \vdots       &   \vdots      \\
	               	                   &        &            &                  &  &               &   x_1x_{n-2}x_{n-1} &x_1x_{n-2}x_{n}\\
	                   &        &            &                  &  &               &                   & x_1x_{n-1}x_{n}\\
\end{array}}
$

We define the sets:
\medskip

$\small{
\begin{array}{ll}
	A_1 = & \{x_1x_{i_2}x_{n}\},\\
	A_2 = & \{x_1x_{i_2}x_{n-1}, x_1x_{i_2+1}x_{n}\},\\
	A_3 = & \{x_1x_{i_2}x_{n-2},x_1x_{i_2+1}x_{n-1}, x_1x_{i_2+2}x_{n}\}, \\
	\vdots &  \\
	A_{n-i_2}= & \{x_1x_{i_2}x_{i_2+1}x_n,x_1x_{i_2+1}x_{i_2+2},\ldots,x_1x_{i_3-1}x_{i_3},\ldots,x_1x_{n-2}x_{n-1},x_1x_{n-1}x_{n}\}, \\
	A_{n-i_2+1} = & \{x_2x_3x_{j_3}\},\\
	A_{n-i_2+2} = & \{x_2x_3x_{j_3-1}\},\\
	\vdots &  \\
	A_{n-i_2+j_3-3}= & \{x_2x_3x_4\}.
\end{array}}
$\\

One can easy prove that the sets $A_1,\ldots A_{n-i_2},A_{n-i_2+1},\ldots A_{n-i_2+j_3-3}$ verify the conditions from Schmitt--Vogel Lemma. This implies that $\ara(I)\leq n-i_2+j_3-3$, thus $\projdim(S/I)\leq n-(i_2-j_3+3)$.
\end{proof}

\begin{Proposition1}\label{depthlast}
Let $u=x_1x_{i_2}x_{i_3}, v=x_2x_3x_{j_3}$, with $j_3\leq i_2-1.$ Then $\depth(S/I)=i_2-j_3+3.$
\end{Proposition1}

\begin{proof}
By Proposition~\ref{depthlast} and the inequality $\projdim(S/I)\leq \ara(I)$, we have $\depth(S/I)\geq i_2-j_3+3.$ It remains to prove the other inequality for depth. To this purpose we show that there exists a facet $F\in \mathcal{F}(\Delta)$ of cardinality $i_2-j_3+3.$ We need to consider several cases.

For $i_2=i_3-1,$ we take $F=\{1,\ldots,i_2\}\setminus\{4,\ldots,j_3\}.$ Then $x_F\not\in I$ since every divisor of $x_F$ of degree $3$ is either larger
than $u$ or smaller than $v.$ Now let $s\in [n]\setminus F.$ If $4\leq s\leq j_3$, then $x_2x_3x_s| x_sx_F$ and $u>_{lex}x_2x_3x_s\geq_{lex}v$, thus $F\cup\{s\}\not\in \Delta$. If $s>i_2,$ then $x_sx_F$ is divisible by $x_1x_{i_2}x_s\in I,$ hence $F\cup\{s\}\not\in\Delta.$

For $i_2<i_3-1$ and $j_3<i_2-1,$ we take $F=\{1,2,\ldots,i_2-1,i_3\}\setminus\{4,\ldots,j_3\}$ and claim that $F\in \mathcal{F}(\Delta).$ Indeed,
$x_F\not\in I$ and for any $s\in [n]\setminus F$ we obtain $F\cup\{s\}\not\in \Delta$ by inspecting each of the following sub cases.

(a) For $4\leq s\leq j_3$, we have $x_2x_3x_s\in I$ and $x_2x_3x_s|x_sx_F.$

(b) For $i_2\leq s\leq i_3-1,$ we have $x_1x_sx_{i_3}\in I$ and $x_1x_sx_{i_3}|x_sx_F.$

(c) For $s> i_3$ we have $x_1x_{i_3}x_s\in I$ and $x_1x_{i_3}x_s| x_sx_F.$

Finally, let us consider $i_2<i_3-1$ and $j_3=i_2-1,$ that is, $i_2-j_3+3=4.$ Then we take $F=\{1,2,3,i_3\}$. It is obvious that $x_F\not\in I.$ We 
show that $F\cup\{s\}\not\in\Delta$ for all $s\in [n]\setminus F.$ Indeed, we have:

(i) If $4\leq s\leq j_3$ then $x_2x_3x_s\in I$ and $x_2x_3x_s| x_sx_F.$

(ii) If $j_3<s\leq i_3-1,$ then $x_1x_sx_{i_3}\in I$ and $x_2x_3x_s|x_sx_F.$

(iii) If $s>i_3,$ then $x_1x_{i_3}x_s\in I$ and $x_1x_{i_3}x_s|x_sx_F.$
\end{proof}

The computation of the depth in this case presents already an additional difficulty compared with the $2-$degree case. 

The regularity of this class of ideals is not yet computed. By computer calculations, we can conjecture that the regularity of a squarefree lexsegment ideal generated in degree $3$ reaches each of the values $3$, $4$ and $5$.

\section{Some open questions}
In this section, we discuss possible questions concerning the depth and the regularity of an arbitrary squarefree lesegment ideal.

Let $I=(L(u,v))\subset S$ be a squarefree lexsegment ideal determined by $u,v\in\Mon_q^s(S)$, $u>_{lex}v$. As before, we consider $\Delta$ to be the simplicial complex whose Stanley--Reisner ideal is equal to $I$. We already noticed that for the depth of squarefree lexsegment ideals there are known only upper and lower bounds. Some questions concerning this invariant arise.
 
By using the results obtained for squarefree lexsegment ideals generated in small degrees, the following questions and conjectures might have positive answers. 

\begin{Conjecture1}
Let $u=x_1x_{i_2}\cdots x_{i_q}, v=x_{j_1}\cdots x_{j_q}$ with the property that $j_1\geq 2$ and $I=(L(u,v))$ the squarefree lexsegment ideal. Then $\depth(S/I)=q-1$ if and only if $x_{i_2-1}x_{i_3-1}\cdots x_{i_{q-1}}x_n\geq v$.
\end{Conjecture1}

\begin{Question1}
Let  $q-1\leq d\leq n-2$ be an integer. Are there squarefree monomials $u,v\in \Mon_q^s(S)$ such that $\depth(S/I)=d?$
\end{Question1}

The initial and final squarefree lexsegment ideals generated in degree $q$ have a linear resolution. Therefore, these ideals have regularity $q.$ A possible approach in the computation of the Castelnuovo--Mumford regularity might be the following one. We may write $I=J+L$ where $J=(L(u,x_1x_{n-q+2}\cdots x_n))$ and $L=(L(x_2x_3\cdots x_{q+1},v))$. Obviously,  $J$ and $L$ have $q-$linear resolutions. It is known \cite{T} that $\reg(J+L)\leq \reg(J) + \reg(L) -1$. Therefore, we have the inequalities $q\leq \reg(I) \leq 2q-1.$ 

By using the exact sequence 
$$0\to J\cap L \to J\oplus L \to I\to 0,$$ one may easily see that $\reg(I) =\reg(J\cap L)-1.$ Unfortunately, to compute the regularity of $J\cap L$ is still difficult. However, one may hope that the following question has a positive answer taking into account computer calculations in several particular cases.

\begin{Question1}
Let $q\leq r\leq 2q-1$ be an arbitrary integer. Are there squarefree monomials $u,v\in \Mon_q^s(S)$ such that $\reg((L(u,v)))=r?$
\end{Question1}

\chapter{Gotzmann lexsegment ideals}
In this chapter, we study the lexsegment ideals in the non--squarefree case. We characterize the componentwise lexsegment ideals which are componentwise linear and the lexsegment ideals generated in one degree which are Gotzmann. The results of this chapter are contained in our paper \cite{OOS}.
\section{Componentwise lexsegment ideals}

We define the componentwise lexsegment ideals and we characterize all the componentwise lexsegment ideals which are componentwise linear.

\begin{Definition1}\label{compwiselex}\rm Let $I$ be a monomial ideal in $S$ and $d$ the least degree of the minimal monomial generators. The ideal $I$ is called \textit{componentwise lexsegment} if, for all $j\geq d$, its degree $j$ component $I_j$ is generated, as $k-$vector space, by the lexsegment set $\mathcal L(x_1^{j-d}u,vx_n^{j-d})$.
\end{Definition1}

Obviously, completely lexsegment ideals are componentwise lexsegment ideals as well.

\begin{Example1}\rm  The ideal $I=(x_1x_3^2,x_2^3,x_1x_2^2x_3)$ is a componentwise lexsegment ideal. Indeed, one may note that $I_3$ is the $k$-vector space spanned by $\mathcal L(x_1x_3^2,x_2^3)$ and $I_4$ is the $k-$vector space generated by $\mathcal L(x_1^2x_3^2,x_2^3x_3)$. Since $\mathcal L(x_1^2x_3^2,x_2^3x_3)$ is a completely lexsegment set, by the characterization of completely lexsegment ideals \cite{DH}, $I_j$ is generated by the lexsegment set $\mathcal L(x_1^{j-2}x_3^2,x_2^3x_3^{j-3})$ for all $j\geq 4$.
\end{Example1}

We characterize all the componentwise lexsegment ideals which are componentwise linear. In the same time, we prove the equivalence of the notions componentwise linear ideal and componentwise linear quotients for this particular class of graded ideals.
 
One may note that we can assume $x_1\mid u$ since otherwise we can study the ideal in a polynomial ring in a smaller number of variables.
 
\begin{Theorem1} Let $I$ be a componentwise lexsegment ideal and $d\geq 1$ the lowest degree of the minimal monomial generators of $I$. Let $u,v\in\Mon_d(S)$, $x_1|u$ be such that $I_{\langle d\rangle}=(\mathcal L(u,v))$. The following conditions are equivalent:

\begin{itemize}
	\item[(a)] $I$ is a componentwise linear ideal.
	\item[(b)] $I_{\langle d\rangle}$ has a linear resolution.
	\item[(c)] $I_{\langle d\rangle}$ has linear quotients.
	\item[(d)] $I$ has componentwise linear quotients.
\end{itemize}
\end{Theorem1}

\begin{proof} (a)$\Rightarrow$(b) Since $I$ is componentwise linear, the statement is straightforward.

(b)$\Rightarrow$(c) This is straightforward by Theorem \ref{lin res lin quot lex}.

(c)$\Rightarrow$(d) We separately treat the case of completely and lexsegment ideals and lexsegment ideals which are not complete. Firstly, let us assume that $I_{\langle d\rangle}$ is a completely lexsegment ideal with linear quotients. Hence $I=I_{\langle d\rangle}$ and $I$ has componentwise linear quotients (Theorem \ref{lin-quot-comp}).

If $I_{\langle d\rangle}=(\mathcal L(u,v))$ is a lexsegment ideal which is not complete, with linear quotients, then $I_{\langle d\rangle}$ has a linear resolution and, by Theorem \ref{noncompletelylex}, $u$ and $v$ must have the form
	\[u=x_1x_{l+1}^{a_{l+1}}\cdots x_n^{a_n}\ \mbox{and}\ v=x_lx_n^{d-1}
\]
for some $l$, $2\leq l<n$. Therefore, $\nu_1(u)=1$ and $\nu_1(v)=0$. Here, for a monomial $m=x_1^{a_1}\cdots x_n^{a_n}$, we denoted by $\nu_i(m)$ the exponent of the variable $x_i$, that is $\nu_i(m)=a_i$.

If we look at the ends of the lexsegment $\mathcal L(x_1u,vx_n)$, we have $\nu_1(x_1u)=2$, $\nu_1(vx_n)=0$ and one may easily see that $(\mathcal L(x_1u,vx_n))$ is a completely lexsegment ideal. By Theorem \ref{completelylex}, $(\mathcal L(x_1u,vx_n))$ has a linear resolution and, using Theorem \ref{lin res lin quot lex}, $(\mathcal L(x_1u,vx_n))$ has linear quotients. Since $(\mathcal L(x_1u,vx_n))$ is a completely lexsegment ideal with a linear resolution, the ideals generated by the shadows of $\mathcal L(x_1u,vx_n)$ are completely lexsegment ideals with linear resolutions, hence they have linear quotients by Theorem \ref{lin res lin quot lex}. Therefore, $I$ has componentwise linear quotients.  

(d)$\Rightarrow$(a) Since any ideal with linear quotients generated in one degree has a linear resolution, the statement follows by comparing the definitions.
\end{proof}

\section{Gotzmann completely lexsegment ideals}

In this section we are going to characterize the completely lexsegment ideals generated in degree $d$ which are Gotzmann.

\begin{Theorem1}\label{compG} Let $u,v\in \Mon_d(S)$, $x_1\mid u$ such that $I=(\mathcal L(u,v))$ is a completely lexsegment ideal of $S$ which is not an initial lexsegment ideal. Let $j$ be the exponent of the variable $x_n$ in $v$ and $a=|\Mon_d(S)\setminus \mathcal L^i(u)|$. The following statements are equivalent:
\begin{itemize}
	\item [(a)] $I$ is a Gotzmann ideal.
	\item [(b)] $a\geq{n+d-1\choose d}-(j+1)$.
\end{itemize}
\end{Theorem1}

\begin{proof} Let $b=|\Mon_d(S)\setminus \mathcal L^i(v)|$ and $w\in \Mon_d(S)$ such that $|\mathcal L(u,v)|=|\mathcal L^i(w)|$. We denote $c=|\Mon_d(S)\setminus \mathcal L^i(w)|$. Then $|\mathcal L^i(w)|=|\mathcal L^i(v)|-|\mathcal L^i(u)|+1=a-b+1$, which yields:
	\[{n+d-1\choose d}-c=a-b+1,
\]
that is
	\begin{eqnarray}c={n+d-1 \choose d}-(a+1)+b. \label{1}
\end{eqnarray}
Since $I$ is completely, $I$ is Gotzmann if and only if 
	\begin{eqnarray}|\mathcal L^i(wx_n)|=|\mathcal L(ux_1,vx_n)|=|\mathcal L^i(vx_n)|-|\mathcal L^i(ux_1)|+1.\label{2}
\end{eqnarray}
Since $x_1\mid u$, we have $|\mathcal L^i(ux_1)|=|\mathcal L^i(u)|$. Therefore, the equality (\ref{2}) is equivalent to
	\[|\mathcal L^i(wx_n)|=|\mathcal L^i(vx_n)|-|\mathcal L^i(u)|+1
\]
that is
	\[|\Mon_{d+1}(S)|-c^{\langle d\rangle}=|\Mon_{d+1}(S)|-b^{\langle d\rangle}-(|\Mon_d(S)|-a)+1.
\]
Here we used the well known formula given in Proposition \ref{cardinality mon-shad}.
Hence $I$ is Gotzmann if and only if
	\begin{eqnarray}c^{\langle d\rangle}=b^{\langle d\rangle}+{n+d-1\choose d}-a-1.\label{3}
\end{eqnarray}
By using (\ref{1}), we obtain
	\[c^{\langle d\rangle}=b^{\langle d\rangle}+c-b,
\]
that is
	\[c^{\langle d\rangle}-c=b^{\langle d\rangle}-b,
\]
which is equivalent to
\begin{eqnarray} c^{(d)}=b^{(d)}.\label{4}
\end{eqnarray}
Let us firstly consider the case $b=0$, that is $v=x_n^d$ and $I$ is the final lexsegment determined by $u$. The equation (\ref{4}) becomes
\begin{eqnarray} c^{(d)}=0.\label{5}
\end{eqnarray}
By Lemma \ref{c^d=0 dnd c<egal b}, $c^{(d)}=0$ if and only if $c\leq d$. 

For the case $b>0$, the monomial $v$ has the form 
	\[v=x_{l_1}\cdots x_{l_{d-j}}x_n^j,
\]
for some $j\geq 0$ and $1\leq l_1\leq\ldots\leq l_{d-j}\leq n-1$. The $d$-binomial expansion of $b$ is
	\[b={n-{l_1}+d-1\choose d}+\ldots+{n-l_{d-j}+j\choose j+1}.
\]
By Lemma \ref{b(d)=c(d)}, the equality (\ref{4}) holds if and only if $j\geq 1$ and $c-b\leq j$. Then we have obtained $c-b\leq j$ for any $b$. By (\ref{1}), this inequality holds if and only if ${n+d-1\choose d}-(a+1)\leq j$, that is
	\[a\geq{n+d-1\choose d}-(j+1).
\]
\end{proof}
\section{Gotzmann lexsegment ideals which are not complete}
The characterization of Gotzmann lexsegment ideals which are not complete involves the Taylor complex. Thus, we recall the construction.

In general, it is not easy to obtain the minimal graded free resolution for any monomial ideal. However, there is a nice non-minimal resolution that was constructed by Diana Taylor in 1960. It generalizes the Koszul complex in a natural way, and nowadays it is called the Taylor complex. 

Let $I$ be a monomial ideal of $S=k[x_1,\ldots,x_n]$ with the minimal monomial generating set $G(I)=\{u_1,\ldots, u_r\}$. The Taylor resolution $(T_{\bullet}(I),d_{\bullet})$ of $I$ is defined as follows. Let $L$ be the free $S-$module with the basis $\{e_1,\ldots,e_r\}$. Then $T_q(I)=\bigwedge\limits^{q+1}L$ for $0\leq q\leq r-1$ and $d_q:T_q(I)\rightarrow T_{q-1}(I)$ for $1\leq q\leq r-1$ is defined as follows
	\[d_q(e_{i_0}\wedge\ldots\wedge e_{i_q})=\sum_{s=0}^q(-1)^s\frac{\lcm(u_{i_0},\ldots,u_{i_q})}{\lcm(u_{i_0},\ldots ,\check{u}_{i_s},\ldots,u_{i_q})}e_{i_0}\wedge\ldots\wedge\check{e}_{i_s}\wedge\ldots\wedge e_{i_q}. 
\]
The augmentation $\varepsilon:T_0\rightarrow I$ is defined by $\varepsilon(e_i)=u_i$ for all $1\leq i\leq q$.

It is known that, in general, the Taylor resolution is not minimal. However, M. Okudaira and Y. Takayama characterized all the monomial ideals with linear resolutions whose Taylor resolutions are minimal \cite{OT}.

\begin{Theorem1}\label{ot}\cite{OT}
Let $I$ be a monomial ideal with a linear resolution. The following conditions are equivalent:
\begin{itemize}
	\item[(i)] The Taylor resolution of $I$ is minimal.
	\item[(ii)] $I=m\cdot (x_{i_1},\ldots,x_{i_l})$ for some $1\leq {i_1}<\ldots<i_l\leq n$ and for a monomial $m$. 
\end{itemize}
\end{Theorem1}

In \cite{HHMT}, the componentwise linear ideals whose Taylor resolutions are minimal are described.
\begin{Theorem1}\label{minTay}\cite{HHMT} Let $I$ be a componentwise linear monomial ideal of $S$. The following conditions are equivalent:
\begin{itemize}
	\item[(i)] The Taylor resolution of $I$ is minimal.
	\item[(ii)] $\max\{m(u):u\in G(I)\}=|G(I)|$.
	\item[(iii)] $I$ is a Gotzmann ideal with $|G(I)|\leq n$.
\end{itemize}
\end{Theorem1}

Now we can complete the characterization of Gotzmann lexsegment ideals which are not complete.

\begin{Theorem1} Let $u=x_t^{a_t}\cdots x_{n}^{a_n}$, $v=x_t^{b_t}\cdots x_{n}^{b_n}$ be two monomials of degree $d$, $u>_{lex} v$, $a_t\neq0$, $t\geq1$ and $I=(\mathcal L(u,v))$ a lexsegment ideal which is not complete. Then $I$ is a Gotzmann ideal in $S$ if and only if $I=m(x_{l},x_{l+1},\ldots,x_{l+p})$ for some $t\leq l\leq n$, some $1\leq p\leq n-l$ and a monomial $m$.
\end{Theorem1}

\begin{proof} If $I=m(x_{l},x_{l+1},\ldots,x_{l+p})$ for some $t\leq l\leq n$, some $1\leq p\leq n-l$ and a monomial $m$, then $I$ is isomorphic to the monomial prime ideal $(x_{l},x_{l+1},\ldots,x_{l+p})$ and the Koszul complex of the sequence $x_{l},x_{l+1},\ldots,x_{l+p}$ is isomorphic to the minimal graded free resolution of $I$. Therefore $I$ has a linear resolution and, by Theorem \ref{ot}, the Taylor resolution of $I$ is minimal. Since any ideal with a linear resolution is componentwise linear, it follows by Theorem \ref{minTay} that $I$ is a Gotzmann ideal. 

Now it remains to prove that, if $I$ is a Gotzmann ideal in $S$, then $I$ has the required form.

Firstly, we prove that $\projdim(S/I)<n$. For this, we study the following cases:

\textit{Case 1: $t=1$, $b_1=0,a_1=1$.} Since $I$ is a lexsegment ideal, which is not complete, but is Gotzmann, $I$ has a linear resolution. Therefore, by Theorem \ref{noncompletelylex}, $u$ and $v$ have the form
	\[u=x_1x_{l+1}^{a_{l+1}}\ldots x_n^{a_n}\ \mbox{and}\ v=x_lx_n^{d-1}
\]
for some $l$, $2\leq l\leq n-1$. Since $x_nu<_{lex}x_1v$, using Proposition \ref{depthzero} we get $\depth(S/I)\neq0$. Hence $\projdim(S/I)<n$.

\textit{Case 2: $t=1$, $0<b_1<a_1$.}  Since $I$ is a lexsegment ideal which is not complete, we must have $b_1=a_1-1$. Now, if $I$ does not have a linear resolution, $I$ is not Gotzmann. The ideal $I$ has a linear resolution if and only if $J=(\mathcal L(u',v'))$ has a linear resolution, where $u'=u/x_1^{b_1}$ and $v'=v/x_1^{b_1}$. One may easy check that $J$ is a lexsegment ideal which is not complete. Therefore $J$ has a linear resolution if and only if $u'$ and $v'$ have the form
	$u'=x_1x_{l+1}^{a_{l+1}}\ldots x_n^{a_n}\ \mbox{and}\ v'=x_lx_n^{d-1}$
for some $l$, $2\leq l\leq n-1$ by Theorem \ref{noncompletelylex} and this implies $u=x_1^{a_1}x_{l+1}^{a_{l+1}}\ldots x_n^{a_n}\ \mbox{and}\ v'=x_1^{a_1-1}x_lx_n^{d-1}$.
Since $x_nu<_{lex}x_1v$, using Proposition \ref{depthzero} we get $\depth(S/I)\neq0$. Hence $\projdim(S/I)<n$. 

\textit{Case 3: $t=1$, $a_1=b_1>0$.} Since $x_nu<_{lex}x_1v$, we have that $\depth(S/I)\neq0$ by Proposition \ref{depthzero}. Hence $\projdim(S/I)<n$.

\textit{Case 4: $t>1$.} We obviously have $x_nu<_{lex}x_1v$ and, by Proposition \ref{depthzero}, $\depth(S/I)\neq0$. Therefore $\projdim(S/I)<n$.

We may conclude that $\projdim(S/I)<n$ in all the cases.

Let $w\in\Mon_d(S)$ be a monomial such that $|\mathcal L(u,v)|=|\mathcal L^i(w)|$. Since $I$ is a Gotzmann ideal, $I^{lex}$ is generated in degree $d$, that is, $I^{lex}=(\mathcal L^i(w))$. By Theorem \ref{betti I and I^lex}, $I$ and $I^{lex}$ have the same Betti numbers. In particular, we have
	\[\projdim(I)=\projdim(I^{lex}).
\]
Since $\projdim(S/I)<n$, we have $\projdim(S/I^{lex})<n$. The ideal $I^{lex}$ is stable in the sense of Eliahou and Kervaire, thus there exists $j<n$ such that $w=x_1^{d-1}x_j$. Therefore, $|\mathcal L(u,v)|=j<n$. By the hypothesis, $I$ is a Gotzmann ideal and $I$ is componentwise linear since it has a linear resolution. By Theorem \ref{minTay}, the Taylor resolution of $I$ is minimal. The conclusion follows by Theorem \ref{ot} and taking into account that $I$ is a lexsegment ideal.
\end{proof}

\addcontentsline{toc}{chapter}{Ideas for future}
\chapter*{Ideas for future} 

In the last years, many authors were interested in describing as much as possible the properties of the (squarefree) lexsegment ideals. Their importance in the commutative combinatorial algebra is given by the major role played by them in the study of the Hilbert function. The initial lexsegment ideals had been used by F.S. Macaulay \cite{Mac} in 1928 in order to determine an upper bound for the Hilbert function of a cyclic graded module. His work represents the first step in studying extremal properties of initial lexsegment ideals. Moreover, A. Bigatti, H. Hulett, and K. Pardue proved that every lexsegment ideal attains maximal Betti numbers among all graded ideals with the same Hilbert function. When computing Betti numbers, initial lexsegment ideals posses the property of having maximal Betti numbers \cite{Bi}, \cite{Hu} and \cite{Pa}. T. Deery \cite{De} proved that certain final lexsegment ideals have minimal Betti numbers for given Hilbert functions.

For arbitrary lexsegment ideals, many homological invariants are known. Passing to the squarefree case, there are many unsolved problems.

Firstly, in \cite{OO}, we succeeded to determine the minimal primary decomposition for completely squarefree lexsegment ideals. It naturally arise the same problem when considering squarefree lexsegment ideals which are not complete. This would allow us to determine the Krull dimension and the multiplicity. As an application of the minimal primary decomposition, it would be interesting in describing the sequentially Cohen--Macaulay squarefree lexsegment ideals which are not complete. This will offer a complete characterization of squarefree lexsegment ideals which are sequentially Cohen--Macaulay, as was already done for the non--squarefree case.

In \cite{EKOT} there are some conjectures and open problems concerning the depth and the regularity of an arbitrary squarefree lexsegment ideal. The depth is already computed in \cite{EOT} and \cite{EKOT}, when the ideal is generated in small degrees, namely in degrees $2$ and $3$. The proofs use different techniques and, even for $3-$degree case, the computation becomes difficult. For lexsegment edge ideals, that is squarefree lexsegment ideals generated in degree $2$, the Castelnuovo--Mumford regularity is given in \cite{EOT}. It will be interesting to derive formulae for the depth and Castelnuovo--Mumford regularity of an arbitrary squarefree lexsegment ideal. Moreover, one may consider the problem of computing the arithmetical rank of arbitrary squarefree lexsegment ideals $I$ and to characterize those for which $\ara(I)=\projdim_S(S/I)$.

For lexsegment ideals it is known that the properties of having a linear resolution and linear quotients coincides, \cite{BEOS}. In analogy with the non--squarefree case, it will be useful to determine a similar result for squarefree lexsegment ideals. Partial results are given in \cite{BEOS} only for squarefree lexsegment ideals which are not complete.

One may also consider the following problem proposed by J. Herzog: Describe the minimal primary decomposition of an arbitrary lexsegment ideal. Some results in this direction were obtained by M. Ishaq \cite{I}, who computed the set of associated prime ideals of a lexsegment ideal. For an associated prime ideal $\pp$ of a lexsegment ideal, it is useful to know the $\pp-$primary ideals.

\backmatter
\fancyhead{}

\end{document}